# Data-driven Multistage Distributionally Robust Linear Optimization with Nested Distance

Rui Gao, Rohit Arora and Yizhe Huang
Department of Information, Risk and Operations Management, University of Texas at Austin, Austin, TX 78705
rui.gao@mccombs.utexas.edu, arorarohit@utexas.edu, yizhe.huang@mccombs.utexas.edu

We study multistage distributionally robust linear optimization, where the uncertainty set is defined as a ball of distribution centered at a scenario tree using the nested distance. The resulting minimax problem is notoriously difficult to solve due to its inherent non-convexity. In this paper, we demonstrate that, under mild conditions, the robust risk evaluation of a given policy can be expressed in an equivalent recursive form. Furthermore, assuming stagewise independence, we derive equivalent dynamic programming reformulations to find an optimal robust policy that is time-consistent and well-defined on unseen sample paths. Our reformulations reconcile two modeling frameworks: the multistage-static formulation (with nested distance) and the multistage-dynamic formulation (with one-period Wasserstein distance). Moreover, we identify tractable cases when the value functions can be computed efficiently using convex optimization techniques.

## 1. Introduction

Distributionally Robust Optimization (DRO) is an emerging paradigm for data-driven decision-making, offering robust solutions that account for data uncertainty. For static problems, much significant progress has been made recently in terms of computation (Delage and Ye 2010, Goh and Sim 2010, Ben-Tal et al. 2013, Wiesemann et al. 2014, Esfahani and Kuhn 2018), regularization (Lam 2016, Duchi and Namkoong 2019, Shafieezadeh-Abadeh et al. 2019, Gao et al. 2024), and statistical guarantees (Delage and Ye 2010, Lam 2019, Duchi et al. 2021, Blanchet et al. 2019, Gao 2023), etc. However, extending the paradigm to the multistage setting continues to be challenging, and the results obtained so far have been limited and not entirely satisfactory.

In multistage problems, scenario trees are commonly used to represent data processes. There are two main data-driven approaches for constructing scenario trees. The first approach is based on Monte Carlo sampling techniques (Shapiro 2003), such as conditional sampling that includes stagewise independent sampling as an important case. The second approach is based on scenario generation and reduction (Dupačová et al. 2000, Høyland and Wallace 2001, Dupačová et al. 2003, Henrion and Römisch 2022), aiming at generating a reasonable number of scenarios that represent the essential characteristics of the uncertainty. Scenario trees serve as discrete approximations of the underlying true stochastic process and are utilized to formulate decision-making problems along the sample paths within the tree. These problems are commonly known as scenario approximation or sample average approximation (SAA). However, it is crucial to acknowledge that policies solved using SAA are not defined for unseen sample paths. In most cases, heuristic policies developed for unseen sample paths lack optimality guarantees (Ben-Tal et al. 2009, Note and Remarks 14.1). This limitation highlights the challenge of data scarcity in multistage problems (Shapiro and Nemirovski 2005). Clearly, there is a need for a distributionally robust model that can handle unseen sample paths in a principled manner.

### 1.1. Challenges in Choosing a Distributionally Robust Model

The formulation of a multistage distributionally robust model can be approached in various ways, however, there is currently no consensus in the literature on the best approach. This lack of consensus stems from the fact that different extensions of the single-stage formulation can result in different frameworks, each with its own interpretation and perceived naturalness. As a result, continued debate and exploration in the field persist, as highlighted in recent studies by Pichler and Shapiro (2021), Shapiro and Pichler (2023).

One straightforward formulation, referred to as the *multistage-static formulation*, is a direct extension of the standard single-stage formulation

$$\inf_{\boldsymbol{x}\in\mathcal{X}} \sup_{\mathbb{P}\in\mathfrak{M}} \mathbb{E}_{\boldsymbol{\xi}_{[T]}\sim\mathbb{P}}\left[\sum_{t=1}^{T} f_t(\boldsymbol{x}_t,\boldsymbol{\xi}_t)\right]. \tag{1}$$

Here $f_t(\boldsymbol{x}_t,\boldsymbol{\xi}_t)$ denotes the per-stage cost associated with a $T$-stage sample path $\boldsymbol{\xi}_{[T]} := (\xi_1,\boldsymbol{\xi}_2,\ldots,\boldsymbol{\xi}_T)\in\Xi_1\times\cdots\times\Xi_T$ under a policy $\boldsymbol{x}=(x_1,\boldsymbol{x}_2,\ldots,\boldsymbol{x}_T)$[1]. The formulation hedges against an uncertainty set $\mathfrak{M}$ of $T$-stage stochastic processes and optimizes the worst-case expected cumulative cost over the set $\mathcal{X}$ of policies satisfying some feasibility and non-anticipativity constraints $x_1\in\mathcal{X}_1$, $\boldsymbol{x}_t\in\mathcal{X}_t(\boldsymbol{x}_{t-1})$, $t=2,\ldots,T$. Typically, the uncertainty set is constructed based on summary statistics such as support and moment information (Bertsimas et al. 2019, Xin and Goldberg 2022), or based on statistical distance such as relative entropy (Hansen and Sargent 2001), Wasserstein distance (Bertsimas et al. 2023, Sturt 2023), and nested distance (Analui and Pflug 2014, Glanzer et al. 2019). The multistage-static formulation (1) is conceptually simple and offers a clear interpretation as a method to mitigate uncertainty in the data process. However, the multistage-static objective function is not explicitly adjusted for dynamics of the decision process (Pichler and Shapiro 2021), making it challenging to utilize dynamic programming recursions. This lack of adjustment may also raise concerns about time inconsistency (Iancu et al. 2015), which is criticized in decision theory for its violation of rational behavior.

An alternative formulation, referred to as the *multistage-dynamic formulation*, has been devised to facilitate dynamic programming recursion and is routinely employed in computational studies. In this formulation, the cost-to-go function takes on a recursive form

$$Q_t(\boldsymbol{x}_{t-1},\boldsymbol{\xi}_{[t]}) = \inf_{\boldsymbol{x}_t\in\mathcal{X}_t(\boldsymbol{x}_{t-1})}\left\{f_t(\boldsymbol{x}_t,\boldsymbol{\xi}_t) + \sup_{\boldsymbol{\xi}_{t+1}\sim\mathbb{P}_{t+1}\in\mathfrak{M}_{t+1}} \mathbb{E}_{\mathbb{P}_{t+1}}\left[Q_{t+1}(\boldsymbol{x}_t,\boldsymbol{\xi}_{[t+1]})\right]\right\},\quad t\in[T], \tag{2}$$

and $Q_{T+1}(\cdot,\cdot)\equiv 0$. Here the uncertainty set[2] $\mathfrak{M}_{t+1}$ can be defined through composite distributionally robust functionals and conditional distributionally robust functionals (Pichler and Shapiro 2021, Shapiro and Pichler 2023) (see also conditional risk mappings (Ruszczyński and Shapiro 2006)). Of particular interest is the stagewise independent setting, where these two functionals are equivalent. In comparison to the multistage-static formulation (1), the multistage-dynamic formulation (2) is generally more computationally friendly. Nevertheless, it should be noted that if not appropriately specified, it can be overly conservative and lack interpretability. For instance, the composition of single-period Average Value at Risk (AVaR) takes conditional tail expectations of conditional tail expectations. The resulting multi-period risk measure does not offer the same straightforward interpretation as the single-period AVaR (Shapiro 2012) and can potentially lead to overly conservative risk assessments (Iancu et al. 2015). Another example is seen in the formulation with composite distributionally robust functionals, where different sample paths at different stages can be associated with different worst-case $T$-stage distributions. This raises concerns about the pessimism of the resulting policy.

If we can find a choice of $\mathfrak{M}$ such that the multistage-static formulation (1) can be equivalently represented in a multistage-dynamic form (2), it opens up the possibility of solving (1) through dynamic programming, yielding a time-consistent optimal robust policy. Conversely, if the multistage-dynamic formulation (2) can be transformed into a multistage-static form (1) with an interpretable uncertainty set $\mathfrak{M}$, it could help mitigate the issue of conservativeness by allowing us to work with a more interpretable multi-period risk measure induced from the multistage-static counterpart. By reconciling these two formulations, we can leverage the advantages of both approaches and obtain a comprehensive understanding of the problem. Unfortunately, for generic multistage linear problems, the only identified instances of $\mathfrak{M}$ in the existing literature are limited to two extreme cases (Shapiro et al. 2021, Remark 33): singleton (non-robust) and entirety (extremely robust).

[1] We use bold font for random variables and regular font for deterministic values like constants or elements in the sample space. In line with the convention in stochastic programming literature, we consider the first-stage data $\boldsymbol{\xi}_1$ as deterministic. Therefore, in the sequel, we will interchangeably use $x_1$ and $\boldsymbol{x}_1$, and $\xi_1$ and $\boldsymbol{\xi}_1$. We simply set $\Xi_1=\{\xi_1\}$.

[2] Note that the uncertainty set can be dependent on the history or some reference (nominal) data process. Here we omit such a dependence for notational simplicity.

In light of the discussion above, an important and open question arises:

*Question 1* : *Is there a modeling choice that can be easily interpreted from the viewpoint of the multistage-static formulation* (1)*, while simultaneously allowing for an equivalent decomposition into the multistage-dynamic formulation* (2) *with interpretable one-period uncertainty sets?*

In a nutshell, we will affirmatively address this question by showing that the multistage-static formulation (1) with *nested distance* is indeed equivalent to the multistage-dynamic formulation (2) with one-period Wasserstein distance.

### 1.2. Challenges in Solving Multistage DRO with Nested Distance

Given our answer to Question 1, the distributionally robust problem (1) with nested distance can be solved using dynamic programming (2). However, before delving into our contributions, let us first examine the computational challenges associated with problem (1) with nested distance. This discussion also provides insights into the challenges involved in establishing the equivalence results.

Solving problem (1) with nested distance is believed to be complex (Pflug and Pichler 2014, Page 239). Indeed, existing algorithmic approaches Analui and Pflug (2014), Glanzer et al. (2019) involve numerically searching over probability distributions in the uncertainty set and are computationally expensive. To reduce the computational burden, these approaches fix scenario values and tree structures, allowing changes only in the probability weights.

In essence, the computational challenges arise from both the inner robust risk evaluation and the outer policy optimization. The inner robust risk evaluation involves optimizing an infinite-dimensional space of probability distributions. In distributionally robust optimization literature, to make this problem more tractable, a common approach is to reformulate it into a finite-dimensional one using duality and conditioning techniques. However, the nested distance uncertainty set becomes non-convex as soon as $T \geq 3$ (Pflug and Pichler 2014, Fig. 1.17). This non-convexity poses a challenge when using duality-based arguments. Meanwhile, the outer policy optimization involves an infinite-dimensional optimization over policies. In contrast to its non-robust scenario approximation counterpart, which only specifies policy values on a finite number of sample paths, solving (1) requires specifying policy values for all possible sample paths of distributions in the uncertainty set. Furthermore, even finding the first-stage optimal decision in the case of $T = 2$ is shown to be NP-hard in general (Hanasusanto and Kuhn 2018, Xie 2020).

On the positive side, there have been some notable advancements in solving two-stage Wasserstein DRO and multistage stochastic programming. For instance, for two-stage Wasserstein DRO, under certain conditions, Hanasusanto and Kuhn (2018) develop co-positive program reformulations using 2-Wasserstein distance and linear program reformulations using 1-Wasserstein distance, and Xie (2020) provides tractable convex program reformulations using $\infty$-Wasserstein distance. Moreover, for linear and convex multistage (non-robust) problems, significant progress has been made recently in Stochastic Dual Dynamic Programming (SDDP) (Birge 1985, Pereira and Pinto 1991); we refer to Lan and Shapiro (2024) for an excellent review. Considering that our problem (1) with nested distance reduces to a two-stage Wasserstein DRO problem when $T = 2$, and its dynamic formulation (2) reduces to its multistage stochastic programming when the radius is zero, we aim to leverage these advancements to enhance the computational tractability of our problem (1) or (2). We pose the following question:

*Question 2* : *Can we identify conditions under which there is a computationally tractable way to find the optimal solution for the formulation* (1) *with nested distance or its equivalent reformulation* (2)*?*

To address this question, we will identify conditions under which the optimal policy can be solved using SDDP and convex optimization. These conditions extend from those for two-stage Wasserstein DRO (Esfahani and Kuhn 2018, Hanasusanto and Kuhn 2018, Xie 2020).

### 1.3. Our Contributions

Our contributions are as follows.

(I) We derive dynamic programming reformulations for problem (1) with nested distance, providing an affirmative answer to Question 1. In Section 3.1, we focus on the inner maximization of (1), which evaluates the worst-case risk for a fixed policy. We show that it can be equivalently decomposed into dynamic programs defined via one-period Wasserstein distance relative to the nominal conditional distribution. In Section 3.2, we study the outer minimization of (1) under the assumption of stagewise independence. We show that the multistage-static formulation (1) with nested distance is equivalent to the multistage-dynamic formulation (2) with one-period Wasserstein distance. To the best of our knowledge, this is the first non-degenerate uncertainty set that reconciles both static and dynamic formulations for generic multistage linear programs.

(II) As a byproduct of our proof for the equivalence results, we derive a time-consistent optimal robust policy that solves both (1) and (2). Notably, our formulations render a policy that is well-defined for every possible sample path of the distributions in the uncertainty set with provable optimality guarantees.

Furthermore, in Section 4, we identify conditions under which the value function can be computed efficiently via convex optimization, providing a positive answer to Question 2. In particular, in Section 4.1 we show that, when the uncertainty appears on the objective only, the robust value function is equivalent to a norm-regularized SAA value function that penalizes large norms on the decision variables. In Section 4.2 we show that, when the uncertainty appears on the right-hand side only, the equivalent reformulation of the value function de-regularizes the SAA value function by encouraging large norms on the dual variables. This leads to a solution that hedges extreme perturbations of the right-hand side when $p = \infty$, and a solution that coincides with the non-robust counterpart when $p = 1$.

### 1.4. Related Literature

The study of multistage DRO with nested distance was first introduced by Analui and Pflug (2014) (see also Pflug and Pichler (2014, Section 7.3)) and further explored in the follow-up work by Glanzer et al. (2019). In comparison, our results differ from theirs in several major ways. First, we aim to develop dynamic programming and convex tractable reformulations. In contrast, their successive programming approach does not consider dynamic programming and has to deal with the non-convexity of subproblems. Second, unlike their algorithms that essentially restrict the support of sample paths, our approach does not impose extra restrictions on the stochastic processes while still maintaining computational efficiency.

When the nominal stochastic process follows a fan-shaped distribution, namely, the conditional distribution of $\widehat{\boldsymbol{\xi}}_{t+1}$ given $\widehat{\boldsymbol{\xi}}_t$ is a Dirac measure for all $t \geq 2$, the multistage-static formulation (1) with nested distance coincides with the formulation using Wasserstein distance (Pflug and Pichler 2012, Remark 6). Consequently, in the case of $p = \infty$, our problem (1) reduces to the multistage-static formulation with $\infty$-Wasserstein distance, which has been studied in Bertsimas et al. (2023). However, their focus is on developing asymptotic optimality guarantees and solving it with affine policies and finite adaptability, without assuming stagewise independence. Different from their work, our goal is to find the optimal policy without restricting it to specific forms, while assuming stagewise independence.

As stagewise independence is a common assumption (Lan and Shapiro (2024), see also Shapiro et al. (2021, Remark 5)), we believe our results are relevant for computational studies.

Below, we provide a review of some relevant literature that is not covered in the previously mentioned works.

*On nested distance.* In stochastic programming literature, the nested distance was first introduced and studied in the seminal works Pflug (2010), Pflug and Pichler (2012). Since then, it has been applied to scenario generation/reduction, stability analysis and approximation of multistage stochastic programming (Pflug and Pichler 2015, Maggioni and Pflug 2016, Kovacevic and Pichler 2015, Chen and Yan 2018, Horejšová et al. 2020). Its statistical properties (Pflug and Pichler 2016, Glanzer et al. 2019, Veraguas et al. 2020) and computational properties (Cabral and da Costa 2017, Pichler and Weinhardt 2022) have also been investigated. In the optimal transport and mathematical finance literature, the nested distance is also called bi-causal transport distance or adapted Wasserstein distance (Backhoff-Veraguas et al. 2020, Backhoff et al. 2022). Incorporating nested distance in multistage distributionally robust optimization was first considered in Analui and Pflug (2014), Pflug and Pichler (2014) and then in Glanzer et al. (2019) for a pricing problem. Unlike the algorithmic approach in these works, our algorithm is more computationally friendly and enjoys similar tractability as its SAA counterpart. A recent paper Yang et al. (2022) studies DRO with side information using causal transport distance, which can be viewed as a convex relaxation of our problem in three stages. The idea of imposing non-anticipativity constraints on the transport plan can be traced back to the Yamada-Watanabe criterion for stochastic differential equations (Yamada and Watanabe 1971) as well as the causal transportation in continuous time (Lassalle 2018) and in discrete time (Backhoff et al. 2017).

*On other multistage distributionally robust models.* In the introduction, we have provided a list of some common uncertainty sets for the multistage-static formulation (1). Now let us briefly review some choices of uncertainty sets for the multistage-dynamic formulation (2). These include the AVaR (Ruszczyński and Shapiro 2006), entropic risk measure (Dowson et al. 2020), $\phi$-divergence (Klabjan et al. 2013, Hanasusanto and Kuhn 2013, Park and Bayraksan 2023, Rahimian et al. 2022), and Wasserstein distance (Shapiro and Pichler 2023). When assuming stagewise independence, both moment-based sets (Shapiro and Xin 2020, Xin and Goldberg 2021, Yu and Shen 2020) and statistical distance-based sets are widely considered. Examples of statistical distances include $\chi^2$-divergence (Philpott et al. 2018), $L_\infty$-norm (Huang et al. 2017), and Wasserstein distance (Duque and Morton 2020, Zhang and Sun 2020, 2022).

In addition to enabling a dynamic programming reformulation and rendering a time-consistent policy that is defined for every unseen sample path, our choice of the nested distance offers several other appealing features from a modeling standpoint. (i) It fully utilizes the entire distributional information, distinguishing it from moment-based and risk-measure-based sets that rely only on partial data information like moments and tail risks. (ii) It provides effective protection against data perturbations that extend beyond the support of the nominal scenario tree. This distinguishes it from divergence-based sets, which impose strict restrictions on the support of relevant distributions (Bayraksan and Love 2015). (iii) It defines a plausible family of stochastic processes, comprised of non-anticipative perturbations of sample paths. This stands in contrast to the Wasserstein distance, which permits perturbations dependent on future information. Further elaboration on this topic can be found in Section 2.2.

*On the computation of multistage DRO with transport distance-based sets.* There are several computational works that directly solve the multistage-dynamic formulation (2) without considering the multistage-static formulation (1). For the formulation (2) with stagewise independent 1-Wasserstein uncertainty sets, Duque and Morton (2020) propose an SDDP algorithm while restricting the support of scenarios on a pre-specified finite set. Zhang and Sun (2020) propose dual dynamic programming algorithms and analyze their complexity while assuming finitely-supported distributions. This issue is resolved in a recent work by Zhang and Sun (2022). However, these works do not provide a policy or the output policy does not come with provable optimality guarantees. In contrast, our formulation yields a policy with robust optimality guarantees for generic multistage robust linear problems, which also justifies the use of some heuristic policies that extend from sample paths within the scenario tree to those outside the tree (Thénié and Vial 2008, Shapiro et al. 2012, Defourny et al. 2013, Keutchayan et al. 2017). For two-stage robust and distributionally robust optimization, the optimality of affine decision rules has been investigated in Bertsimas and Goyal (2012), El Housni and Goyal (2021), Georghiou et al. (2021). Optimal robust policies have also been studied in other multistage DRO problems, such as robust control (Bertsimas et al. 2010, Iancu et al. 2013, Taşkesen et al. 2024), robust optimal stopping (Sturt 2023), robust contextual optimization Zhang et al. (2023). These techniques, however, do not extend to our settings straightforwardly. Finally, it is worth mentioning that in this paper, we are content with a reformulation that is as tractable as multistage SAA, although the latter alone may have many computational issues (Shapiro et al. 2021, Chapter 5.8) that are beyond the scope of this paper.

*On the concept of time consistency.* Various concepts of time consistency have been discussed in economics literature (Strotz et al. 1955, Hansen and Sargent 2001, Epstein and Schneider 2003, Etner et al. 2012), in mathematical finance literature (Wang 1999, Föllmer and Schied 2002, Artzner et al. 2007, Roorda and Schumacher 2007, Cheridito and Kupper 2009) and in robust control literature (Iyengar 2005, Nilim and El Ghaoui 2005, Wiesemann et al. 2013). Our notion of time consistency is more aligned with the stochastic programming literature (Ruszczyński and Shapiro 2006, Shapiro 2012, 2016, Shapiro and Xin 2020, Xin and Goldberg 2021, Pichler et al. 2022). In general, the multistage-static formulation (1) does not have a time-consistent optimal policy due to a lack of dynamic programming representation (Pichler and Shapiro 2021). For certain classes of problems in inventory control, it has been shown that the multistage-static formulation (1) with some moment-based uncertainty sets has a time-consistent optimal policy under certain conditions (Shapiro 2012, Xin and Goldberg 2021, 2022). In contrast, our result on time consistency holds for generic multistage linear programs.

The rest of the paper is organized as follows. In Section 2, we introduce the distributionally robust linear multistage program and provide a quick overview of the nested distance used for constructing the ambiguity set. In Section 3, we develop dynamic programming reformulations for both risk evaluation and policy optimization. In Section 4, we specialize our result into cases that admit tractable convex reformulations. In Section 5, we apply our result to the portfolio selection problem using an SDDP algorithm. We conclude the paper in Section 6. All proofs are deferred to the Appendices.

## 2. Multistage Distributionally Robust Optimization with Nested Distance

In this section, we present the distributionally robust formulation and discuss the nested distance via examples.

### 2.1. Multistage Stochastic Programming and its Distributionally Robust Counterpart

Consider a $T$-stage stochastic linear optimization problem

$$\begin{aligned} \inf_{x_1, \boldsymbol{x}_2, \ldots, \boldsymbol{x}_T} \quad & \mathbb{E}_{\mathbb{P}}\left[c_1^\top x_1 + \boldsymbol{c}_2^\top \boldsymbol{x}_2 + \cdots + \boldsymbol{c}_T^\top \boldsymbol{x}_T\right], \\ \text{s.t.} \quad & A_1 x_1 = b_1,\ x_1 \geq 0, \\ & \boldsymbol{B}_t \boldsymbol{x}_{t-1} + \boldsymbol{A}_t \boldsymbol{x}_t = \boldsymbol{b}_t,\ \boldsymbol{x}_t \geq 0, \quad t = 2, \ldots, T, \end{aligned}$$

where $\boldsymbol{\xi}_t := (\boldsymbol{A}_t, \boldsymbol{B}_t, \boldsymbol{c}_t, \boldsymbol{b}_t) \in \Xi_t \subset \mathbb{R}^{d_t}$, $t \in [T] := \{1, \ldots, T\}$, are data vectors and matrices, some or all of which may be random. For $t \in [T]$, we denote by $\boldsymbol{\xi}_{[t]} := (\boldsymbol{\xi}_1, \ldots, \boldsymbol{\xi}_t) \in \Xi_{[t]} := \Xi_1 \times \cdots \Xi_t$ the history of the data process up to time $t$. Let $\mathcal{P}(\Xi_{[t]})$ be the set of probability distributions on $\Xi_{[t]}$, $t \in [T]$. For any distribution $\mathbb{P} \in \mathcal{P}(\Xi_{[T]})$ of a stochastic process $\boldsymbol{\xi}_{[T]}$, we denote by $\mathbb{P}_{[t]}$ the marginal distribution of $\boldsymbol{\xi}_{[t]}$ under $\mathbb{P}$, and by $\mathbb{P}_{t|\boldsymbol{\xi}_{[t-1]}}$ the conditional distribution of $\boldsymbol{\xi}_t$ given the history $\boldsymbol{\xi}_{[t-1]}$.

The minimization is performed over the set of non-anticipative policies $\boldsymbol{x}_t = \boldsymbol{x}_t(\boldsymbol{\xi}_{[t]})$, each of which is measurable with respect to $\sigma(\boldsymbol{\xi}_{[t]})$, the $\sigma$-algebra induced by $\boldsymbol{\xi}_{[t]}$. To ease notations, we denote the feasible regions

$$\mathcal{X}_1 := \{x_1 \geq 0 : \ A_1 x_1 = b_1\},$$
$$\mathcal{X}_t(x_{t-1}, \boldsymbol{\xi}_t) := \{x_t \geq 0 : \ \boldsymbol{A}_t x_t = \boldsymbol{b}_t - \boldsymbol{B}_t x_{t-1}\}, \quad t = 2, \ldots, T.$$

We will always assume $\mathcal{X}_1 \neq \emptyset$. The multistage problem above admits a dynamic programming formulation

$$Q_t(x_{t-1}, \boldsymbol{\xi}_{[t]}) = \inf_{x_t \in \mathcal{X}_t(x_{t-1}, \boldsymbol{\xi}_t)} \left\{ c_t^\top x_t + \mathbb{E}_{\mathbb{P}_{t+1|\boldsymbol{\xi}_{[t]}}} [Q_{t+1}(x_t, \boldsymbol{\xi}_{[t+1]})] \right\}, \quad t \in [T],$$

with $Q_{T+1} \equiv 0$ and $\mathcal{X}_1(x_0, \xi_1) \equiv \mathcal{X}_1$. We represent the set of non-anticipative and feasible policies as

$$\mathcal{X} := \left\{ \left(x_1, \boldsymbol{x}_2(\cdot), \ldots, \boldsymbol{x}_T(\cdot)\right) : \ \boldsymbol{x}_t \in \mathcal{X}_t(\boldsymbol{x}_{t-1}), \ t \in [T] \right\}.$$

Here the shorthand notation $\boldsymbol{x}_t \in \mathcal{X}_t(\boldsymbol{x}_{t-1})$ is interpreted as $\boldsymbol{x}_t(\xi_{[t]}) \in \mathcal{X}_t(\boldsymbol{x}_{t-1}(\xi_{[t-1]}), \xi_t)$ for all $\xi_{[t]} \in \Xi_{[t]}$.

Quite often in practice, the data-generating distribution of the random process $\boldsymbol{\xi}_{[T]}$ is not known exactly. A common approach is to replace the underlying data-generating distribution with a scenario tree $\widehat{\mathbb{P}}$, which is typically constructed using conditional sampling or scenario reduction. Let us denote by $\widehat{\boldsymbol{\xi}}_{[T]}$ the stochastic process with a finitely-supported distribution $\widehat{\mathbb{P}}$ and by $\widehat{\Xi}_t$ the support of $\widehat{\boldsymbol{\xi}}_t$, $t \in [T]$. To account for the distributional uncertainty, we consider the following multistage distributionally robust optimization

$$\inf_{\boldsymbol{x} \in \mathcal{X}} \sup_{\mathbb{P} \in \mathfrak{M}} \mathbb{E}_{\mathbb{P}} \left[ \sum_{t \in [T]} \boldsymbol{c}_t^\top \boldsymbol{x}_t(\boldsymbol{\xi}_{[t]}) \right], \tag{$\mathsf{P}_{\text{static}}$}$$

and the distributional uncertainty set $\mathfrak{M}$ specifies a set of $T$-stage distributions to hedge against. In particular, we consider the following uncertainty set

$$\mathfrak{M} := \left\{ \mathbb{P} \in \mathcal{P}(\Xi_{[T]}) : \ \mathsf{D}_p(\widehat{\mathbb{P}}, \mathbb{P}) \leq \vartheta \right\}, \tag{3}$$

where $\vartheta > 0$ is the radius of the uncertainty set, and $\mathsf{D}_p$ is the *p-nested distance* proposed by Pflug (2010), Pflug and Pichler (2012). As will be elaborated on in Section 2.2, the choice of the nested distance takes account of the information evolution in the multistage problem. We will also consider a soft robust formulation when $p \in [1, \infty)$

$$\inf_{\boldsymbol{x} \in \mathcal{X}} \sup_{\mathbb{P} \in \mathcal{P}(\Xi_{[T]})} \left\{ \mathbb{E}_{\mathbb{P}} \left[ \sum_{t \in [T]} \boldsymbol{c}_t^\top \boldsymbol{x}_t(\boldsymbol{\xi}_{[t]}) \right] - \lambda \mathsf{D}_p^p(\widehat{\mathbb{P}}, \mathbb{P}) \right\}, \tag{$\mathsf{P}_{\text{static}}$-soft}$$

for some fixed $\lambda \in (0, \infty)$.

### 2.2. Nested Distance

Similar to the Wasserstein distance, the nested distance is based on an optimal transport problem. But in addition to the marginal constraints on the transport plan as in the definition of Wasserstein distance, it also requires that the transport plan should be non-anticipative with respect to the filtration $\sigma(\widehat{\boldsymbol{\xi}}_{[t]}) \otimes \sigma(\boldsymbol{\xi}_{[t]})$. In other words, the nested distance calculates the minimum cost needed to transport probability mass from $\widehat{\mathbb{P}}$ to $\mathbb{P}$ among a set of non-anticipative transport plans.

Let $\mathsf{d}(\cdot, \cdot)$ be a metric on $\Xi_{[T]}$. For any $\widehat{\mathbb{P}}, \mathbb{P} \in \mathcal{P}(\Xi_{[T]})$, we denote by $\Gamma(\widehat{\mathbb{P}}, \mathbb{P})$ the set of joint distributions on $\Xi_{[T]}^2$ with marginals $\widehat{\mathbb{P}}$ and $\mathbb{P}$. For a joint distribution $\gamma \in \Gamma(\widehat{\mathbb{P}}, \mathbb{P})$, we use $\gamma_{\widehat{\boldsymbol{\xi}}_t | (\widehat{\boldsymbol{\xi}}_{[t-1]}, \boldsymbol{\xi}_{[t-1]})}$ to denote the conditional distribution of $\widehat{\boldsymbol{\xi}}_t$ given $(\widehat{\boldsymbol{\xi}}_{[t-1]}, \boldsymbol{\xi}_{[t-1]})$ under $\gamma$. Recall that $\widehat{\mathbb{P}}_{t|\widehat{\boldsymbol{\xi}}_{[t-1]}}$ denotes the conditional distribution of $\widehat{\boldsymbol{\xi}}_t$ given $\widehat{\boldsymbol{\xi}}_{[t-1]}$ under $\widehat{\mathbb{P}}$.

Definition 1 (Nested Distance). Define the set of non-anticipative transport plans

$$\Gamma_{bc}(\widehat{\mathbb{P}}, \mathbb{P}) = \left\{ \gamma \in \Gamma(\widehat{\mathbb{P}}, \mathbb{P}) : \ \gamma_{\widehat{\xi}_{t+1} | (\widehat{\xi}_{[t]}, \xi_{[t]})} = \widehat{\mathbb{P}}_{t+1|\widehat{\xi}_{[t]}}, \ \gamma_{\xi_{t+1} | (\widehat{\xi}_{[t]}, \xi_{[t]})} = \mathbb{P}_{t+1|\xi_{[t]}}, \ \forall t = 1, \dots, T-1 \right\}. \quad (4)$$

The nested distance $\mathsf{D}_p(\widehat{\mathbb{P}}, \mathbb{P})$ between $\widehat{\mathbb{P}}$ and $\mathbb{P}$ is defined as

$$\mathsf{D}_p(\widehat{\mathbb{P}}, \mathbb{P}) := \begin{cases} \left( \inf_{\gamma \in \Gamma_{bc}(\widehat{\mathbb{P}}, \mathbb{P})} \mathbb{E}_{(\widehat{\xi}_{[T]}, \xi_{[T]}) \sim \gamma} \left[ \mathsf{d}(\widehat{\xi}_{[T]}, \xi_{[T]})^p \right] \right)^{1/p}, & p \in [1, \infty), \\ \inf_{\gamma \in \Gamma_{bc}(\widehat{\mathbb{P}}, \mathbb{P})} \gamma\text{-ess sup } \mathsf{d}(\widehat{\xi}_{[T]}, \xi_{[T]}), & p = \infty. \end{cases} \quad (5)$$

◊

The non-anticipativity constraints (4) can be equivalently stated as follows. Under the joint distribution $\gamma$,

$$\xi_{[t]} \perp \widehat{\xi}_{t+1} \,|\, \widehat{\xi}_{[t]}, \quad (6a)$$

$$\widehat{\xi}_{[t]} \perp \xi_{t+1} \,|\, \xi_{[t]}, \quad (6b)$$

namely, $\xi_{[t]}$ and $\widehat{\xi}_{t+1}$ are *conditionally independent* given $\widehat{\xi}_{[t]}$; and $\widehat{\xi}_{[t]}$ and $\xi_{t+1}$ are conditionally independent given $\xi_{[t]}$. Suppose there exists a transport map $\mathbb{T} = (\mathbb{T}_1, \dots, \mathbb{T}_T)$ from $\widehat{\mathbb{P}}$ to $\mathbb{P}$, then (6a) implies that where $\widehat{\xi}_{[t]}$ is transported (i.e., $\xi_{[t]} = \mathbb{T}_{[t]}(\widehat{\xi}_{[T]})$) is independent of the future $\widehat{\xi}_{t+1}$. Thereby, $\mathbb{T}$ satisfies (6a) if and only if it is of the form $\mathbb{T}(\widehat{\xi}_{[T]}) = (\mathbb{T}_1(\widehat{\xi}_1), \mathbb{T}_2(\widehat{\xi}_{[2]}), \dots, \mathbb{T}_t(\widehat{\xi}_{[t]}), \dots, \mathbb{T}_T(\widehat{\xi}_{[T]}))$, $\forall \widehat{\xi}_{[T]} \in \widehat{\Xi}$. Similarly, the condition (6b) indicates that where $\xi_{[t]}$ is transported should not be dependent on the future information $\xi_{t+1}$. Thereby, if $\mathbb{T}_t$, $t = 1, \dots, T$, are invertible, then $\mathbb{T}$ satisfies (6b) as well (Backhoff et al. 2017, Remark 3.4). The equivalent definition (6) of non-anticipative transport plans (4) provides a convenient way to check whether a transport plan is non-anticipative or not.

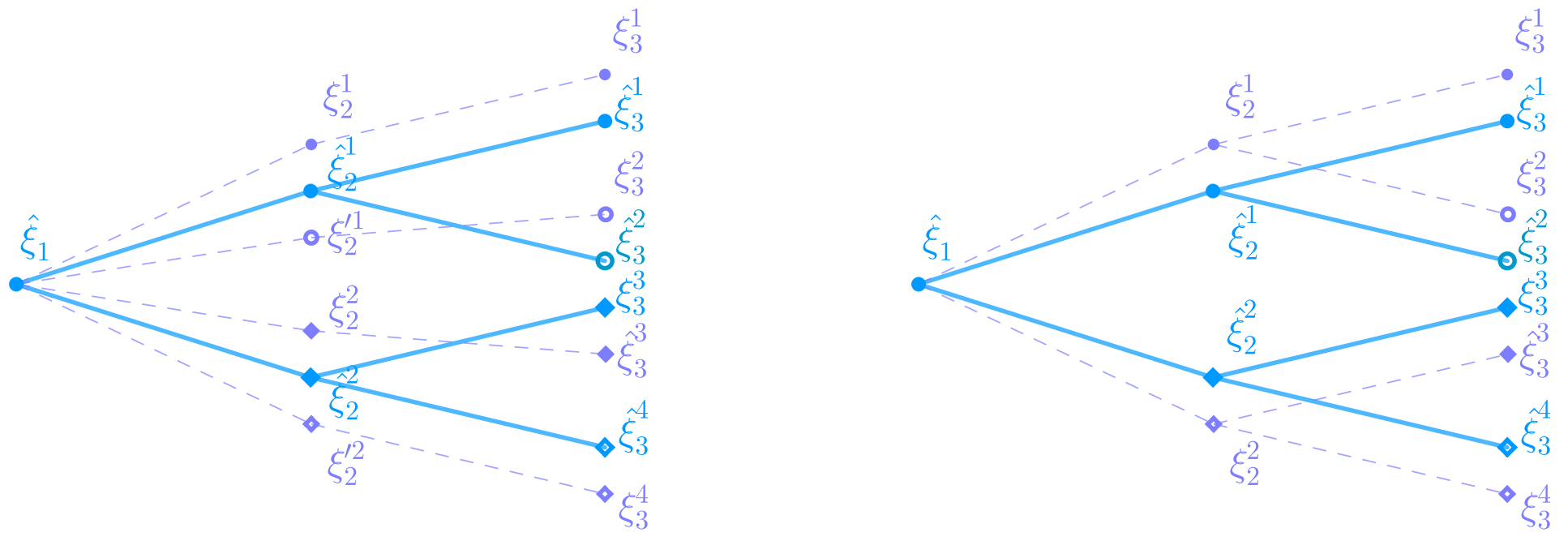

**Figure 1** Illustration of anticipative (left) and non-anticipative (right) transport plans.

In Figure 1, we illustrate the difference between an anticipative transport plan (left) and a non-anticipative transport plan (right). In both figures, the four solid sample paths are transported to the four dashed sample paths. In the left figure, the four solid sample paths from top to bottom are transported to the four dashed sample paths from top to bottom respectively. In particular, $(\widehat{\xi}_2^1, \widehat{\xi}_3^1)$ is transported to $(\xi_2^1, \xi_3^1)$ whereas $(\widehat{\xi}_2^1, \widehat{\xi}_3^2)$ is transported to $(\xi_2'^1, \xi_3^2)$. This means where $\widehat{\xi}_2^1$ is transported (namely, $\xi_2^1$ or $\xi_2'^1$) depends on the realization in the third stage (namely, $\widehat{\xi}_3^1$ or $\widehat{\xi}_3^2$). In contrast, in the right figure, $\widehat{\xi}_2^1$ is always transported to $\xi_2^1$ and $\xi_2^1$ is always transported to $\widehat{\xi}_2^1$, regardless of the value in the third stage.

In the literature, a transport plan from $\widehat{\mathbb{P}}$ to $\mathbb{P}$ is termed *causal* if it satisfies the non-anticipativity constraint (6a), and is termed *bi-causal* if it satisfies both non-anticipativity constraints (6). When only the first set of constraints in (4) is imposed, the resulting distance is called *causal transport distance* (Backhoff et al. 2017), denoted as $\mathsf{C}_p(\widehat{\mathbb{P}}, \mathbb{P})$. If we replace $\Gamma_{bc}$ by $\Gamma$, then (5) becomes the defining expression for the Wasserstein distance. Note that $\Gamma_{bc}(\widehat{\mathbb{P}}, \mathbb{P})$ is always a non-empty subset of $\Gamma(\widehat{\mathbb{P}}, \mathbb{P})$, containing at least the independent transport plan, namely the product distribution with marginals $\widehat{\mathbb{P}}$ and $\mathbb{P}$.

Let us illustrate these concepts with the following examples.

Example 1 (Non-anticipative data perturbations). In many applications, the data process $\boldsymbol{\xi}_{[T]}$ often adheres to a causal relationship, as represented by the following causal diagram

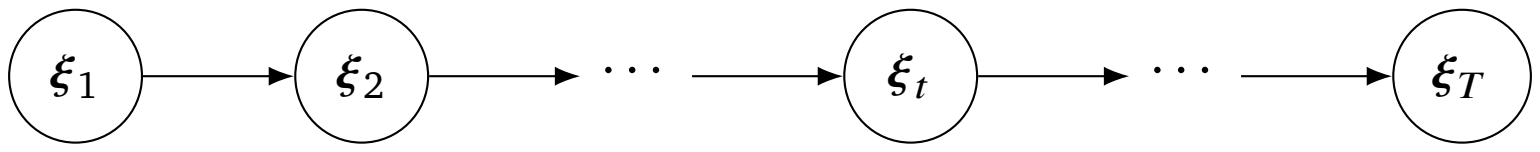

This relationship can be observed in various scenarios, such as the demand process for a product or the return rate process of financial assets. The data uncertainty of $\boldsymbol{\xi}_t$ can arise directly from errors like sampling or measurement inaccuracies in $\boldsymbol{\xi}_t$. Alternatively, it might be indirectly influenced by errors in statistical modeling or data processing propagated from historical data $\boldsymbol{\xi}_{[t-1]}$. Consequently, it is logical to consider data perturbations that exhibit historical dependencies. However, such perturbations should not be dependent on future uncertainties, as these are typically unknown. This rationale provides a justification for the non-anticipativity constraints in Definition 1. ♣

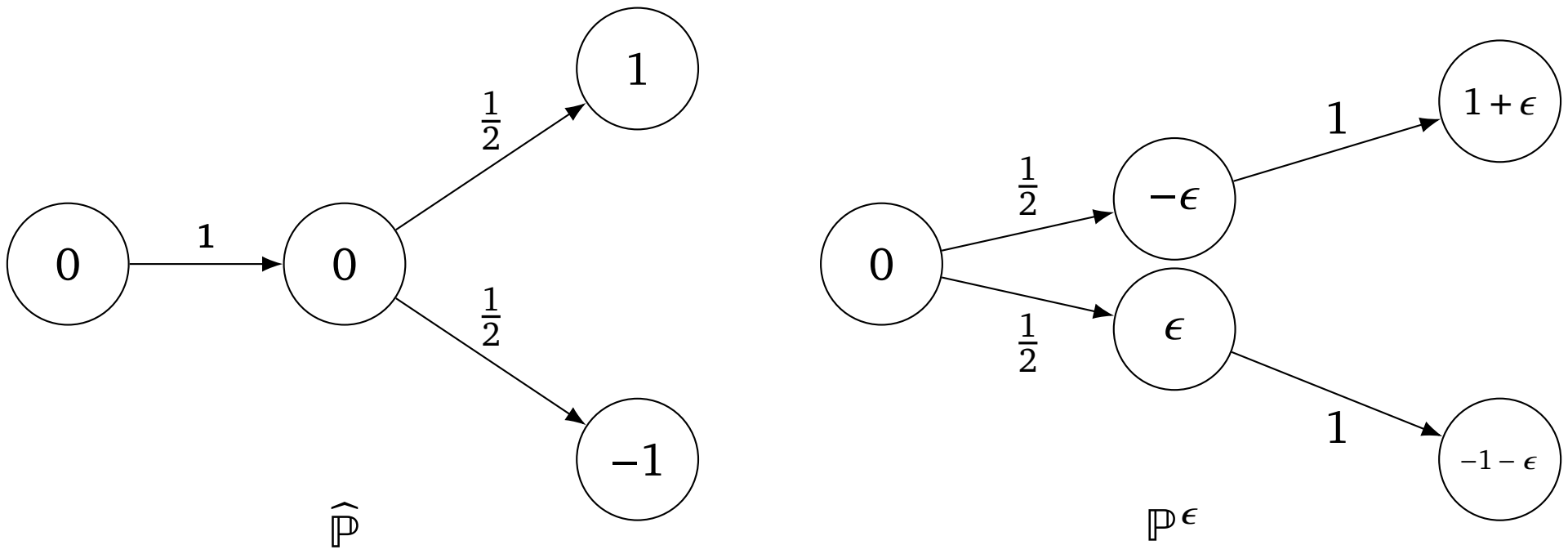

**Figure 2** The $p$-Wasserstein distance between $\widehat{\mathbb{P}}$ and $\mathbb{P}$ is $O(\epsilon)$ whereas the nested distance $\mathsf{D}_1(\widehat{\mathbb{P}}, \mathbb{P}) = O(1)$.

Example 2 (Tree structures in the Wasserstein and nested distance ball). (cf. Heitsch et al. (2006, Example 2.6), Pflug and Pichler (2012, Example 17)) Consider two stochastic processes $\widehat{\mathbb{P}}$ and $\mathbb{P}^\epsilon$ represented by the two scenario trees plotted in Figure 2:

$$\widehat{\mathbb{P}} = \frac{1}{2}\boldsymbol{\delta}_{(0,0,1)} + \frac{1}{2}\boldsymbol{\delta}_{(0,0,-1)},$$
$$\mathbb{P}^\epsilon = \frac{1}{2}\boldsymbol{\delta}_{(0,-\epsilon,1+\epsilon)} + \frac{1}{2}\boldsymbol{\delta}_{(0,\epsilon,-1-\epsilon)}.$$

The two processes $\widehat{\mathbb{P}}$ and $\mathbb{P}^\epsilon$ have different evolution of information. For $\widehat{\mathbb{P}}$, conditional on observing $\widehat{\boldsymbol{\xi}}_2 = 0$, $\widehat{\boldsymbol{\xi}}_3$ takes $\pm 1$ with equal probability; while for $\mathbb{P}$ with any $\epsilon > 0$, conditional on observing $\boldsymbol{\xi}_2$, the value of $\boldsymbol{\xi}_3$ is certain: $\mathbb{P}^\epsilon_{\xi_3 = 1+\epsilon \,|\, \xi_2 = -\epsilon} = \mathbb{P}^\epsilon_{\widehat{\xi}_3 = -1-\epsilon \,|\, \widehat{\xi}_2 = \epsilon} = 1$.

Suppose $\mathsf{d}(\widehat{\xi}_{[3]}, \xi_{[3]}) = \|\widehat{\xi}_{[3]} - \xi_{[3]}\|_p$. Then we have $\mathcal{W}_p(\widehat{\mathbb{P}}, \mathbb{P}^\epsilon) = 2^{1/p}\epsilon$ for all $p \in [1, \infty]$. On the other hand, a causal transport plan $\gamma$ from $\widehat{\mathbb{P}}$ to $\mathbb{P}^\epsilon$ should satisfy $\mathbb{P}^\epsilon_{\widehat{\xi}_3 | \widehat{\xi}_2 = 0} = \frac{1}{2} = \gamma_{\widehat{\xi}_3 | \widehat{\xi}_2 = 0, \xi_2}$. Hence the only feasible transport plan is the independent product distribution $\gamma(\widehat{\boldsymbol{\xi}}_{[3]}, \boldsymbol{\xi}_{[3]}) = \frac{1}{4}$ for all pairs of $\widehat{\boldsymbol{\xi}}_{[3]}, \boldsymbol{\xi}_{[3]}$, which is in fact bi-causal. Thus we have $\mathsf{D}_p(\widehat{\mathbb{P}}, \mathbb{P}^\epsilon) = \left(\frac{1}{4} \cdot 2(2\epsilon^p + \epsilon^p + (2+\epsilon)^p)\right)^{1/p} > 2^{\frac{p-1}{p}}$. As such, the nested-distance ball centered at $\widehat{\mathbb{P}}$ with radius 1 would not contain $\mathbb{P}^\epsilon$ with any $\epsilon > 0$.

Consider a function

$$Z(\xi_{[3]}) = |\xi_2 - \xi_3|.$$

Let $\mathfrak{M}_W$ be the $\infty$-Wasserstein ball centered at $\widehat{\mathbb{P}}$ with radius $\vartheta = \epsilon$, where $\epsilon \in (0,1)$. Using the duality for Wasserstein DRO (e.g., Zhang et al. (2022)), it is easy to check that the Wasserstein robust risk equals

$$\sup_{\mathbb{P} \in \mathfrak{M}_W} \mathbb{E}_{\mathbb{P}}[Z(\boldsymbol{\xi}_{[3]})] = \mathbb{E}_{\widehat{\mathbb{P}}}\left[\sup_{|\xi_2 - \widehat{\xi}_2| \le \epsilon, |\xi_3 - \widehat{\xi}_3| \le \epsilon} |\xi_2 - \xi_3|\right] = 1 + 2\epsilon,$$

with the worst-case distribution being $\mathbb{P}^\epsilon$. On the other hand, let $\mathfrak{M}_t$ be the $\infty$-Wasserstein ball centered at $\widehat{\mathbb{P}}_t$ with radius $\vartheta = \epsilon$, $t = 2, 3$. For the nested distance uncertainty set $\mathfrak{M}$, we will show that

$$\sup_{\mathbb{P} \in \mathfrak{M}} \mathbb{E}_{\mathbb{P}}[Z(\boldsymbol{\xi}_{[3]})] = \sup_{\mathbb{P}_2 \in \mathfrak{M}_2} \mathbb{E}_{\mathbb{P}_2}\left[\sup_{\mathbb{P}_3 \in \mathfrak{M}_3} \mathbb{E}_{\mathbb{P}_3}\left[|\boldsymbol{\xi}_2 - \boldsymbol{\xi}_3|\right]\right].$$

For any $\xi_2$, it holds that

$$\sup_{\mathbb{P}_3 \in \mathfrak{M}_3} \mathbb{E}_{\mathbb{P}_3}[|\xi_2 - \boldsymbol{\xi}_3|] = \mathbb{E}_{\widehat{\mathbb{P}}_3}[|\xi_2 - \widehat{\boldsymbol{\xi}}_3|] + \epsilon,$$

and

$$\sup_{|\xi_2 - 0| \le \epsilon} \mathbb{E}_{\widehat{\mathbb{P}}_3}\left[|\xi_2 - \widehat{\boldsymbol{\xi}}_3|\right] = \frac{1}{2} \sup_{|\xi_2| \le \epsilon} |\xi_2 - 1| + |\xi_2 + 1| = 1.$$

It follows that

$$\sup_{\mathbb{P} \in \mathfrak{M}} \mathbb{E}_{\mathbb{P}}[Z(\boldsymbol{\xi}_{[3]})] = \sup_{|\xi_2 - 0| \le \epsilon} \mathbb{E}_{\widehat{\mathbb{P}}_3}\left[|\xi_2 - \widehat{\boldsymbol{\xi}}_3|\right] + \epsilon = 1 + \epsilon,$$

with the worst-case distribution being

$$\frac{1}{2}\boldsymbol{\delta}_{(0,\xi,1+\epsilon)} + \frac{1}{2}\boldsymbol{\delta}_{(0,\xi,-1-\epsilon)}.$$

where $\xi \in [-\epsilon, \epsilon]$. Observe that this distribution has the same information structure as $\widehat{\mathbb{P}}$. ♣

Example 3 (Non-causal, causal, and bi-causal transport plans). In Figure 3 we plot two three-stage scenario trees, with labels along the edges indicating the conditional probabilities of realizing a scenario given their parent nodes. Specifically, define the sample paths

$$\widehat{\xi}^1_{[3]} = (\hat{a}, \hat{b}, \hat{d}),\ \widehat{\xi}^2_{[3]} = (\hat{a}, \hat{b}, \hat{e}),\ \widehat{\xi}^3_{[3]} = (\hat{a}, \hat{c}, \hat{f}),$$
$$\xi^1_{[3]} = (a, b, d),\ \xi^2_{[3]} = (a, c, e),\ \xi^3_{[3]} = (a, c, f),$$

where a sample path is represented by a triple of nodes. Then the two trees represent probability distributions on the three sample paths

$$\widehat{\mathbb{P}} = \frac{1}{6}\boldsymbol{\delta}_{\widehat{\xi}^1_{[3]}} + \frac{1}{6}\boldsymbol{\delta}_{\widehat{\xi}^2_{[3]}} + \frac{2}{3}\boldsymbol{\delta}_{\widehat{\xi}^3_{[3]}}, \quad \mathbb{P} = \frac{1}{2}\boldsymbol{\delta}_{\xi^1_{[3]}} + \frac{1}{4}\boldsymbol{\delta}_{\xi^2_{[3]}} + \frac{1}{4}\boldsymbol{\delta}_{\xi^3_{[3]}},$$

where $\boldsymbol{\delta}_\xi$ indicates a Dirac mass at a sample path $\xi$.

Consider the following three transport plans between the two trees, represented by a joint distribution with marginals $\widehat{\mathbb{P}}$ and $\mathbb{P}$:

$$\gamma = \begin{pmatrix} 1/6 & 0 & 0 \\ 0 & 1/6 & 0 \\ 1/3 & 1/12 & 1/4 \end{pmatrix}, \quad \gamma_c = \begin{pmatrix} 1/24 & 1/8 & 0 \\ 1/24 & 1/8 & 0 \\ 5/12 & 0 & 1/4 \end{pmatrix}, \quad \gamma_{bc} = \begin{pmatrix} 1/24 & 1/16 & 1/16 \\ 1/24 & 1/16 & 1/16 \\ 5/12 & 1/8 & 1/8 \end{pmatrix}.$$

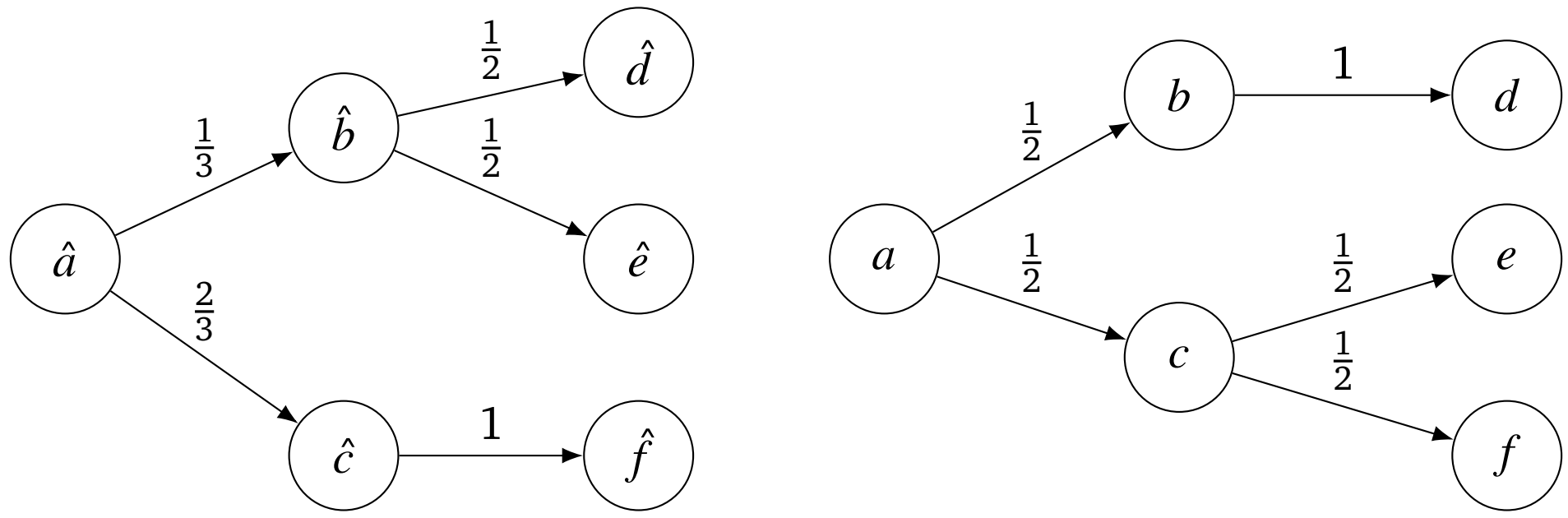

**Figure 3** Two three-stage scenario trees

In each matrix, its element in the $i$-th row and the $j$-th column represents the probability mass transported from the path $\widehat{\xi}^i_{[3]}$ to the path $\xi^j_{[3]}$, $i, j = 1, 2, 3$. We have the following observations.

(I) The transport plan $\gamma$ is not causal, as both causal constraints in (6) are violated. Indeed, $\widehat{\boldsymbol{\xi}}_2 = \widehat{b}$ is transported to $b$ if $\widehat{\boldsymbol{\xi}}_3 = \hat{d}$ and to $c$ if $\widehat{\boldsymbol{\xi}}_3 = \hat{e}$. This means that the future value of $\widehat{\boldsymbol{\xi}}_3$ affects the value of $\boldsymbol{\xi}_2$. On the other hand, $\boldsymbol{\xi}_2 = c$ is transported to $\hat{c}$ if $\boldsymbol{\xi}_3 = f$ and is split to $\hat{c}$ and $\hat{b}$ if $\boldsymbol{\xi}_3 = e$.

(II) The transport plan $\gamma_c$ is causal from $\widehat{\mathbb{P}}$ to $\mathbb{P}$, but not from $\mathbb{P}$ to $\widehat{\mathbb{P}}$. Indeed, regardless of the value of $\widehat{\boldsymbol{\xi}}_3$, $\hat{b}$ splits $1/24$ probability mass to $b$ and $1/8$ probability mass to $c$, and $\hat{c}$ splits $5/12$ probability mass to $b$ and $1/4$ probability mass to $c$. On the other hand, (6b) is violated, as $\boldsymbol{\xi}_2 = c$ is transported to $\hat{b}$ if $\boldsymbol{\xi}_3 = e$ and to $\hat{c}$ if $\boldsymbol{\xi}_3 = f$.

(III) The transport plan $\gamma_{bc}$ is bi-causal. Indeed, regardless of the value of $\widehat{\boldsymbol{\xi}}_3$, $\hat{b}$ splits $1/24$ probability mass to $b$ and $1/8$ probability mass to $c$, and $\hat{c}$ splits $5/12$ probability mass to $b$ and $1/4$ probability mass to $c$. Furthermore, regardless of the value of $\boldsymbol{\xi}_3$, $b$ splits $1/12$ probability mass to $\hat{b}$ and $5/12$ probability mass to $c$, and $c$ splits $1/8$ probability mass to $\hat{b}$ and $1/8$ probability mass to $\hat{c}$. ♣

**Setup**. In the rest of this paper, we will assume $\Xi_t$ is a non-empty subset of some normed space $(\mathbb{R}^{d_t}, \|\cdot\|)$, $t \in [T]$. Note that the norm in each stage can be chosen differently, but we omit such dependence for the ease of notation. We set

$$\mathsf{d}(\widehat{\xi}_{[T]}, \xi_{[T]}) = \begin{cases} \left( \sum_{t \in [T]} \|\widehat{\xi}_t - \xi_t\|^p \right)^{1/p}, & p \in [1, \infty), \\ \max_{t \in [T]} \|\widehat{\xi}_t - \xi_t\|, & p = \infty, \end{cases} \quad \forall\, \widehat{\xi}_{[T]}, \xi_{[T]} \in \Xi_{[T]},$$

and we use overloaded notation $\mathsf{d}(\widehat{\xi}_{[t]}, \xi_{[t]})$ for the distance between sub-sample paths.

### 2.3. Statistical Guarantees

Using recent concentration results on the nested distance (Backhoff et al. 2022, Acciaio and Hou 2022), in this subsection, we develop a finite-sample guarantee for (1). Suppose we are given $N$ i.i.d. sample paths from the underlying true distribution $\mathbb{P}^*$, let us choose the nominal process $\widehat{\mathbb{P}}$ as the adapted empirical measure constructed in (Acciaio and Hou 2022, Definition 2.5). Essentially, this is a finite scenario based on nearest-neighbor clustering of the empirical scenario tree. An illustration is given in Figure 4.

The subplot (a) shows a fan-shaped empirical scenario tree with an equally weighted distribution on 4 sample paths; the relative location of each node is consistent with the value of the outcomes associated with the each node. We cluster each node to one of the centroid ±5 based on the nearest-neighbor rule. For example, nodes with positive values are clustered into the center at 5. The resulting clustered scenario tree is shown in the subplot (b), with the edge probabilities calculated based on the clustering. The scenario tree (b) is called the adapted empirical measure. It has been shown that this measure converges to the underlying true distribution in nested distance as $N$ goes to infinity (even when the true distribution is continuous).

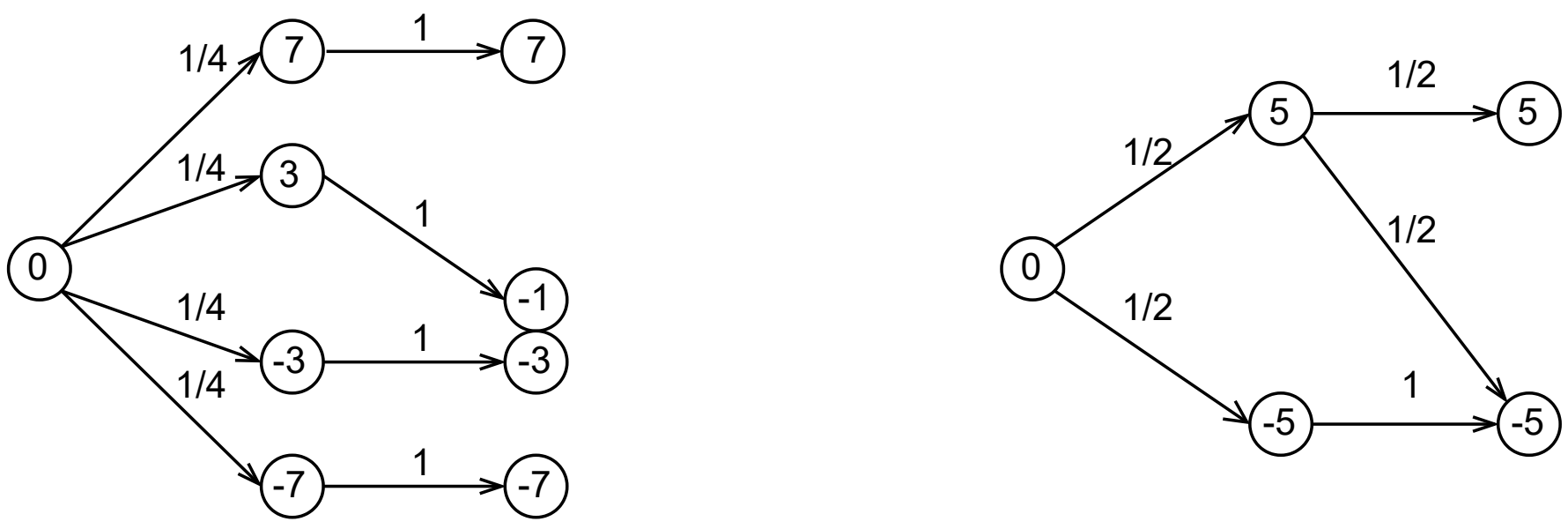

(a) (b)

**Figure 4** Illustration of adapted empirical measure (b) constructed from an empirical scenario tree (a)

PROPOSITION 1. *Let $p = 1$ and $d_t = d$, $t = 2, \dots, T$. Define $D(d) = d$ when $d \geq 3$ and $D(d) = d+1$ when $d = 1, 2$. Assume that*

(I) *There exists $\eta \geq 2$ such that $\mathbb{E}_{\mathbb{P}^*}[\exp(\|\boldsymbol{\xi}\|^{\eta})] < \infty$.*
(II) *There exists $L > 0$ such that for every $\boldsymbol{\xi}_{[t]}, \boldsymbol{\xi}'_{[t]} \in \Xi_{[t]}$,*

$$\mathcal{W}_1(\mathbb{P}^*_{t+1|\boldsymbol{\xi}_{[t]}}, \mathbb{P}^*_{t+1|\boldsymbol{\xi}'_{[t]}}) \leq L\|\boldsymbol{\xi}_{[t]} - \boldsymbol{\xi}'_{[t]}\|.$$

(III) *For $t \in [T-1]$, there exist $\mu_t : \Xi_{[t]} \to \mathbb{R}^d$, $\varsigma_t : \Xi_{[t]} \to \mathbb{R}_+$ and $\nu : \Xi_{[t]} \to \mathcal{P}(\mathbb{R}^d)$, such that for every $\xi_{[t]} \in \Xi_{[t]}$,*

$$\mathbb{P}^*_{t+1|\xi_{[t]}} \sim \mu_t(\xi_{[t]}) + \varsigma_t(\xi_{[t]})\varepsilon_{[t]},$$

*where $\varsigma_t(\xi_{[t]}) \leq k\|\xi_{[t]}\|^q$ for some $q \geq 0$ and $k > 0$, $\sup_{\xi_{[t]} \in \Xi_{[T]}} \mathbb{E}_{\nu(\xi_{[t]})}[\exp(\|\boldsymbol{\xi}_{t+1}\|^{\eta})] < \infty$, and $\varepsilon_{[t]} \sim \nu(\xi_{[t]})$.*

*Then there exists $c, C, K > 0$ such that for any $\beta \in (0, C/e^c)$, by setting $\widehat{\mathbb{P}}$ as the adapted empirical measure and*

$$\vartheta_N(\beta) = \begin{cases} \sqrt{\log(C/\beta)/c}\,N^{-1/2} + KN^{-\frac{1}{D(d)(T-1)}}, & q = 0, \\ (\log(C/\beta)/c)^{q+1}N^{-1/2} + KN^{-\frac{1}{D(d)(T-1)}}, & q > 0, \end{cases}$$

*it holds with probability at least $1 - \beta$ that*

$$\inf_{\boldsymbol{x} \in \mathcal{X}} \mathbb{E}_{\mathbb{P}^*}\left[\sum_{t \in [T]} \boldsymbol{c}_t^\top \boldsymbol{x}_t(\boldsymbol{\xi}_{[t]})\right] \leq \inf_{\boldsymbol{x} \in \mathcal{X}} \sup_{\mathbb{P} \in \mathfrak{M}} \mathbb{E}_{\mathbb{P}}\left[\sum_{t \in [T]} \boldsymbol{c}_t^\top \boldsymbol{x}_t(\boldsymbol{\xi}_{[t]})\right].$$

The assumption are adopted from Acciaio and Hou (2022, Theorem 2.19). Assumption (I) requires the ground truth $\mathbb{P}^*$ has finite exponential moment; Assumption (II) requires that the underlying transition probability kernel is $L$-Lipschitz in the history; and Assumption (III) means that the transition probability kernels have uniform exponential moments with polynomial growth rate. The proof follows directly from the concentration of the adapted empirical measure; see Appendix EC.1. This bound exhibits a non-parametric rate $O(N^{-\frac{1}{dT}})$, which is unavoidable without additional assumptions, given that we are optimizing over infinite-dimensional policies on an (infinite-dimensional) policies on an $O(dT)$-dimensional space. In Section 3.2, we will show that the dependence on $T$ can be further improved under the assumption of stagewise independence.

# 3. Dynamic Programming Reformulations

In this section, we develop dynamic programming reformulations for ($\mathsf{P}_{\mathsf{static}}$). In Section 3.1, we focus on the inner maximization of ($\mathsf{P}_{\mathsf{static}}$), which evaluates the worst-case risk of a fixed policy, and we discuss the outer minimization over policies in Section 3.2.

## 3.1. Robust Risk Evaluation

We first consider $p=\infty$. The following theorem provides a dynamic programming equivalent reformulation for the inner supremum in ($\mathsf{P}_{\mathsf{static}}$).

THEOREM 1. *Let $p=\infty$. Define $V_{T+1}^{\widehat{\boldsymbol{\xi}}_{[T+1]}} \equiv 0$ and*

$$V_t^{\widehat{\boldsymbol{\xi}}_{[t]}}(\boldsymbol{\xi}_{[t]}) := \boldsymbol{c}_t^\top \boldsymbol{x}_t(\boldsymbol{\xi}_{[t]}) + \mathbb{E}_{\widehat{\xi}_{t+1}\sim \widehat{\mathbb{P}}_{t+1|\widehat{\boldsymbol{\xi}}_{[t]}}}\left[\sup_{\xi_{t+1}\in\Xi_{t+1}:\|\xi_{t+1}-\widehat{\xi}_{t+1}\|\le\vartheta} V_{t+1}^{\widehat{\boldsymbol{\xi}}_{[t+1]}}(\boldsymbol{\xi}_{[t]},\xi_{t+1})\right],\quad t\in[T]. \tag{7}$$

*Assume $V_t^{\widehat{\boldsymbol{\xi}}_{[t]}}(\boldsymbol{\xi}_{[t-1]},\cdot)$ is lower semi-continuous. Then the inner supremum of* ($\mathsf{P}_{\mathsf{static}}$) *equals $V_1^{\widehat{\xi}_1}(\xi_1)$.*

In Theorem 1, the value function $V_t^{\widehat{\boldsymbol{\xi}}_{[t]}}(\cdot)$ is defined with respect to the nominal realization $\widehat{\boldsymbol{\xi}}_{[t]}$. It assesses the current cost along with the next-stage risk relative to the nominal conditional distribution $\widehat{\mathbb{P}}_{t+1|\widehat{\boldsymbol{\xi}}_{[t]}}$. More specifically, the risk $V_t^{\widehat{\boldsymbol{\xi}}_{[t]}}(\boldsymbol{\xi}_{[t]})$ at stage $t$ is broken down into two components: (i) the current-stage cost, $c_t^\top \boldsymbol{x}_t(\boldsymbol{\xi}_{[t]})$; (ii) the risk-to-go, which evaluates the worst-case value of the risk function $V_{t+1}^{\widehat{\boldsymbol{\xi}}_{[t+1]}}(\boldsymbol{\xi}_{[t]},\cdot)$ in a $\vartheta$-neighborhood of $\widehat{\boldsymbol{\xi}}_{t+1}$, and then averaged over the nominal conditional distribution $\widehat{\boldsymbol{\xi}}_{t+1}\sim\widehat{\mathbb{P}}_{t|\widehat{\boldsymbol{\xi}}_{[t]}}$. The continuity assumption on $V_t^{\widehat{\boldsymbol{\xi}}_{[t]}}(\cdot)$ means that given the history $\boldsymbol{\xi}_{[t-1]}$ and the nominal sample path $\widehat{\boldsymbol{\xi}}_{[t]}$, the value function is continuous with respect to the current-stage uncertainty $\boldsymbol{\xi}_t$. This is a mild assumption and is satisfied for the optimal value function, as will be shown in Section 3.2. When $\vartheta=0$, (7) reduces to the standard Bellman recursion. When $T=2$, the two-stage DRO with $\infty$-nested distance becomes the two-stage DRO with $\infty$-Wasserstein distance, and the result in Theorem 1 is consistent with Xie (2020, Theorem 2).

The proof idea of Theorem 1 can be summarized as follows. Due to the non-convexity of the nested distance uncertainty set, directly dualizing the inner supremum in ($\mathsf{P}_{\mathsf{static}}$) is not a viable approach. To overcome this challenge, we first consider a convex relaxation of ($\mathsf{P}_{\mathsf{static}}$), replacing the nested distance uncertainty set with the causal distance uncertainty set $\mathfrak{M}^{\mathsf{C}} := \{\mathbb{P}\in\mathcal{P}(\Xi):\ \mathsf{C}_\infty(\widehat{\mathbb{P}},\mathbb{P})\le\vartheta\}$. Using the conditional independence in the causal constraint (6a) and the tower property of conditional expectation, we are able to derive a dynamic programming reformulation for the relaxed problem. Next, we show that this convex relaxation is, in fact, tight. By modifying the worst-case distribution within the causal distance ball, it is shown that there exists a distribution whose nested distance to $\widehat{\mathbb{P}}$ is approximately equal to the causal distance, and it yields an objective value that is approximately equal to the worst-case risk over the causal distance ball, thanks to the continuity assumption on $V_t^{\widehat{\boldsymbol{\xi}}_{[t]}}$. Thereby, the expression holds for the nested distance ball as well. A complete proof can be found in Appendix EC.2.1. We remark that (Pflug and Pichler 2014, Section 7.3.1) also considers a convex relaxation of the nested distance uncertainty set, but their construction of the convex hull is based on compounding finite trees and does not lead to a dynamic programming recursion.

Next, we consider the case $p\in[1,\infty)$. We have the following result that establishes the nested reformulations of both (1) and its soft variant.

THEOREM 2. *Let $p \in [1,\infty)$. Set $V_{T+1}^{\widehat{\xi}_{[T+1]}} \equiv 0$ and*

$$V_t^{\widehat{\xi}_{[t]}}(\boldsymbol{\xi}_{[t]}) := \boldsymbol{c}_t^\top \boldsymbol{x}_t(\boldsymbol{\xi}_{[t]}) + \mathbb{E}_{\widehat{\xi}_{t+1}\sim\widehat{\mathbb{P}}_{t+1|\widehat{\xi}_{[t]}}}\left[\sup_{\xi_{t+1}\in\Xi_{t+1}}\left\{V_{t+1}^{\widehat{\xi}_{[t+1]}}(\boldsymbol{\xi}_{[t]},\xi_{t+1}) - \lambda\|\xi_{t+1}-\widehat{\boldsymbol{\xi}}_{t+1}\|^p\right\}\right], \quad t\in[T]. \tag{8}$$

*Assume $V_t^{\widehat{\xi}_{[t]}}(\boldsymbol{\xi}_{[t-1]},\cdot)$ is lower semi-continuous. Then the inner supremum of* ($\mathsf{P}_{\text{static}}$-soft) *equals $V_1^{\widehat{\xi}_1}(\xi_1)$ and the inner supremum of* ($\mathsf{P}_{\text{static}}$) *equals*

$$\min_{\lambda\geq 0}\left\{\lambda\vartheta^p + c_1^\top x_1 + \mathbb{E}_{\widehat{\mathbb{P}}_2}\left[\sup_{\xi_2\in\Xi_2}\left\{V_2^{\widehat{\xi}_2}(\xi_1,\xi_2) - \lambda\|\xi_2-\widehat{\boldsymbol{\xi}}_2\|^p\right\}\right]\right\}.$$

The reformulation (8) of the soft problem is often convenient to work with when it comes to time consistency; see Section 3.2. The interpretation of the risk-to-go function is similar to the previous case of $p=\infty$, except that the risk is now measured by the $p$-Wasserstein distance penalty. The proof idea is similar as well, except that the approximately worst-case distribution is constructed differently. A detailed proof can be found in Appendix EC.2.2.

Below, we discuss how the formulations derived in Theorems 1 and 2 are related to other formulations in the literature.

REMARK 1 (CONNECTION WITH DYNAMIC RISK MEASURES). The dynamic recursions (7) and (8) are conceptually related to the *conditional risk mappings* introduced in Ruszczyński and Shapiro (2006). Specifically, a coherent conditional risk mapping can be represented as

$$\rho_{t|\xi_{[t-1]}}[Z(\boldsymbol{\xi}_{[t]})] = \sup_{\mathbb{P}_t\in\mathfrak{M}_t}\mathbb{E}_{\xi_t\sim\mathbb{P}_t}[Z(\boldsymbol{\xi}_{[t-1]},\boldsymbol{\xi}_t)],$$

where $Z$ is any measurable function on $\Xi_{[t]}$, and a convex conditional risk mapping can be represented as

$$\rho_{t|\xi_{[t-1]}}[Z(\boldsymbol{\xi}_{[t]})] = \sup_{\mathbb{P}_t\in\mathcal{P}(\Xi_t)}\left\{\mathbb{E}_{\xi_t\sim\mathbb{P}_t}[Z(\boldsymbol{\xi}_{[t-1]},\boldsymbol{\xi}_t)] - \lambda J_t(\mathbb{P}_t)\right\},$$

where $J_t$ is a convex functional on $\mathcal{P}(\Xi_t)$. In our setting, let us define a mapping $\rho_{t|\xi_{[t-1]}}^{\widehat{\xi}_{[t-1]}}$ via

$$\begin{aligned}&\rho_{t|\xi_{[t-1]}}^{\widehat{\xi}_{[t-1]}}[Z^{\widehat{\xi}_{[t-1]}}(\boldsymbol{\xi}_{[t]})]\\ &:= \begin{cases}\sup_{\gamma_t\in\Gamma(\widehat{\mathbb{P}}_{t|\widehat{\xi}_{[t-1]}},\cdot)}\left\{\mathbb{E}_{(\widehat{\xi}_t,\xi_t)\sim\gamma_t}\left[Z^{\widehat{\xi}_{[t]}}(\boldsymbol{\xi}_{[t]})\right] : \gamma_t\text{-ess sup}\|\widehat{\boldsymbol{\xi}}_t-\boldsymbol{\xi}_t\|\leq\vartheta\right\}, & p=\infty,\\ \sup_{\gamma_t\in\Gamma(\widehat{\mathbb{P}}_{t|\widehat{\xi}_{[t-1]}},\cdot)}\left\{\mathbb{E}_{(\widehat{\xi}_t,\xi_t)\sim\gamma_t}\left[Z^{\widehat{\xi}_{[t]}}(\boldsymbol{\xi}_{[t]})\right] - \lambda\mathbb{E}_{(\widehat{\xi}_t,\xi_t)\sim\gamma_t}\left[\|\widehat{\boldsymbol{\xi}}_t-\boldsymbol{\xi}_t\|^p\right]\right\}, & p\in[1,\infty),\end{cases}\end{aligned} \tag{9}$$

where $\Gamma(\widehat{\mathbb{P}}_{t|\widehat{\xi}_{[t-1]}},\cdot)$ is the set of joint distributions on $\Xi_t^2$ whose first marginal distribution is $\widehat{\mathbb{P}}_{t|\widehat{\xi}_{[t-1]}}$. With this definition, (7) and (8) can be rewritten as

$$V_t^{\widehat{\xi}_{[t]}}(\boldsymbol{\xi}_{[t]}) = \boldsymbol{c}_t^\top\boldsymbol{x}_t(\boldsymbol{\xi}_{[t]}) + \rho_{t+1|\xi_{[t]}}^{\widehat{\xi}_{[t]}}[V_{t+1}^{\widehat{\xi}_{[t+1]}}(\boldsymbol{\xi}_{[t+1]})].$$

It is important to note that in the original definition of a conditional risk mapping, $\rho_{t|\xi_{[t-1]}}$ considers only the filtration generated by $\boldsymbol{\xi}_{[t]}$ However, in our case, it extends to involve the filtration generated by the nominal stochastic process $\widehat{\boldsymbol{\xi}}_{[T]}$. Note that in (9), when $p=\infty$, the uncertainty set is equivalent to a one-period $\infty$-Wasserstein ball $\mathfrak{M}_t^{\widehat{\xi}_{[t-1]}} = \{\mathbb{P}_t\in\mathcal{P}(\Xi_t) : \mathcal{W}_\infty(\widehat{\mathbb{P}}_{t|\widehat{\xi}_{[t-1]}},\mathbb{P}_t)\leq\vartheta\}$ centered at the nominal conditional distribution $\widehat{\mathbb{P}}_{t+1|\widehat{\xi}_{[t]}}$ with radius $\vartheta$, and $\rho_{t|\xi_{[t-1]}}^{\widehat{\xi}_{[t-1]}}$ defined a coherent risk measure; when $p\in[1,\infty)$, the penalty term is equivalent to a one-period $p$-Wasserstein distance penalty $J_t^{\widehat{\xi}_{[t-1]}}(\mathbb{P}_t) = \mathcal{W}_p(\widehat{\mathbb{P}}_{t|\widehat{\xi}_{[t-1]}},\mathbb{P}_t)^p$, and $\rho_{t|\xi_{[t-1]}}^{\widehat{\xi}_{[t-1]}}$ defines a convex risk measure. ♣

Remark 2 (Stagewise independent nominal tree). Suppose the nominal scenario tree $\widehat{\mathbb{P}}$ is stagewise independent, that is, $\widehat{\mathbb{P}}_{t+1|\widehat{\boldsymbol{\xi}}_{[t]}} = \widehat{\mathbb{P}}_{t+1}$ for all $t = 1, \ldots, T-1$. Then (7) and (8) become respectively

$$V_t(\boldsymbol{\xi}_{[t]}) := \begin{cases} \boldsymbol{c}_t^\top \boldsymbol{x}_t(\boldsymbol{\xi}_{[t]}) + \mathbb{E}_{\widehat{\xi}_{t+1}\sim\widehat{\mathbb{P}}_{t+1}}\left[\sup_{\xi_{t+1}\in\Xi_{t+1}:\|\xi_{t+1}-\widehat{\xi}_{t+1}\|\le\vartheta} V_{t+1}(\boldsymbol{\xi}_{[t]},\xi_{t+1})\right], & p=\infty, \\ \boldsymbol{c}_t^\top \boldsymbol{x}_t(\boldsymbol{\xi}_{[t]}) + \mathbb{E}_{\widehat{\xi}_{t+1}\sim\widehat{\mathbb{P}}_{t+1}}\left[\sup_{\xi_{t+1}\in\Xi_{t+1}}\left\{V_{t+1}(\boldsymbol{\xi}_{[t]},\xi_{t+1}) - \lambda\|\xi_{t+1}-\widehat{\boldsymbol{\xi}}_{t+1}\|^p\right\}\right], & p\in[1,\infty). \end{cases} \tag{10}$$

Set $\mathfrak{M}_t = \{\mathbb{P}_t \in \mathcal{P}(\Xi_t) : \mathcal{W}_\infty(\widehat{\mathbb{P}}_t, \mathbb{P}_t) \le \vartheta\}$, and

$$\rho_t[\cdot] := \begin{cases} \sup_{\mathbb{P}_t\in\mathfrak{M}_t} \mathbb{E}_{\mathbb{P}_t}[\cdot], & p=\infty, \\ \sup_{\mathbb{P}_t\in\mathcal{P}(\Xi_t)}\left\{\mathbb{E}_{\mathbb{P}_t}[\cdot] - \lambda\mathcal{W}_p(\widehat{\mathbb{P}}_t,\mathbb{P}_t)^p\right\}, & p\in[1,\infty), \end{cases} \tag{11}$$

Then using the duality for Wasserstein DRO (e.g., Zhang et al. (2022)), (10) can be rewritten as

$$V_t(\boldsymbol{\xi}_{[t]}) = \boldsymbol{c}_t^\top \boldsymbol{x}_t(\boldsymbol{\xi}_{[t]}) + \rho_{t+1}\left[V_{t+1}(\boldsymbol{\xi}_{[t]},\cdot)\right]. \tag{12}$$

In Section 3.2, we will discuss this setting in more detail. ♣

Remark 3 (Nominal tree with trivial filtration). Suppose $\widehat{\mathbb{P}}$ is a scenario fan, that is, $\widehat{\mathbb{P}}_{t+1|\widehat{\boldsymbol{\xi}}_{[t]}}$ is a Dirac measure for all $t = 1, \ldots, T-1$. Then by expanding (7) and (8), we obtain that

$$\sup_{\mathbb{P}\in\mathfrak{M}} \mathbb{E}_{\mathbb{P}}\left[\sum_{t\in[T]} \boldsymbol{c}_t^\top \boldsymbol{x}_t(\boldsymbol{\xi}_{[t]})\right] = c_1^\top x_1 + \mathbb{E}_{\widehat{\mathbb{P}}_2}\left[\sup_{\substack{\xi_t\in\Xi_t:\|\xi_t-\widehat{\xi}_t\|\le\vartheta \\ t=2,\ldots,T}} \sum_{t=2}^T \boldsymbol{c}_t^\top \boldsymbol{x}_t(\boldsymbol{\xi}_{[t]})\right]$$

and

$$\sup_{\mathbb{P}\in\mathcal{P}(\Xi)}\left\{\mathbb{E}_{\mathbb{P}}\Big[\sum_{t\in[T]} \boldsymbol{c}_t^\top \boldsymbol{x}_t(\boldsymbol{\xi}_{[t]})\Big] - \lambda \mathrm{D}_p^p(\widehat{\mathbb{P}},\mathbb{P})\right\} = c_1^\top x_1 + \mathbb{E}_{\widehat{\mathbb{P}}_2}\left[\sup_{\substack{\xi_t\in\Xi_t \\ t=2,\ldots,T}}\left\{\sum_{t=2}^T \boldsymbol{c}_t^\top \boldsymbol{x}_t(\boldsymbol{\xi}_{[t]}) - \lambda\sum_{t=2}^T\|\xi_t - \widehat{\boldsymbol{\xi}}_t\|^p\right\}\right].$$

Notably, they coincide with the dual formulation of Wasserstein DRO. Now suppose the underlying stochastic process has a density and let $\widehat{\mathbb{P}}$ be its empirical counterpart. Then with probability one, $\widehat{\mathbb{P}}$ is a scenario fan. As a result, under mild conditions on the value function, the statistical consistency of the worst-case risk under the nested distance follows from that of Wasserstein DRO (Esfahani and Kuhn 2018, Bertsimas et al. 2023), even though the nested distance itself does not have statistical consistency (Pflug and Pichler 2016, Proposition 1). ♣

Remark 4 (Rectangularity). Let $Z^{\boldsymbol{x}}(\boldsymbol{\xi}_{[T]})$ be the cumulative cost associated with a policy $\boldsymbol{x}$. Then (12) implies that

$$\sup_{\mathbb{P}\in\mathfrak{M}}[Z^{\boldsymbol{x}}] = \sup_{\mathbb{P}_2\in\mathfrak{M}_2}\mathbb{E}_{\mathbb{P}_2}\left[\sup_{\mathbb{P}_3\in\mathfrak{M}_3}\mathbb{E}_{\mathbb{P}_3}\Big[\cdots\sup_{\mathbb{P}_T\in\mathfrak{M}_T}\mathbb{E}_{\mathbb{P}_T}[Z^{\boldsymbol{x}}]\cdots\Big]\right].$$

If we set $\mathfrak{M}^\times := \left\{\widehat{\mathbb{P}}_1\times\mathbb{P}_2\times\mathbb{P}_3\times\cdots\times\mathbb{P}_T : \mathbb{P}_t\in\mathfrak{M}_t, t=2,\ldots,T\right\}$ and $\mathfrak{Z}$ be the set of cumulative cost functions of all feasible policies whose value function satisfies the continuity assumption in Theorem 1, then the above equality indicates that the $\infty$-nested distance set $\mathfrak{M}$ is a rectangular set associated with the product of one-period $\infty$-Wasserstein distance sets, $\mathfrak{M}^\times$, and $\mathfrak{Z}$, in the sense of Shapiro (2016). ♣

### 3.2. Time-Consistent Policy Optimization under Stagewise Independence

In numerical studies on multistage DRO, the following problem is often considered

$$\begin{aligned} \inf_{x_1\in\mathcal{X}_1} c_1^\top x_1 + \inf_{\boldsymbol{x}_2\in\mathcal{X}_2(x_1)} \rho_2\Bigg[c_2^\top \boldsymbol{x}_2(\boldsymbol{\xi}_2) + \inf_{\boldsymbol{x}_3\in\mathcal{X}_3(\boldsymbol{x}_2)} \rho_{3|\boldsymbol{\xi}_{[2]}}^{\widehat{\boldsymbol{\xi}}_{[2]}}\Big[c_3^\top \boldsymbol{x}_3(\boldsymbol{\xi}_{[3]}) + \cdots \\ + \inf_{\boldsymbol{x}_t\in\mathcal{X}_T(\boldsymbol{x}_{T-1})} \rho_{T|\boldsymbol{\xi}_{[T-1]}}^{\widehat{\boldsymbol{\xi}}_{[T-1]}}\left[c_T^\top \boldsymbol{x}_T(\boldsymbol{\xi}_{[T]})\right]\cdots\Big]\Bigg]. \end{aligned} \tag{13}$$

In general, this is just a lower bound of our considered problem (see Remark 1)

$$\inf_{\boldsymbol{x}_t\in\mathcal{X}_t(\boldsymbol{x}_{t-1}),\forall t\in[T]} c_1^\top x_1+\rho_2\left[c_2^\top \boldsymbol{x}_2(\boldsymbol{\xi}_2)+\rho_{3|\boldsymbol{\xi}_{[2]}}^{\widehat{\boldsymbol{\xi}}_{[2]}}\left[c_3^\top \boldsymbol{x}_3(\boldsymbol{\xi}_{[3]})+\cdots+\rho_{T|\boldsymbol{\xi}_{[T-1]}}^{\widehat{\boldsymbol{\xi}}_{[T-1]}}\left[c_T^\top \boldsymbol{x}_T(\boldsymbol{\xi}_{[T]})\right]\cdots\right]\right], \tag{14}$$

as it involves an exchange of minimization and expectation. Indeed, (13) allows different policy values under different realizations of the nominal process, whereas in (14) the policy does not depend on the realization of the nominal process.

Nevertheless, we will show that the two problems are equivalent when the nominal scenario tree $\widehat{\mathbb{P}}$ is stagewise independent.

Assumption 1 **(Stagewise independence)**. *The nominal scenario tree $\widehat{\mathbb{P}}$ is stagewise independent, namely, $\widehat{\mathbb{P}}_{t|\widehat{\boldsymbol{\xi}}_{[t-1]}}=\widehat{\mathbb{P}}_t$ for all $t\in[T]$.*

In the remainder of the paper, we make the following assumption (c.f. Shapiro et al. (2021, Definition 3.1)).

Assumption 2 **(Relatively complete recourse)**. *For every $t=2,\dots,T$ and every sequence of feasible decisions $(\boldsymbol{x}_1,\dots,\boldsymbol{x}_{t-1})$, the set $\mathcal{X}_t(\boldsymbol{x}_{t-1},\xi_t)$ is non-empty for every $\xi_t\in\Xi_t$.*

Without this assumption, the worst-case risk is always infinite for $p\in[1,\infty)$. And for $p=\infty$, $\Xi_t$ can be replaced with $\bigcup_{\widehat{\xi}_t\in\operatorname{supp}\widehat{\mathbb{P}}_t}\{\xi_t\in\Xi_t:\|\xi_t-\widehat{\xi}_t\|\le\vartheta\}$.

The following assumption is needed for applying existing Wasserstein DRO duality in our analysis.

Assumption 3 **(Sufficiently expensive recourse)**. *For any feasible $x_1$, we have*

$$\inf_{\boldsymbol{x}_2,\boldsymbol{x}_3,\dots\boldsymbol{x}_T\in\mathcal{X}}\mathbb{E}_{\widehat{\mathbb{P}}}\left[\sum_{t=2}^{T}\boldsymbol{c}_t^\top\boldsymbol{x}_t(\widehat{\boldsymbol{\xi}}_{[t]})\right]>-\infty.$$

In other words, this subproblem started from the second stage is dual feasible. For a two-stage problem, this assumption reduces to

$$\inf_{\boldsymbol{x}_2\in\mathcal{X}_2(\boldsymbol{x}_1,\widehat{\xi}_2)}\boldsymbol{c}_2^\top\boldsymbol{x}_2(\boldsymbol{\xi}_{[2]})>-\infty,$$

for any feasible $x_1$, which is the same as the sufficiently expensive recourse assumption in previous works (Hanasusanto and Kuhn 2018, Xie 2020). Since $\widehat{\mathbb{P}}$ is finitely supported, by Assumption 3, it is easy to show by induction that for any $s\in[T-1]$, for any feasible $x_{[s]}$, we have

$$\inf_{\boldsymbol{x}_{s+1},\dots\boldsymbol{x}_T\in\mathcal{X}}\mathbb{E}_{\widehat{\mathbb{P}}}\left[\sum_{t=s+1}^{T}\boldsymbol{c}_t^\top\boldsymbol{x}_t(\widehat{\boldsymbol{\xi}}_{[t]})\right]>-\infty.$$

This result helps establish strong duality in the proof of Theorem 3 and the corollaries in Section 4.

Define the dynamic programming formulation

$$\begin{aligned}
\mathcal{Q}_{T+1}&:=0,\\
\mathcal{Q}_t(x_{t-1},\boldsymbol{\xi}_{[t-1]})&:=\inf_{\boldsymbol{x}_t(\boldsymbol{\xi}_{[t]})\in\mathcal{X}_t(x_{t-1},\cdot)}\rho_t\left[\boldsymbol{c}_t^\top\boldsymbol{x}_t(\boldsymbol{\xi}_{[t]})+\mathcal{Q}_{t+1}\big(\boldsymbol{x}_t(\boldsymbol{\xi}_{[t]}),\boldsymbol{\xi}_{[t]}\big)\right],\\
\mathcal{Q}_1&:=\inf_{x_1\in\mathcal{X}_1}\left\{c_1^\top x_1+\mathcal{Q}_2(x_1)\right\},
\end{aligned}\tag{$\mathsf{P}_{\text{dynamic}}$}$$

where $\rho_t$ is defined in (11).

In the next theorem, we show that this dynamic programming formulation is equivalent to ($\mathsf{P}_{\text{static}}$) when $p=\infty$ and ($\mathsf{P}_{\text{static}}$-soft) when $p\in[1,\infty)$. The proof can be found in Appendix EC.3.

Theorem 3. *Suppose Assumptions 1, 2 and 3 hold. Then the optimal value of the dynamic program* ($\mathsf{P}_{\mathsf{dynamic}}$) *equals the optimal value of* ($\mathsf{P}_{\mathsf{static}}$) *when* $p=\infty$ *and the optimal value of* ($\mathsf{P}_{\mathsf{static}}$-soft) *when* $p\in[1,\infty)$. *Moreover, set* $Q_{T+1}\equiv 0$ *and*

$$Q_t(x_{t-1},\widehat{\boldsymbol{\xi}}_t) := \begin{cases} \sup_{\xi_t\in\Xi_t:\|\xi_t-\widehat{\xi}_t\|\le\vartheta} \inf_{x_t\in\mathcal{X}_t(x_{t-1},\xi_t)} \left\{c_t^\top x_t + \mathbb{E}_{\widehat{\mathbb{P}}_{t+1}}\left[Q_{t+1}(x_t,\widehat{\boldsymbol{\xi}}_{t+1})\right]\right\}, & p=\infty,\\ \sup_{\xi_t\in\Xi_t} \inf_{x_t\in\mathcal{X}_t(x_{t-1},\xi_t)} \left\{c_t^\top x_t + \mathbb{E}_{\widehat{\mathbb{P}}_{t+1}}\left[Q_{t+1}(x_t,\widehat{\boldsymbol{\xi}}_{t+1})\right] - \lambda\|\xi_t-\widehat{\xi}_t\|^p\right\}, & p\in[1,\infty), \end{cases} \tag{15}$$

*for* $t=2,\dots,T$. *Then it holds that* $\mathcal{Q}_t(x_{t-1},\boldsymbol{\xi}_{[t-1]}) = \mathbb{E}_{\widehat{\mathbb{P}}_t}[Q_t(x_{t-1},\widehat{\boldsymbol{\xi}}_t)]$.

Theorem 3 shows that expected risk-to-go function $\mathcal{Q}_t$ defined in ($\mathsf{P}_{\mathsf{dynamic}}$), which involves functional optimization over $\boldsymbol{x}_t(\boldsymbol{\xi}_{[t-1]},\cdot)$, can be evaluated through the expectation of the risk-to-go function $Q_t$, which involves only a finite-dimensional problem. Note that due to the stagewise independence Assumption 1, the risk-to-go function $Q_t(x_{t-1},\cdot)$ depends only on the current-stage uncertainty and the expected risk-to-go function $\mathcal{Q}_t(x_{t-1},\cdot)$ is a constant function. Thus, we can omit the second argument of $\mathcal{Q}_t$ and denote it as $\mathcal{Q}_t(x_{t-1})$. Also note that when $p\in$ only the soft-formulation

Remark 5 (Optimal time-consistent policy for ($\mathsf{P}_{\mathsf{dynamic}}$)). Suppose the optimal policy for ($\mathsf{P}_{\mathsf{dynamic}}$) exists. Then an optimal policy $\boldsymbol{x}^\star=(\boldsymbol{x}_1^\star,\dots,\boldsymbol{x}_T^\star)$ is given recursively by

$$\begin{aligned} \boldsymbol{x}_1^\star &\in \operatorname*{arg\,min}_{x_1\in\mathcal{X}_1} \left\{c_1^\top x_1 + \mathbb{E}_{\widehat{\mathbb{P}}_2}[Q_2(x_1)]\right\},\\ \boldsymbol{x}_t^\star(\xi_{[t]}) &\in \operatorname*{arg\,min}_{x_t\in\mathcal{X}_t(\boldsymbol{x}_{t-1}^\star(\xi_{[t-1]}),\xi_t)} \left\{c_t^\top x_t + \mathbb{E}_{\widehat{\mathbb{P}}_{t+1}}\left[Q_{t+1}(x_t,\widehat{\boldsymbol{\xi}}_{t+1})\right]\right\}, \quad \xi_{[t]}\in\Xi_{[t]}, t=2,\dots,T. \end{aligned} \tag{16}$$

This policy $\boldsymbol{x}^\star$ is time-consistent, satisfying the robust Bellman equation (10). It is important to note that $\boldsymbol{x}^\star$ is defined on the entire sample space, marking a significant difference from its non-robust counterpart, which yields a solution defined solely on the sample paths within the nominal scenario tree. The optimal robust value function associated the dynamic formulation satisfies lower semi-continuity condition in Theorems 1 and 2 (see Proposition EC.5 in Appendix EC.3), thereby we conclude that $\boldsymbol{x}^\star$ has the same worst-case expected cost in ($\mathsf{P}_{\mathsf{static}}$) and ($\mathsf{P}_{\mathsf{static}}$-soft), respectively, which means that $\boldsymbol{x}^\star$ is the optimal policy for these problems. ♣

The next example demonstrates that the formulation with Wasserstein distance is not time consistent even under stagewise independence.

Example 4 (Time inconsistency of formulations with Wasserstein distance). Suppose $T=3$, $c_1=c_2=0$, $c_3=1$, $\Xi_2=\mathbb{R}_+$, $\Xi_3=\mathbb{R}$, and let $\mathcal{X}_2(x_1,\xi_2)=\{x_2=(x_{21},x_{22})\in\mathbb{R}_+^2 : x_{21}-x_{22}=\xi_2\}$, $\mathcal{X}_3(x_2,\xi_3)=\{x_3\ge 0:\ x_3\ge x_{21}-x_{22}-\xi_3,\ x_3\ge \xi_3-x_{21}+x_{22}\}$. Let $\widehat{\mathbb{P}}=\widehat{\mathbb{P}}_2\times\widehat{\mathbb{P}}_3$, where $\widehat{\mathbb{P}}_2=\boldsymbol{\delta}_0$, $\widehat{\mathbb{P}}_3=\frac{1}{2}\boldsymbol{\delta}_{-1}+\frac{1}{2}\boldsymbol{\delta}_1$. Let $\vartheta\in(0,1)$. Define $\mathfrak{M}_W$ as the $\infty$-Wasserstein ball centered at $\widehat{\mathbb{P}}$ with radius $\vartheta$, and $\mathfrak{M}_t$ as the $\infty$-Wasserstein ball centered at $\widehat{\mathbb{P}}_t$ with radius $\vartheta$. The formulation (1) with $\infty$-Wasserstein ball reads

$$\inf_{\boldsymbol{x}_2\in\mathcal{X}_2(x_1),\boldsymbol{x}_3\in\mathcal{X}_3(\boldsymbol{x}_2)} \sup_{\mathbb{P}\in\mathfrak{M}_W} \mathbb{E}_{\mathbb{P}}[\boldsymbol{x}_3(\boldsymbol{\xi}_3)], \tag{17}$$

and the min-max-min-max dynamic formulation ($\mathsf{P}_{\mathsf{dynamic}}$) with $\infty$-Wasserstein distance reads

$$\inf_{\boldsymbol{x}_2\in\mathcal{X}_2(x_1)} \sup_{\mathbb{P}_2\in\mathfrak{M}_2} \mathbb{E}_{\mathbb{P}_2}\left[\inf_{\boldsymbol{x}_3\in\mathcal{X}_3(\boldsymbol{x}_2)} \sup_{\mathbb{P}_3\in\mathfrak{M}_3} \mathbb{E}_{\mathbb{P}_3}\left[\boldsymbol{x}_3(\boldsymbol{\xi}_{[3]})\right]\right]. \tag{18}$$

We have that

$$V_3(x_2,\xi_3) = \inf_{x_3\in\mathcal{X}_3(x_2,\xi_3)} x_3 = |x_{21}-x_{22}-\xi_3|.$$

Using Wasserstein DRO duality (Zhang et al. 2022), (17) is equivalent to

$$\inf_{\boldsymbol{x}_2} \mathbb{E}_{\widehat{\mathbb{P}}}\left[\sup_{|\xi_2-\widehat{\xi}_2|\le\vartheta,|\xi_3-\widehat{\xi}_3|\le\vartheta} |\boldsymbol{x}_{21}(\xi_2)-\boldsymbol{x}_{22}(\xi_2)-\xi_3|\right] = 2+\vartheta.$$

On the other hand, (18) is equivalent to

$$\inf_{\boldsymbol{x}_2} \sup_{\mathbb{P}_2 \in \mathfrak{M}_2} \left[ \sup_{\mathbb{P}_3 \in \mathfrak{M}_3} \mathbb{E}_{\mathbb{P}_3}[V_3(\boldsymbol{x}_2(\boldsymbol{\xi}_2), \xi_3)] \right] = \mathbb{E}_{\widehat{\mathbb{P}}_2} \left[ \sup_{|\xi_2 - \widehat{\xi}_2| \le \vartheta} \mathbb{E}_{\widehat{\mathbb{P}}_3} \left[ |\xi_2 - \widehat{\boldsymbol{\xi}}_3| \right] + \vartheta \right] = 1 + \vartheta.$$

This demonstrates that (17) is not time consistent in the sense of Shapiro (2016). ♣

Thanks to Theorem 3, using statistical guarantees for single-stage Wasserstein DRO (Kuhn et al. 2019), we can improve the bound in the finite sample guarantee when there is stagewise independence. Suppose that for each stage $t = 2, \ldots, T$, we are given $N$ i.i.d. samples that forms an empirical distribution $\widehat{\mathbb{P}}_t$.

Proposition 2. *Let $p \in [1, \infty)$ and $d_t = d$, $t = 2, \ldots, T$. Suppose the true distribution $\mathbb{P}^*$ has stagewise independent marginals $\mathbb{P}_t^*$. Assume there exists $\eta > p$ such that $\mathbb{E}_{\mathbb{P}_t^*}[\exp(\|\boldsymbol{\xi}_t\|^\eta)] < \infty$. Let $\beta \in (0, 1)$. Then there exists $c, C > 0$ such that for any $N \ge \frac{\log(C\beta^{-1})}{c}$, by setting $\widehat{\mathbb{P}} = \times_{t=2}^T \widehat{\mathbb{P}}_t$ and*

$$\vartheta_N(\beta) = (T-1)^{1/p} \left( \frac{\log(C(T-1)/\beta)}{cN} \right)^{\min(\frac{p}{d}, \frac{1}{2})},$$

*it holds with probability at least $(1 - \beta)$ that*

$$\inf_{\boldsymbol{x} \in \mathcal{X}} \mathbb{E}_{\mathbb{P}^*} \left[ \sum_{t \in [T]} \boldsymbol{c}_t^\top \boldsymbol{x}_t(\boldsymbol{\xi}_{[t]}) \right] \le \inf_{\boldsymbol{x} \in \mathcal{X}} \sup_{\mathbb{P} \in \mathfrak{M}} \mathbb{E}_{\mathbb{P}} \left[ \sum_{t \in [T]} \boldsymbol{c}_t^\top \boldsymbol{x}_t(\boldsymbol{\xi}_{[t]}) \right].$$

The proof for the proposition is provided in Appendix EC.4. Compared with Proposition 1 that suffers from curse of dimensionality in both $d$ and $T$, the bound in Proposition 2 has only a nearly-linear dependence on $T$ when $p = 1$. Besides, Proposition 2 applies for all $p \in [1, \infty)$.

## 4. Tractable Policy Optimization

In Section 3.2, we have derived dynamic programming reformulation ($\mathsf{P}_{\text{dynamic}}$) to solve for the optimal policy. Observe that the maximization over $\xi_t$ in (15) may not always be tractable. Nevertheless, in Section 4, we will explore cases where tractable solutions are possible, by identifying cases where the risk-to-go function $Q_t$ can be computed efficiently. Note that when $T = 2$, (1) (or ($\mathsf{P}_{\text{dynamic}}$)) reduces to two-stage Wasserstein DRO, for which exact tractable reformulations have been established when $p \in \{1, \infty\}$ and the uncertainty appears in the objective or right-hand side (Esfahani and Kuhn 2018, Hanasusanto and Kuhn 2018, Xie 2020). Below, we will show that these results can be extended to multistage problems. All results assume Assumptions 1, 2, 3 in the previous section hold.

### 4.1. Objective Uncertainty Only

We first consider problems with objective uncertainty only, in which case we identify $\boldsymbol{\xi}_t$ with $\boldsymbol{c}_t$, and the constraint set $\mathcal{X}_t(x_{t-1}, \xi_t) = \mathcal{X}_t(x_{t-1}) = \{x_t \ge 0 : A_t x_t = b_t - B_t x_{t-1}\}$ is deterministic once $x_{t-1}$ is given. The following theorem shows an equivalent expression of (15).

Corollary 1. *Suppose $\Xi_t = (\mathbb{R}^{d_t}, \|\cdot\|)$ for $t = 2, \ldots, T$. Then $Q_t$ defined in (15) can be computed as*

$$Q_t(x_{t-1}, \widehat{\boldsymbol{c}}_t) = \begin{cases} \inf_{x_t \in \mathcal{X}_t(x_{t-1}), \|x_t\|_* \le \lambda} \left\{ \widehat{\boldsymbol{c}}_t^\top x_t + \mathbb{E}_{\widehat{\mathbb{P}}_{t+1}} \left[ Q_{t+1}(x_t, \widehat{\boldsymbol{c}}_{t+1}) \right] \right\}, & p = 1, \\ \inf_{x_t \in \mathcal{X}_t(x_{t-1})} \left\{ \widehat{\boldsymbol{c}}_t^\top x_t + (1 - 1/p)(\frac{1}{p\lambda})^{\frac{1}{p-1}} \|x_t\|_*^{\frac{p}{p-1}} + \mathbb{E}_{\widehat{\mathbb{P}}_{t+1}} \left[ Q_{t+1}(x_t, \widehat{\boldsymbol{c}}_{t+1}) \right] \right\}, & p \in (1, \infty), \\ \inf_{x_t \in \mathcal{X}_t(x_{t-1})} \left\{ \widehat{\boldsymbol{c}}_t^\top x_t + \vartheta \|x_t\|_* + \mathbb{E}_{\widehat{\mathbb{P}}_{t+1}} \left[ Q_{t+1}(x_t, \widehat{\boldsymbol{c}}_{t+1}) \right] \right\}, & p = \infty. \end{cases}$$

Corollary 1 shows that when the uncertainty appears only in the objective, (10) is equivalent to a scenario approximation problem with norm regularization on each $\boldsymbol{x}_t$. The regularization is a hard constraint when $p=1$ and a soft penalty when $p\in(1,\infty]$. The proof is based on an additional duality argument; see EC.5.1 for details.

We remark that the optimal robust policy may not be unique, as demonstrated in the following result.

Proposition 3. *Set*

$$\widehat{\boldsymbol{x}}_t(\widehat{\boldsymbol{c}}_{[t]}) \in \begin{cases} \arg\min_{x_t\in\mathcal{X}_t(\widehat{\boldsymbol{x}}_{t-1}(\widehat{\boldsymbol{c}}_{[t-1]})),\|x_t\|_*\le\lambda} \left\{\widehat{\boldsymbol{c}}_t^\top x_t + \mathbb{E}_{\widehat{\mathbb{P}}_{t+1}}\left[Q_{t+1}(x_t,\widehat{\boldsymbol{c}}_{t+1})\right]\right\}, & p=1,\\ \arg\min_{x_t\in\mathcal{X}_t(\widehat{\boldsymbol{x}}_{t-1}(\widehat{\boldsymbol{c}}_{[t-1]}))} \left\{\widehat{\boldsymbol{c}}_t^\top x_t + (1-\frac{1}{p})(\frac{1}{p\lambda})^{\frac{1}{p-1}}\|x_t\|_*^{\frac{p}{p-1}} + \mathbb{E}_{\widehat{\mathbb{P}}_{t+1}}\left[Q_{t+1}(x_t,\widehat{\boldsymbol{c}}_{t+1})\right]\right\}, & p\in(1,\infty),\\ \arg\min_{x_t\in\mathcal{X}_t(\widehat{\boldsymbol{x}}_{t-1}(\widehat{\boldsymbol{c}}_{[t-1]}))} \left\{\widehat{\boldsymbol{c}}_t^\top x_t + \vartheta\|x_t\|_* + \mathbb{E}_{\widehat{\mathbb{P}}_{t+1}}\left[Q_{t+1}(x_t,\widehat{\boldsymbol{c}}_{t+1})\right]\right\}, & p=\infty, \end{cases}$$

*and set recursively* $\widehat{\boldsymbol{c}}_1^c := c_1$ *and*

$$\widehat{\boldsymbol{c}}_t^c := \begin{cases} \underset{\widehat{\boldsymbol{c}}_t\in\operatorname{supp}\widehat{\mathbb{P}}_t}{\arg\min} \left\{\widehat{\boldsymbol{c}}_t^\top\widehat{\boldsymbol{x}}_t(\widehat{\boldsymbol{c}}_{[t-1]}^c,\widehat{\boldsymbol{c}}_t) + \mathbb{E}_{\widehat{\mathbb{P}}_{t+1}}\left[Q_{t+1}\big(\widehat{\boldsymbol{x}}_t(\widehat{\boldsymbol{c}}_{[t-1]}^c,\widehat{\boldsymbol{c}}_t),\widehat{\boldsymbol{c}}_{t+1}\big)\right] + \lambda\|\widehat{\boldsymbol{c}}_t - c_t\|^p\right\}, & p\in[1,\infty),\\ \underset{\widehat{\boldsymbol{c}}_t\in\operatorname{supp}\widehat{\mathbb{P}}_t:\|\widehat{\boldsymbol{c}}_t-c_t\|\le\vartheta}{\arg\min} \widehat{\boldsymbol{c}}_t^\top\widehat{\boldsymbol{x}}_t(\widehat{\boldsymbol{c}}_{[t-1]}^c,\widehat{\boldsymbol{c}}_t) + \vartheta\|\widehat{\boldsymbol{x}}_t(\widehat{\boldsymbol{c}}_{[t-1]}^c,\widehat{\boldsymbol{c}}_t)\|_* + \mathbb{E}_{\widehat{\mathbb{P}}_{t+1}}\left[Q_{t+1}\big(\widehat{\boldsymbol{x}}_t(\widehat{\boldsymbol{c}}_{[t-1]}^c,\widehat{\boldsymbol{c}}_t),\widehat{\boldsymbol{c}}_{t+1}\big)\right], & p=\infty. \end{cases}$$

*Then the policy* $\bar{\boldsymbol{x}}_t(c_{[t]}) := \widehat{\boldsymbol{x}}_t(\widehat{\boldsymbol{c}}_{[t]}^c)$ *is also optimal for* (15).

The policy $\bar{\boldsymbol{x}} = (\bar{\boldsymbol{x}}_1,\ldots,\bar{\boldsymbol{x}}_T)$ is defined for every sample path in $\Xi_{[T]}$ when $p\in[1,\infty)$ and in $\{c_{[T]}\in\Xi_{[T]} : \|c_{[T]} - \widehat{\boldsymbol{c}}_{[T]}\|_\infty \le \vartheta\}$ when $p=\infty$. The sample path $\widehat{\boldsymbol{c}}_{[T]}^c$ represents the best in-sample path, with regard to the norm-regularized cost-to-go, within a $\vartheta$-neighborhood of $c_{[T]}$ when $p=\infty$, or within a $\lambda$-soft neighborhood of $c_{[T]}$ when $p\in[1,\infty)$. Notably, the computation of the policy $\bar{\boldsymbol{x}}$ requires knowledge of only the optimal robust policy values on sample paths from the nominal scenario tree, contrasting with the policy $\boldsymbol{x}^\star$ defined in (16) that requires the entire cost-to-go function $Q_t(\cdot,\widehat{\boldsymbol{c}}_t)$. Its optimality is obtained by verifying that the worst-case risk of $\bar{\boldsymbol{x}}_t$ does not exceed the risk-to-go as defined by $Q_t$. For a detailed proof, please refer to Appendix EC.5.1.

### 4.2. Right-hand Side Uncertainty Only

Next, we consider problems with right-hand side uncertainty only. To ease the presentation, we consider either $\boldsymbol{\xi}_t = \boldsymbol{b}_t$ or $\boldsymbol{\xi}_t = \boldsymbol{B}_t$. The following result provides an equivalent reformulation of (15) for $p=\infty$.

Corollary 2. *Suppose* $p=\infty$, *and* $\boldsymbol{\xi}_t = \boldsymbol{b}_t$, $\Xi_t = (\mathbb{R}^{d_t},\|\cdot\|)$, *or* $\boldsymbol{\xi}_t = \boldsymbol{B}_t$, $\Xi_t = (\mathbb{R}^{d_t\times m_t},\|\cdot\|)$, $t=2,\ldots,T$. *Set* $\psi_t(x_t) := c_t^\top x_t + \mathbb{E}_{\widehat{\mathbb{P}}_{t+1}}\left[Q_{t+1}(x_t,\widehat{\boldsymbol{\xi}}_{t+1})\right] + \mathbb{I}\{x_t\ge 0\}$. *Then* $Q_t$ *defined in* (15) *can be computed as*

$$\begin{aligned} Q_t(x_{t-1},\widehat{\boldsymbol{b}}_t) &= \max_{y_t\in\mathbb{R}^{d_t}} \left\{(\widehat{\boldsymbol{b}}_t - B_t x_{t-1})^\top y_t + \vartheta\|y_t\|_* - \psi_t^*(A_t^\top y_t)\right\}, && \text{if } \boldsymbol{\xi}_t = \boldsymbol{b}_t,\\ Q_t(x_{t-1},\widehat{\boldsymbol{B}}_t) &= \max_{y_t\in\mathbb{R}^{d_t}} \left\{(b_t - \widehat{\boldsymbol{B}}_t x_{t-1})^\top y_t + \vartheta\|y_t x_{t-1}^\top\|_* - \psi_t^*(A_t^\top y_t)\right\}, && \text{if } \boldsymbol{\xi}_t = \boldsymbol{B}_t. \end{aligned} \tag{19}$$

*When* $\|\cdot\| = \|\cdot\|_1$, *it holds that*

$$Q_t(x_{t-1},\widehat{\boldsymbol{b}}_t) = \max_{j\in[d_t],\delta\in\{1,-1\}} \inf_{x_t\ge 0} \left\{c_t^\top x_t + \mathbb{E}_{\widehat{\mathbb{P}}_{t+1}}\left[Q_{t+1}(x_t,\widehat{\boldsymbol{b}}_{t+1})\right] : A_t x_t = \widehat{\boldsymbol{b}}_t - B_t x_{t-1} + \vartheta\delta e_j\right\},$$

*where* $e_j$ *is the* $j$*-th unit vector. When* $\|\cdot\| = \|\cdot\|_{\text{op}}$, *where* $\|B\|_{\text{op}} = \sup_{\|v\|\le 1}\|Bv\|_1$, *it holds that*

$$Q_t(x_{t-1},\widehat{\boldsymbol{B}}_t) = \max_{j\in[d_t],\delta\in\{1,-1\}} \inf_{x_t\ge 0} \left\{c_t^\top x_t + \mathbb{E}_{\widehat{\mathbb{P}}_{t+1}}\left[Q_{t+1}(x_t,\widehat{\boldsymbol{B}}_{t+1})\right] : A_t x_t = b_t - \widehat{\boldsymbol{B}}_t x_{t-1} + \vartheta\delta\|x_{t-1}\|e_j\right\}.$$

The risk-to-go function $Q_t$ encourages a large norm of the dual variable $y_t$. It also penalizes a large norm on the primal variable $x_t$ when the uncertainty is present in $\boldsymbol{B}_t$. Note that solving (19) involves maximizing a convex norm function, which can be hard in general (Xie 2020, Proposition 6). Nevertheless, it becomes tractable when the dual norm entails an inf-norm, which can be represented as a component-wise maximum absolute value; see Appendix EC.5.3. The resulting reformulation of $Q_t$ involves solving $2d_t$ problems in total, each of which perturbs the constraints of the scenario approximation problem by a unit vector with a magnitude proportional to $\vartheta$. When $T=2$, this is consistent with Xie (2020, Theorem 3). We will apply this result to a dynamic portfolio selection problem in Section 5. The next result provides an equivalent reformulation of (15) for $p=1$, whose proof can be found in Appendix EC.5.4.

COROLLARY 3. *Suppose $\Xi_t = (\mathbb{R}^{d_t}, \|\cdot\|)$, $t = 2,\dots,T$, and $p = 1$. Set $\mathcal{Y}_{T+1} \equiv \{0\}$, and for $t \in [T]$, define recursively*

$$\mathcal{S}_t := \left\{ y_t \in \mathbb{R}^{d_t} : \exists \boldsymbol{y}_{t+1} \in \mathcal{Y}_{t+1} \text{ s.t. } A_t^\top y_t + B_{t+1}^\top \mathbb{E}_{\widehat{\mathbb{P}}_{t+1}}[\boldsymbol{y}_{t+1}(\widehat{\boldsymbol{\xi}}_{t+1})] \le c_t \right\},$$

*where $\mathcal{Y}_t$ is the space of functions from $\widehat{\Xi}_t$ to $\mathcal{S}_t$. Set $\widehat{Q}_{T+1}(x_t, \widehat{\boldsymbol{\xi}}_{[t+1]}) \equiv 0$ and*

$$\widehat{Q}_t(x_{t-1}, \widehat{\boldsymbol{\xi}}_t) = \inf_{x_t \in \mathcal{X}_t(x_{t-1}, \widehat{\boldsymbol{\xi}}_t)} \left\{ c_t^\top x_t + \mathbb{E}_{\widehat{\mathbb{P}}_{t+1}} \left[ \widehat{Q}_{t+1}(x_t, \widehat{\boldsymbol{\xi}}_{t+1}) \right] \right\}, \quad t = 2,\dots,T.$$

*Then*

$$\begin{aligned}
Q_t(x_{t-1}, \widehat{\boldsymbol{b}}_t) &= \widehat{Q}_t(x_{t-1}, \widehat{\boldsymbol{b}}_t) + \infty \cdot \mathbf{1}\left\{ \lambda < \max_{s=t,\dots,T} \max_{y \in \mathcal{S}_s} \|y\|_* \right\}, && \text{if } \boldsymbol{\xi}_t = \boldsymbol{b}_t,\\
Q_t(x_{t-1}, \widehat{\boldsymbol{B}}_t) &= \widehat{Q}_t(x_{t-1}, \widehat{\boldsymbol{B}}_t) + \infty \cdot \mathbf{1}\left\{ \lambda < \max_{s=t,\dots,T} \max_{y \in \mathcal{S}_s} \|x_t y^\top\|_* \right\}, && \text{if } \boldsymbol{\xi}_t = \boldsymbol{B}_t,
\end{aligned}$$

*and the optimal value of (1) equals*

$$\begin{aligned}
&\inf_{x_1 \in \mathcal{X}_1} \left\{ c_1^\top x_1 + \mathbb{E}_{\widehat{\mathbb{P}}_2}[\widehat{Q}_2(x_1, \widehat{\boldsymbol{b}}_2)] \right\} + \vartheta \cdot \max_{t=2,\dots,T} \max_{y \in \mathcal{S}_t} \|y\|_*, && \text{if } \boldsymbol{\xi}_t = \boldsymbol{b}_t,\\
&\inf_{x_1 \in \mathcal{X}_1} \left\{ c_1^\top x_1 + \mathbb{E}_{\widehat{\mathbb{P}}_2}[\widehat{Q}_2(x_1, \widehat{\boldsymbol{B}}_2)] \right\} + \vartheta \cdot \max_{t=2,\dots,T} \max_{y \in \mathcal{S}_t} \|x_t y^\top\|_*, && \text{if } \boldsymbol{\xi}_t = \boldsymbol{B}_t.
\end{aligned} \tag{20}$$

Note that $\widehat{Q}_t$ is the cost-to-go function for the (non-robust) scenario approximation problem. The first term of (20) is the optimal value of the scenario approximation problem, while the second term of (20) is a linear function of $\vartheta$, whose value is independent of the policy $\boldsymbol{x}$. Consequently, (20) share the same optimal policy values as the scenario approximation problem on each sample path from the nominal scenario tree. We would like to emphasize that the optimal robust policy for (10), as in (16), is defined for all sample paths in $\Xi_{[T]}$, whereas the optimal policy for the scenario approximation problem is only defined for the sample paths within the scenario tree. In this respect, the robust formulation induces a safe way to extend the optimal solution to the scenario approximation problem across the entire sample space and justifies the heuristic policy in the literature (Shapiro et al. 2012, Keutchayan et al. 2017, Ding et al. 2019, Zhang and Sun 2022).

Corollary 3 generalizes the results for two-stage Wasserstein DRO (Hanasusanto and Kuhn 2018, Duque et al. 2022). Below we give an example in the context of the multistage newsvendor problem. A similar observation was made in Esfahani and Kuhn (2018) for the static newsvendor problem.

EXAMPLE 5 (NEWSVENDOR). Consider a multistage distributionally robust newsvendor model. Let $\boldsymbol{x}_t$ be the inventory level after having ordered in stage $t$ but before the demand $\boldsymbol{\xi}_t$ in that stage is realized. Let $c_t$, $c_t^b$ and $c_t^h$ be the ordering, back-order penalty and holding costs per unit in stage $t$, respectively. The multistage newsvendor is given by

$$\begin{aligned}
\inf_{\boldsymbol{x}_t, t \in [T]} \ \sup_{\mathbb{P} \in \mathfrak{M}} \ & \mathbb{E}_{\mathbb{P}} \left[ \sum_{t \in [T]} c_t(\boldsymbol{x}_t - \boldsymbol{x}_{t-1} + \boldsymbol{\xi}_{t-1}) + c_t^h(\boldsymbol{x}_{t-1} - \boldsymbol{\xi}_{t-1})_+ + c_t^b(\boldsymbol{\xi}_{t-1} - \boldsymbol{x}_{t-1})_+ \right]\\
\text{s.t.} \ & \boldsymbol{x}_t \ge \boldsymbol{x}_{t-1} - \boldsymbol{\xi}_{t-1},\\
& \boldsymbol{x}_t \ge 0.
\end{aligned}$$

With additional auxiliary variables $z_t$, we can rewrite it as a multistage linear program with right-hand uncertainty

$$
\begin{aligned}
\inf_{\boldsymbol{x}_t, z_t, t\in[T]} \quad & \sup_{\mathbb{P}\in\mathfrak{M}} \mathbb{E}_{\mathbb{P}}\left[c_1\boldsymbol{x}_1 + z_2 + \cdots + z_T\right] \\
\text{s.t.} \quad & -\boldsymbol{x}_t + \boldsymbol{x}_{t-1} \le \boldsymbol{\xi}_{t-1}, \\
& c_t\boldsymbol{x}_t - z_t + (c_t^h - c_t)\boldsymbol{x}_{t-1} \le (c_t^h - c_t)\boldsymbol{\xi}_{t-1}, \\
& c_t\boldsymbol{x}_t - z_t - (c_t + c_t^b)\boldsymbol{x}_{t-1} \le -(c_t + c_t^b)\boldsymbol{\xi}_{t-1}, \\
& \boldsymbol{x}_t \ge 0.
\end{aligned}
$$

Then Corollary 3 indicates that the Wasserstein robust solution and the non-robust solution coincide on sample paths within the nominal scenario tree. ♣

# 5. Application in Dynamic Portfolio Selection

## 5.1. Problem Formulation

We consider a portfolio selection problem of an investor who seeks to minimize the dis-utility of the terminal wealth. Given some initial wealth $W_1$, she invests in $n$ assets. The monetary value of all $n$ investments are represented using a vector, $x_t \in \mathbb{R}^n$. The return rate at time period $t$ is modeled by a random variable $\boldsymbol{\xi}_t \in \mathbb{R}_t$. At each stage before the terminal, she may re-balance her wealth by taking long-only positions across the $n$ investments. Suppose the investor's terminal dis-utility function is given by $U(W_T) := \max(-a_0 - r_0 W_T, -a_1 - r_1 W_T)$, where $a_0, a_1, r_0, r_1$ are used to encode the investor's preferences. We can write the robust portfolio selection problem as

$$
\begin{aligned}
\inf_{\boldsymbol{x}_1,\ldots,\boldsymbol{x}_{T-1}\ge 0, \boldsymbol{x}_T} \quad & \max_{\mathbb{P}\in\mathfrak{M}} \mathbb{E}_{\mathbb{P}}[U(\boldsymbol{x}_T)] \\
\text{s.t.} \quad & \mathbf{1}^\top \boldsymbol{x}_1 = W_1, \\
& \mathbf{1}^\top \boldsymbol{x}_t = \boldsymbol{\xi}_t^\top \boldsymbol{x}_{t-1}, \quad t = 2, \ldots, T-1, \\
& \boldsymbol{x}_T = \boldsymbol{\xi}_T^\top \boldsymbol{x}_{T-1}.
\end{aligned}
\tag{21}
$$

Using Corollary 2, we obtain the following reformulation, whose proof is given in Appendix EC.6.1.

Corollary 4. *Using the setup in Corollary 2, the dynamic programming reformulation of (21) is given by*

$$
\begin{aligned}
Q_T(x_{T-1}, \widehat{\boldsymbol{\xi}}_{[T]}) &:= \max_{\substack{z\in\{0,1\}\\ \delta\in\{1,-1\}}} \max\Big( -a_0 - r_0 \widehat{\boldsymbol{\xi}}_T^\top x_{T-1} + \vartheta(1-z)\delta\|x_{T-1}\|,\ -a_1 - r_1 \widehat{\boldsymbol{\xi}}_T^\top x_{T-1} + \vartheta z\delta\|x_{T-1}\| \Big), \\
Q_t(x_{t-1}, \widehat{\boldsymbol{\xi}}_{[t]}) &:= \max_{\delta\in\{1,-1\}} \inf_{\substack{\mathbf{1}^\top x_t = \widehat{\boldsymbol{\xi}}_t^\top x_{t-1} + \vartheta\delta\|x_{t-1}\|,\\ x_t\ge 0}} \mathbb{E}_{\widehat{\mathbb{P}}_{t+1|\widehat{\boldsymbol{\xi}}_{[t]}}} \left[ Q_{t+1}(x_t, \widehat{\boldsymbol{\xi}}_{[t+1]}) \right],\ t = 2, \ldots, T-1, \\
Q_1 &:= \inf_{\substack{\mathbf{1}^\top x_1 = W_1\\ x_1 \ge 0}} \mathbb{E}_{\widehat{\mathbb{P}}_{\widehat{\xi}_2}} \left[ Q_2(x_1, \widehat{\boldsymbol{\xi}}_2) \right].
\end{aligned}
$$

## 5.2. Experiment Setup

Let $a_0 = 0$, $a_1 = \frac{1}{2}W_1$, $r_0 = 1$, $r_1 = \frac{1}{2}$, which corresponds to the dis-utility

$$
U(W_T) := \begin{cases} -W_T, & W_T \le W_1, \\ -W_1 - \frac{1}{2}(W_T - W_1), & W_T > W_1. \end{cases}
$$

Suppose the investor starts with an initial wealth of $W_1 = 10000$ units, and can invest the wealth across the following 5 assets, iShares MSCI Emerging Markets ETF (EEM), iShares 20+ Year Treasury Bond

ETF (TLT), Schwab US TIPS ETF (SCHP), SPDR S&P Oil & Gas Equipment & Services ETF (XES), and ProShares UltraShort Financials ETF (SKF). We simulate the monthly asset returns using a log-normal distribution with mean and covariance estimated using adjusted closing prices from January 1, 2018 to June 30, 2021; see Appendix EC.6.2 for estimation results. We assume stage-wise independence, for the convenience of the out-of-sample test.

We describe our out-of-sample testing procedure as follows. The training dataset is a $T$-stage scenario tree, where at stage $t = 2, \ldots, T$ there are $\widehat{N}_t$ independent scenarios. The testing dataset is another $T$-stage scenario tree, independent from the training tree, where at stage $t = 2, \ldots, T$ there are $N_t$ independent scenarios. The reformulated problem in Corollary 4 can be viewed as a regularized SAA problem, so we solve it using a modified SDDP algorithm; see Algorithm 1 in Appendix EC.6.3 for details. Since the SDDP algorithm does not provide a policy but only the first-stage decision, we need to resolve the remaining subproblems to obtain subsequent robust decisions.

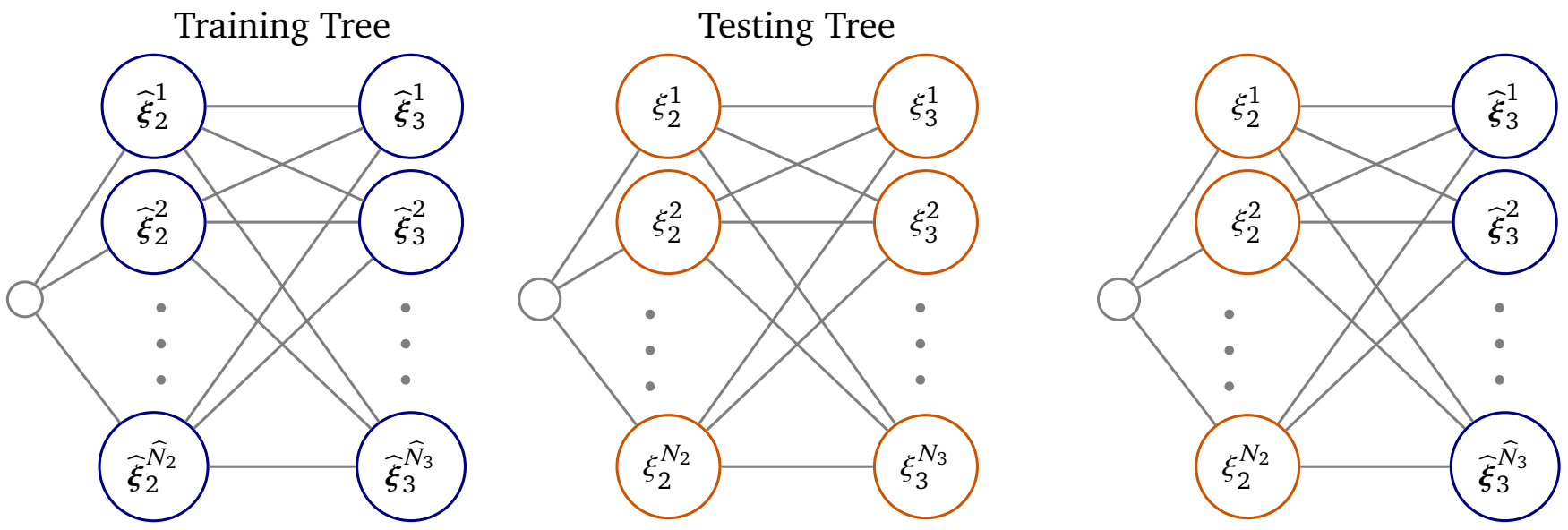

**Figure 5** Illustration of the out-of-sample testing procedure for a 3-stage problem. In the last tree, the second-stage decision is computed using the training scenarios in the third stage after observing a testing scenario in the second stage.

To describe our out-of-sample testing procedure, consider a 3-stage problem as an example, illustrated in Figure 5. First, we use the entire training tree $\{(\widehat{\boldsymbol{\xi}}_t^1, \ldots, \widehat{\boldsymbol{\xi}}_t^{N_t})\}_{t=2,\ldots,T}$ to obtain a first-stage robust decision $x_1^{\text{rob}}$ and evaluate the first-stage cost. Next, at stage 2, we observe a sample $\xi_2^{i_2}$ from the second-stage scenarios in the testing tree $(\xi_2^1, \ldots, \xi_2^{N_2})$, and use all scenarios in the remaining stages (i.e., stage 3 in this setting) from the training tree $(\widehat{\boldsymbol{\xi}}_3^1, \ldots, \widehat{\boldsymbol{\xi}}_3^{\widehat{N}_3})$ to solve for the second-stage robust decision $x_2^{\text{rob}}$ and evaluate its second-stage cost. At stage 3, we draw a sample $\xi_3^{i_3}$ from the third-stage scenarios in the testing tree $(\xi_3^1, \ldots, \xi_2^{N_3})$, and use all scenarios in the remaining stages from the training tree (i.e., none in this setting) to solve for the third-stage robust decision $x_3^{\text{rob}}$ and evaluate its third-stage cost. After going through all stages, we sum up the per-stage costs in all stages, which gives a realization of the out-of-sample cost d with the testing sample path $(\xi_2^{i_2}, \xi_3^{i_3})$. To get an estimate of the expected out-of-sample cost, we follow the procedure above to sample $M$ testing paths and average over them. We refer to Algorithm 2 in Appendix EC.6.3 for a pseudo-code for general $T$-stage problems.

## 5.3. Numerical Results

**5.3.1. Comparison With Sample Average Approximation** In our first experiment, we set $T \in \{3, 4, 5\}$, $\widehat{N}_2 = \cdots = \widehat{N}_T \in \{2, 5, 10\}$, $N_2 = \cdots = N_T = 30$, $M = 25$. We are interested in how the out-of-sample performance depends on different parameter values by varying $\vartheta \in \{0.1, 0.3, 0.5\}$. The benchmark is chosen as the sample average approximation counterpart of our robust formulation. We repeat the above out-of-sample testing procedure 30 times, each of which has an independent instance, and we report the resulting boxplots in Figure 6.

The left column of Figure 6 shows the out-of-sample expected utility of the optimal SAA solution and the optimal robust solution for different choices of $T$ and $\widehat{N}_t$. We have the following observations:

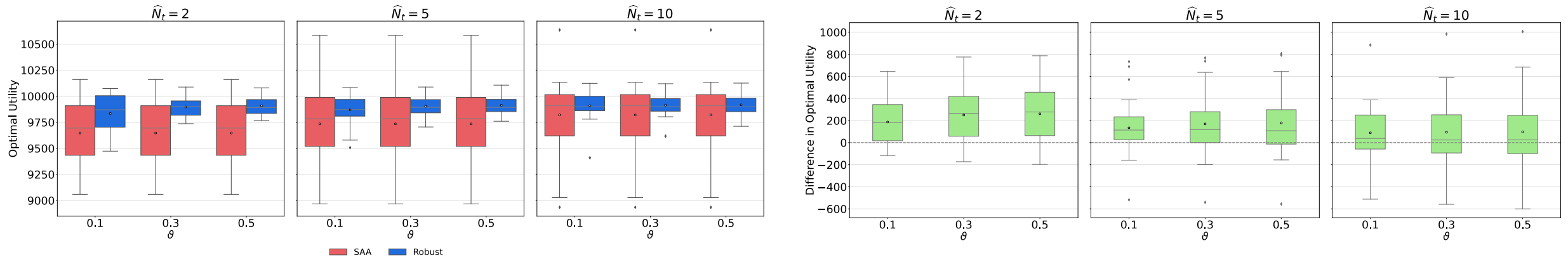

(a) $T = 3$

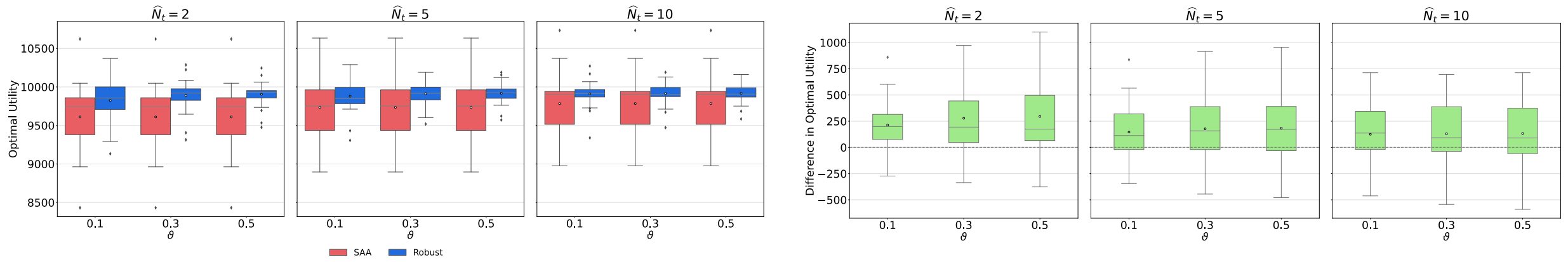

(b) $T = 4$

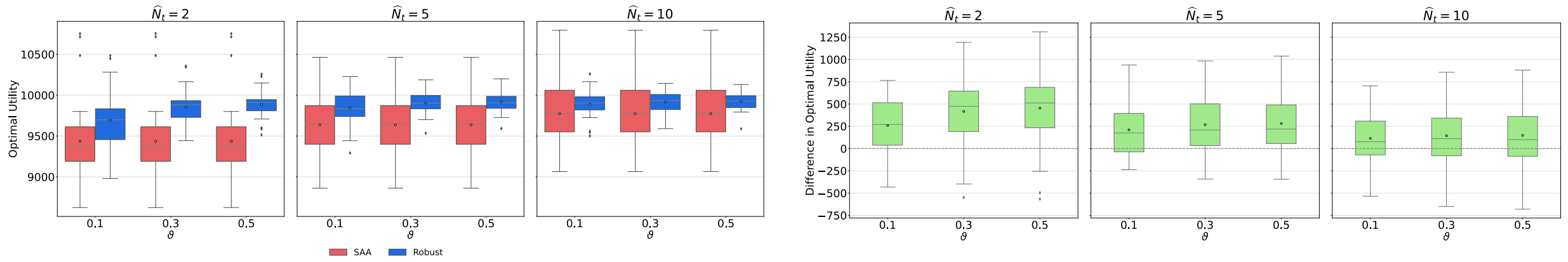

(c) $T = 5$

**Figure 6** Out-of-sample expected optimal utility yielding from robust and SAA formulations (left column) and their differences (right column)

(I) As $\widehat{N}_t$ increases, both the SAA solution and robust solutions achieve a higher out-of-sample expected utility. This makes sense because a larger sample size yields a more faithful representation of the underlying stochastic process.

(II) The average out-of-sample performance (as indicated by the circles) of the robust solution is consistently better than that of the SAA solution, and the variability of the out-of-sample performance (as indicated by the box length) of the robust solutions is consistently smaller than that of the SAA solution. This shows the practical importance of having a robust formulation to achieve better out-of-sample performance.

(III) When $\widehat{N}_t = 2, 5$, a large radius $\vartheta = 0.5$ has the best out-of-sample performance; whereas when $\widehat{N}_t = 10$, a large radius does not have clear advantage anymore.

These observations are consistent with our intuition and hold for all choices of $T$.

To further investigate the impact of the sample size and the radius on the out-of-sample performance, on the right column of Figure 6, we plot the instance-wise difference between SAA and robust solutions. A positive value means the robust solution performs better than the SAA solution out-of-sample for a particular instance. We have the following observations:

(I) The performance of the robust solution has a clear advantage over the SAA solution when the sample size $\widehat{N}_t$ is small, and the advantage diminishes as the sample size becomes larger.

(II) The best radius decreases as the sample size $\widehat{N}_t$ increases and increases as $T$ increases. This makes sense because the distributional uncertainty reduces when more sample paths are observed and amplifies when there are more stages.

These observations are also consistent with our intuition and validate the robust approach. In summary, the numerical results demonstrate the clear advantage of our robust formulation as compared to the sample average approximation.

**5.3.2. Comparison With Average Value at Risk** In this experiment, we compare the out-of-sample performance between the DRO model with nested distance and another commonly used approach, Average Value at Risk (AVaR) (e.g., (Lan and Shapiro 2024)). The AVaR model can be solved by the SDDP algorithm similarly. Set $T \in \{3,5,7,9\}$, $\widehat{N}_2 = \cdots = \widehat{N}_T \in \{2,5,10\}$, $N_2 = \cdots = N_T = 20$, $M = 20$. For the nested distance, set $\vartheta \in \{0.001, 0.003, 0.01, 0.03, 0.1, 0.3, 0.5\}$; for AVaR, set $\alpha \in \{0.05, 0.1, 0.2, 0.5\}$. These hyper-parameters are tuned according to the best out-of-sample performance. We repeat the procedure 30 times. Figure 7 shows the difference in the out-of-sample expected utility between the robust model and the SAA model, for two robust approaches. We have the following observations:

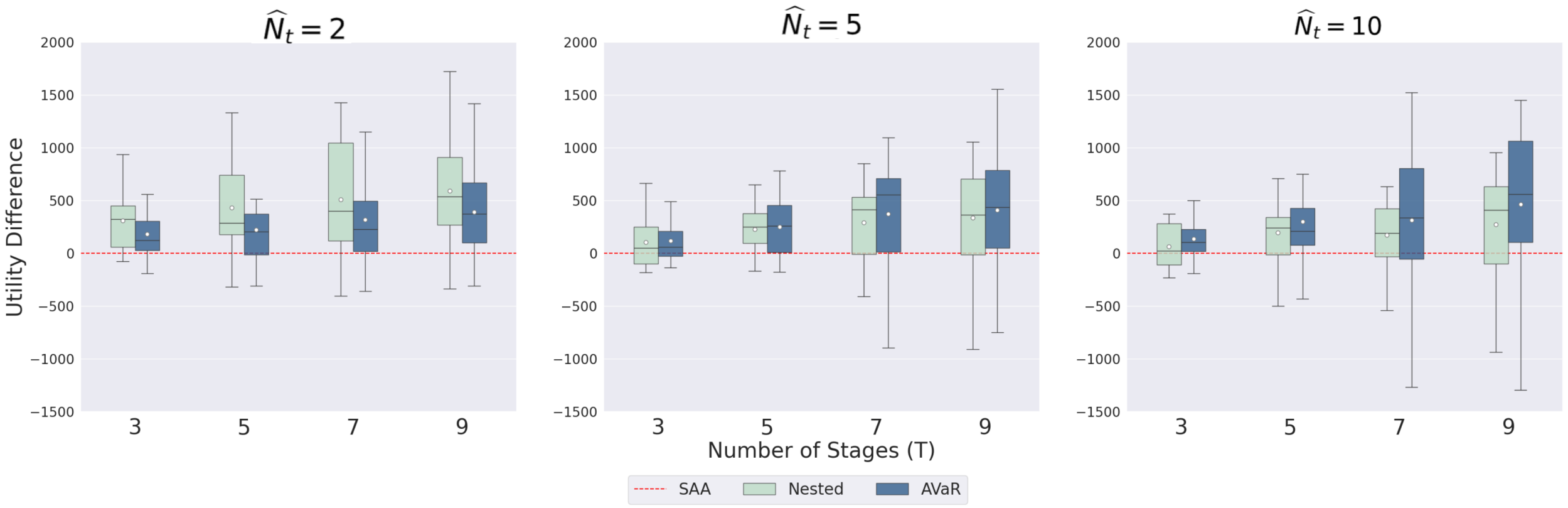

**Figure 7** Out-of-sample expected optimal utility difference yielding from nested and AVaR formulations

(I) Both approaches consistently outperform the non-robust approach in terms of the average out-of-sample performance.

(II) When the training sample size is small, the model with nested distance has higher average and median out-of-sample expected utilities. On the other hand, the model with AVaR has better out-of-sample performance when the training sample size becomes larger.

In summary, the model with nested distance outperforms the model with AVaR when the sample size is small. For larger sample sizes, although the model with nested distance has a lower average expected utility, it offers better protection against worst-case scenarios.

## 6. Concluding Remarks

In this paper, we develop reformulations for distributionally robust optimization with nested distance. These reformulations unveil equivalence between static and dynamic formulations of multistage distributionally robust problems and can be viewed as sample average approximation with norm regularization. For future work, it is interesting to study multistage problems with general convex objectives and constraints and extensions to other sequential decision-making frameworks such as stochastic control, as well as improve the finite-sample performance guarantees of the robust solution.

# Supplementary Material

## EC.1. Proof for Proposition 1

Using (Acciaio and Hou 2022, Theorem 2.19(i)), when $q=0$, there exists constants $c, C, K>0$ such that for all $\delta \geq 0$ and $N \in \mathbb{N}$, it holds that

$$\mathbb{P}\left\{\mathsf{D}_1(\mathbb{P}^*, \widehat{\mathbb{P}}^N) \geq \delta + KN^{-\frac{1}{D(d)T}}\right\} \leq Ce^{-cN\delta^2}.$$

Let $\beta$ equals the right-hand side, solving for $\vartheta_N(\beta) = \delta + KN^{-\frac{1}{D(d)T}}$ yields

$$\vartheta_N(\beta) = (\frac{\log(C/\beta)}{cN})^{\frac{1}{2}} + KN^{-\frac{1}{D(d)T}}.$$

Using the same cited result, when $q>0$, there exists constants $c, C, K>0$ such that for all $\delta \geq 0$ and $N \in \mathbb{N}$ such that $\sqrt{N}\delta \geq 1$, it holds that

$$\mathbb{P}\left\{\mathsf{D}_1(\mathbb{P}^*, \widehat{\mathbb{P}}^N) \geq \vartheta + KN^{-\frac{1}{D(d)T}}\right\} \leq Ce^{-c\sqrt{N}^{\frac{1}{q+1}}\delta^{\frac{1}{q+1}}}.$$

Let $\beta$ equals the right hand side, solving for $\vartheta_N(\beta)$ yields

$$\vartheta_N(\beta) = (\frac{\log(C/\beta)}{c})^{q+1} N^{-1/2} + KN^{-\frac{1}{D(d)T}}.$$

Note that $\sqrt{N}\delta \geq 1$ whenever $\beta \leq \frac{C}{e^c}$. Since $D(d) \geq 2$, we have $D(d)T \geq 2$. Thereby, with probability no less than $1-\beta$, the true distribution falls in the ambiguity set $\mathfrak{M}$ of radius $\vartheta_N(\beta)$. Consequently, with probability at least $1-\beta$, for every feasible policy $\boldsymbol{x} \in \mathcal{X}$, we have

$$\mathbb{E}_{\mathbb{P}^*}\left[\sum_{t\in[T]} \boldsymbol{c}_t^\top \boldsymbol{x}_t(\boldsymbol{\xi}_{[t]})\right] \leq \sup_{\mathbb{P}\in\mathfrak{M}} \mathbb{E}_{\mathbb{P}}\left[\sum_{t\in[T]} \boldsymbol{c}_t^\top \boldsymbol{x}_t(\boldsymbol{\xi}_{[t]})\right].$$

Taking the infimum over $\boldsymbol{x}$ yields the result. □

## EC.2. Proofs for Section 3.1

We consider a relaxation

$$\sup_{\mathbb{P}\in\mathfrak{M}^{\mathsf{C}}} \mathbb{E}_{\mathbb{P}}\left[\sum_{t\in[T]} \boldsymbol{c}_t^\top \boldsymbol{x}_t(\boldsymbol{\xi}_{[t]})\right], \tag{EC.1}$$

where the nested distance uncertainty set (3) is relaxed to the causal transport distance uncertainty set

$$\mathfrak{M}^{\mathsf{C}} = \left\{\mathbb{P} \in \mathcal{P}(\Xi) : \mathsf{C}_p(\widehat{\mathbb{P}}, \mathbb{P}) \leq \vartheta\right\}.$$

We define

$$\mathfrak{M}^{\mathsf{C}}_{[t]} := \left\{\mathbb{P}_{[t]} \in \mathcal{P}(\Xi_{[t]}) : \mathsf{C}_p(\widehat{\mathbb{P}}_{[t]}, \mathbb{P}_{[t]}) \leq \vartheta\right\}$$

and denote by $\Gamma^{\mathsf{C}}_{[t]}$ the set of all causal transport plans from $\widehat{\mathbb{P}}_{[t]}$ to the distributions in $\mathfrak{M}^{\mathsf{C}}_{[t]}$.

### EC.2.1. Proof of Theorem 1

We have the following result for the relaxed problem (EC.1).

Proposition EC.1. *Let $p=\infty$. Then for any feasible policy $(\boldsymbol{x}_1,\dots,\boldsymbol{x}_T)$, it holds that*

$$\begin{aligned}&\sup_{\mathbb{P}\in\mathfrak{M}^{\mathrm{C}}}\mathbb{E}_{\mathbb{P}}\left[\sum_{s\in[T]}\boldsymbol{c}_s^\top\boldsymbol{x}_s(\boldsymbol{\xi}_{[s]})\right]\\&=\sup_{\gamma_{[t]}\in\Gamma^{\mathrm{C}}_{[t]}}\mathbb{E}_{\gamma_{[t]}}\left[\sum_{s\in[t]}\boldsymbol{c}_s^\top\boldsymbol{x}_s(\boldsymbol{\xi}_{[s]})+\mathbb{E}_{\widehat{\mathbb{P}}_{t+1|\widehat{\xi}_{[t]}}}\left[\sup_{\xi_{t+1}\in\Xi_{t+1}:\|\xi_{t+1}-\widehat{\xi}_{t+1}\|\le\vartheta}V_{t+1}^{\widehat{\xi}_{[t+1]}}(\boldsymbol{\xi}_{[t]},\xi_{t+1})\right]\right].\end{aligned}$$

*Proof of Proposition EC.1.* We prove this by induction. By definition, the case of $t=T$ holds trivially from the definition of the causal transport distance. Suppose we have proved for the case of $t$, where $t=2,\dots,T$, now we prove the result also holds for $t-1$. Set

$$V:=\sup_{\mathbb{P}\in\mathfrak{M}^{\mathrm{C}}}\mathbb{E}_{\mathbb{P}}\left[\sum_{s\in[T]}\boldsymbol{c}_s^\top\boldsymbol{x}_s(\boldsymbol{\xi}_{[s]})\right]$$

and

$$\bar{V}_{t+1}^{\widehat{\xi}_{[t]}}(\boldsymbol{\xi}_{[t]}):=\mathbb{E}_{\widehat{\xi}_{t+1}\sim\widehat{\mathbb{P}}_{t+1|\widehat{\xi}_{[t]}}}\left[\sup_{\xi_{t+1}\in\Xi_{t+1}:\|\xi_{t+1}-\widehat{\xi}_{t+1}\|\le\vartheta}V_{t+1}^{\widehat{\xi}_{[t+1]}}(\boldsymbol{\xi}_{[t]},\xi_{t+1})\right].$$

Then the induction hypothesis implies that

$$V=\sup_{\gamma_{[t]}\in\Gamma^{\mathrm{C}}_{[t]}}\mathbb{E}_{\gamma_{[t]}}\left[\sum_{s\in[t]}\boldsymbol{c}_s^\top\boldsymbol{x}_s(\boldsymbol{\xi}_{[s]})+\bar{V}_{t+1}^{\widehat{\xi}_{[t]}}(\boldsymbol{\xi}_{[t]})\right].$$

For any causal transport plan $\gamma_{[t]}$ from $\widehat{\mathbb{P}}_{[t]}$ to $\mathbb{P}_{[t]}$, by the tower property of conditional expectation, it holds that

$$\begin{aligned}&\mathbb{E}_{\gamma_{[t]}}\left[\sum_{s\in[t]}\boldsymbol{c}_s^\top\boldsymbol{x}_s(\boldsymbol{\xi}_{[s]})+\bar{V}_{t+1}^{\widehat{\xi}_{[t]}}(\boldsymbol{\xi}_{[t]})\right]\\&=\mathbb{E}_{\gamma_{[t-1]}}\left[\sum_{s\in[t-1]}\boldsymbol{c}_s^\top\boldsymbol{x}_s(\boldsymbol{\xi}_{[s]})\right]+\mathbb{E}_{\gamma_{[t]}}\left[\boldsymbol{c}_t^\top\boldsymbol{x}_t(\boldsymbol{\xi}_{[t]})+\bar{V}_{t+1}^{\widehat{\xi}_{[t]}}(\boldsymbol{\xi}_{[t]})\right]\\&=\mathbb{E}_{\gamma_{[t-1]}}\left[\sum_{s\in[t-1]}\boldsymbol{c}_s^\top\boldsymbol{x}_s(\boldsymbol{\xi}_{[s]})\right]+\mathbb{E}_{\gamma_{[t-1]}}\left[\mathbb{E}_{\gamma_{\widehat{\xi}_t|(\widehat{\xi}_{[t-1]},\xi_{[t-1]})}}\left[\mathbb{E}_{\gamma_{\xi_t|(\widehat{\xi}_{[t]},\xi_{[t-1]})}}\left[\boldsymbol{c}_t^\top\boldsymbol{x}_t(\boldsymbol{\xi}_{[t]})+\bar{V}_{t+1}^{\widehat{\xi}_{[t]}}(\boldsymbol{\xi}_{[t]})\right]\right]\right]\\&=\mathbb{E}_{\gamma_{[t-1]}}\left[\sum_{s\in[t-1]}\boldsymbol{c}_s^\top\boldsymbol{x}_s(\boldsymbol{\xi}_{[s]})\right]+\mathbb{E}_{\gamma_{[t-1]}}\left[\mathbb{E}_{\widehat{\mathbb{P}}_{t|\widehat{\xi}_{[t-1]}}}\left[\mathbb{E}_{\gamma_{\xi_t|(\widehat{\xi}_{[t]},\xi_{[t-1]})}}\left[\boldsymbol{c}_t^\top\boldsymbol{x}_t(\boldsymbol{\xi}_{[t]})+\bar{V}_{t+1}^{\widehat{\xi}_{[t]}}(\boldsymbol{\xi}_{[t]})\right]\right]\right],\end{aligned}\tag{EC.2}$$

where, in the last equations above, we have used the non-anticipativity of the causal transport plan. Note that given the distribution $\widehat{\mathbb{P}}_{[t]}$, the causal transport plan $\gamma_{[t]}$ is determined completely by $\gamma_{[t-1]}$ and $\gamma_{\xi_t|(\widehat{\xi}_{[t]},\xi_{[t-1]})}$, and the constraint

$$\gamma_{[t]}\text{-ess sup}\max_{s\in[t]}\|\widehat{\boldsymbol{\xi}}_s-\boldsymbol{\xi}_s\|\le\vartheta$$

is equivalent to

$$\gamma_{[t-1]}\text{-ess sup}\max_{s\in[t-1]}\|\widehat{\boldsymbol{\xi}}_s-\boldsymbol{\xi}_s\|\le\vartheta$$

and

$$\gamma_{\xi_t|(\widehat{\xi}_{[t]},\xi_{[t-1]})}\text{-ess sup}\,\|\widehat{\boldsymbol{\xi}}_t-\boldsymbol{\xi}_t\|\leq\vartheta,\quad\forall\widehat{\boldsymbol{\xi}}_{[t]}\in\operatorname{supp}\widehat{\mathbb{P}}_{[t]},\boldsymbol{\xi}_{[t-1]}\in\operatorname{supp}\mathbb{P}_{[t]}.$$

Thereby maximizing over $\gamma_{[t]}\in\Gamma^{\mathsf{C}}_{[t]}$ is equivalent to maximizing over $\gamma_{[t-1]}\in\Gamma^{\mathsf{C}}_{[t-1]}$ and $\gamma_{\xi_t|(\widehat{\xi}_{[t]},\xi_{[t-1]})}\in\mathfrak{M}_t(\widehat{\boldsymbol{\xi}}_t):=\{\mathbb{P}_t\in\mathcal{P}(\Xi_t):\mathbb{P}_t\text{-ess sup}\|\xi_t-\widehat{\boldsymbol{\xi}}_t\|\leq\vartheta\}$. Observe that

$$\begin{aligned}
&\sup_{\gamma_{\xi_t|(\widehat{\xi}_{[t]},\xi_{[t-1]})}\in\mathfrak{M}_t(\widehat{\boldsymbol{\xi}}_t)}\mathbb{E}_{\gamma_{\xi_t|(\widehat{\xi}_{[t]},\xi_{[t-1]})}}\left[\boldsymbol{c}_t^\top\boldsymbol{x}_t(\boldsymbol{\xi}_{[t]})+\bar{V}_{t+1}^{\widehat{\boldsymbol{\xi}}_{[t]}}(\boldsymbol{\xi}_{[t]})\right]\\
=&\sup_{\xi_t\in\Xi_t:\|\xi_t-\widehat{\boldsymbol{\xi}}_t\|\leq\vartheta}\left\{c_t^\top\boldsymbol{x}_t(\boldsymbol{\xi}_{[t-1]},\xi_t)+\bar{V}_{t+1}^{\widehat{\boldsymbol{\xi}}_{[t]}}(\boldsymbol{\xi}_{[t-1]},\xi_t)\right\}\\
=&\sup_{\xi_t\in\Xi_t:\|\xi_t-\widehat{\boldsymbol{\xi}}_t\|\leq\vartheta}V_t^{\widehat{\boldsymbol{\xi}}_{[t]}}(\boldsymbol{\xi}_{[t-1]},\xi_t).
\end{aligned}$$

Hence, together with (EC.2), it follows that

$$V=\sup_{\gamma_{[t-1]}\in\Gamma^{\mathsf{C}}_{[t-1]}}\mathbb{E}_{\gamma_{[t-1]}}\left[\sum_{s\in[t-1]}\boldsymbol{c}_s^\top\boldsymbol{x}_s(\boldsymbol{\xi}_{[s]})+\mathbb{E}_{\widehat{\mathbb{P}}_{t|\widehat{\xi}_{[t-1]}}}\left[\sup_{\xi_t\in\Xi_t:\|\xi_t-\widehat{\boldsymbol{\xi}}_t\|\leq\vartheta}V_t^{\widehat{\boldsymbol{\xi}}_{[t]}}(\boldsymbol{\xi}_{[t-1]},\xi_t)\right]\right],$$

which completes the induction. □

*Proof of Theorem 1.* We prove the theorem by showing that for every feasible policy $(x_1,\boldsymbol{x}_2,\cdots,\boldsymbol{x}_T)$, it holds that

$$\sup_{\mathbb{P}\in\mathfrak{M}}\mathbb{E}_{\mathbb{P}}\left[\sum_{t\in[T]}\boldsymbol{c}_t^\top\boldsymbol{x}_t(\boldsymbol{\xi}_{[t]})\right]=\sup_{\mathbb{P}\in\mathfrak{M}^{\mathsf{C}}}\mathbb{E}_{\mathbb{P}}\left[\sum_{t\in[T]}\boldsymbol{c}_t^\top\boldsymbol{x}_t(\boldsymbol{\xi}_{[t]})\right].$$

Since $\mathfrak{M}\subset\mathfrak{M}^{\mathsf{C}}$, the left-hand side is less than or equal to the right-hand side. It remains to prove the other direction. Using Proposition EC.1, the right-hand side above equals

$$\begin{aligned}
&c_1^\top x_1+\mathbb{E}_{\widehat{\mathbb{P}}_{\widehat{\xi}_2}}\Bigg[\sup_{\xi_2\in\Xi_2:\|\widehat{\xi}_2-\xi_2\|\leq\vartheta}\Bigg\{c_2^\top x_2(\xi_2)+\mathbb{E}_{\widehat{\mathbb{P}}_{3|\widehat{\xi}_2}}\Bigg[\sup_{\xi_3\in\Xi_3:\|\widehat{\xi}_3-\xi_3\|\leq\vartheta}\Bigg\{c_3^\top x_3(\xi_{[3]})+\cdots+\\
&\qquad\mathbb{E}_{\widehat{\mathbb{P}}_{T|\widehat{\xi}_{[T-1]}}}\Bigg[\sup_{\xi_T\in\Xi_T:\|\widehat{\xi}_T-\xi_T\|\leq\vartheta}c_T^\top x_T(\xi_{[T]})\Bigg]\cdots\Bigg\}\Bigg]\Bigg\}\Bigg].
\end{aligned}$$

Denote $f_t(\xi_{[t]})=c_t^\top x_t(\xi_{[t]})$. Let $\mathbb{T}_1^{\epsilon_1}$ be the identity transport map. For $t=2,\ldots,T$, for any $\epsilon_t>0$ that is sufficiently small, there exists a transport map $\mathbb{T}_t^{\epsilon_t}$, where $\mathbb{T}_t^{\epsilon_t}:\widehat{\Xi}_{[t]}\to\Xi_t$, such that $\|\mathbb{T}_t^{\epsilon_t}(\widehat{\boldsymbol{\xi}}_{[t]})-\widehat{\boldsymbol{\xi}}_t\|\leq\vartheta$, and

$$\begin{aligned}
&f_t(\mathbb{T}_{[t]}^{\epsilon_t}(\widehat{\boldsymbol{\xi}}_{[t]}))+\mathbb{E}_{\widehat{\mathbb{P}}_{t+1|\widehat{\xi}_{[t]}}}\Bigg[\sup_{\xi_{t+1}\in\Xi_{t+1}:\|\widehat{\xi}_{t+1}-\xi_{t+1}\|\leq\vartheta}\Bigg\{f_{t+1}\Big(\mathbb{T}_{[t]}^{\epsilon_t}(\widehat{\boldsymbol{\xi}}_{[t]}),\xi_{t+1}\Big)+\cdots+\\
&\qquad\mathbb{E}_{\widehat{\mathbb{P}}_{T|\widehat{\xi}_{[T-1]}}}\Bigg[\sup_{\xi_T\in\Xi_T:\|\widehat{\xi}_T-\xi_T\|\leq\vartheta}f_T\Big(\mathbb{T}_{[t]}^{\epsilon_t}(\widehat{\boldsymbol{\xi}}_{[t]}),\xi_{[t+1,T]}\Big)\Bigg]\cdots\Bigg\}\Bigg]\\
\geq&\sup_{\xi_t\in\Xi_t:\|\widehat{\xi}_t-\xi_t\|\leq\vartheta}V_t^{\widehat{\boldsymbol{\xi}}_{[t]}}\Big(\mathbb{T}_{[t-1]}^{\epsilon_{[t-1]}}(\widehat{\boldsymbol{\xi}}_{[t-1]}),\xi_t\Big)-\epsilon_t.
\end{aligned}$$

Moreover, $\mathbb{T}_t^{\epsilon_t}$ can be chosen to be bijective. Indeed, if $\mathbb{T}_t^{\epsilon_t}(\widehat{\boldsymbol{\xi}}_{[t]})=\mathbb{T}_t^{\epsilon_t}(\widehat{\boldsymbol{\xi}}'_{[t]})$ for some $\widehat{\boldsymbol{\xi}}_{[t]}\neq\widehat{\boldsymbol{\xi}}'_{[t]}$ we can modify the mapping to ensure bijectivity. Since $V_t^{\widehat{\boldsymbol{\xi}}_{[t]}}(\mathbb{T}_{[t-1]}^{\epsilon_{[t-1]}}(\widehat{\boldsymbol{\xi}}_{[t-1]}),\cdot)$ is lower semicontinuous, there always exists another point $\xi_t$ such that: $\|\xi_t - \widehat{\boldsymbol{\xi}}_{[t]}\| \le \vartheta$; it is distinct from any value in $\{\mathbb{T}_t^{\epsilon_t}(\widehat{\xi}_{[t]}) : \widehat{\xi}_{[t]} \in \operatorname{supp} \Xi_{[t]}\}$; and setting $\mathbb{T}_t^{\epsilon_t}(\widehat{\boldsymbol{\xi}}_{[t]}) := \xi_t$ satisfies the above inequality, provided that $\epsilon_t$ is sufficiently small. We can apply this procedure to all coinciding values sequentially to redefine $\mathbb{T}_t^{\epsilon_t}(\cdot)$ as a bijective function. It follows that $\mathbb{T}^\epsilon = (\mathbb{T}_1^{\epsilon_1}, \ldots \mathbb{T}_T^{\epsilon_T})$ is a bi-causal transport map as $(\mathbb{T}^\epsilon)^{-1} := ((\mathbb{T}_1^\epsilon)^{-1}, \ldots, (\mathbb{T}_T^\epsilon)^{-1})$ is causal (Backhoff et al. 2017, Remark 3.4). Define

$$\mathbb{P}^\epsilon := \mathbb{T}^\epsilon_\# \widehat{\mathbb{P}}.$$

It follows that

$$\mathsf{D}_\infty(\widehat{\mathbb{P}}, \mathbb{P}^\epsilon) = \sup_{\widehat{\xi}_{[T]} \in \operatorname{supp} \widehat{\mathbb{P}}} \|\widehat{\boldsymbol{\xi}}_{[T]} - \mathbb{T}^\epsilon(\widehat{\boldsymbol{\xi}}_{[T]})\| \le \vartheta.$$

Hence $\mathbb{P}^\epsilon$ is a feasible distribution. In addition, for $t = 2, \ldots, T$, we have

$$\begin{aligned}
&\mathbb{E}_{\widehat{\mathbb{P}}_{[t-1]}}\left[\sum_{s\in[t-1]} \boldsymbol{c}_s^\top \boldsymbol{x}_s(\mathbb{T}_s^\epsilon(\widehat{\boldsymbol{\xi}}_{[s]})) + \mathbb{E}_{\widehat{\mathbb{P}}_{t|\widehat{\boldsymbol{\xi}}_{[t-1]}}}\left[\sup_{\xi_t\in\Xi_t:\|\widehat{\boldsymbol{\xi}}_t-\xi_t\|\le\vartheta} V_t^{\widehat{\boldsymbol{\xi}}_{[t]}}\left(\mathbb{T}_{[t-1]}^{\epsilon[t-1]}(\widehat{\boldsymbol{\xi}}_{[t-1]}), \xi_t\right)\right]\right] - \epsilon_t \\
&\le \mathbb{E}_{\widehat{\mathbb{P}}_{[t-1]}}\left[\sum_{s\in[t-1]} \boldsymbol{c}_s^\top \boldsymbol{x}_s(\mathbb{T}_s^\epsilon(\widehat{\boldsymbol{\xi}}_{[s]})) + \mathbb{E}_{\widehat{\mathbb{P}}_{t|\widehat{\boldsymbol{\xi}}_{[t-1]}}}\left[f_t(\mathbb{T}_{[t]}^{\epsilon_t}(\widehat{\boldsymbol{\xi}}_{[t]})) + \right.\right. \\
&\qquad \mathbb{E}_{\widehat{\mathbb{P}}_{t+1|\widehat{\boldsymbol{\xi}}_{[t]}}}\left[\sup_{\xi_{t+1}\in\Xi_{t+1}:\|\widehat{\boldsymbol{\xi}}_{t+1}-\xi_{t+1}\|\le\vartheta}\left\{f_{t+1}\left(\mathbb{T}_{[t]}^{\epsilon_t}(\widehat{\boldsymbol{\xi}}_{[t]}), \xi_{t+1}\right) + \cdots + \right.\right. \\
&\qquad \left.\left.\left.\left.\mathbb{E}_{\widehat{\mathbb{P}}_{T|\widehat{\boldsymbol{\xi}}_{[T-1]}}}\left[\sup_{\xi_T\in\Xi_T:\|\widehat{\boldsymbol{\xi}}_T-\xi_T\|\le\vartheta} f_T\left(\mathbb{T}_{[t]}^{\epsilon_t}(\widehat{\boldsymbol{\xi}}_{[t]}), \xi_{[t+1,T]}\right)\right]\cdots\right\}\right]\right]\right] \\
&= \mathbb{E}_{\widehat{\mathbb{P}}_{[t]}}\left[\sum_{s\in[t]} \boldsymbol{c}_s^\top \boldsymbol{x}_s(\mathbb{T}_s^\epsilon(\widehat{\boldsymbol{\xi}}_{[s]})) + \mathbb{E}_{\widehat{\mathbb{P}}_{t+1|\widehat{\boldsymbol{\xi}}_{[t]}}}\left[\sup_{\xi_{t+1}\in\Xi_{t+1}:\|\widehat{\boldsymbol{\xi}}_{t+1}-\xi_{t+1}\|\le\vartheta}\left\{f_{t+1}\left(\mathbb{T}_{[t]}^{\epsilon_t}(\widehat{\boldsymbol{\xi}}_{[t]}), \xi_{t+1}\right) + \cdots + \right.\right.\right. \\
&\qquad \left.\left.\left.\mathbb{E}_{\widehat{\mathbb{P}}_{T|\widehat{\boldsymbol{\xi}}_{[T-1]}}}\left[\sup_{\xi_T\in\Xi_T:\|\widehat{\boldsymbol{\xi}}_T-\xi_T\|\le\vartheta} f_T\left(\mathbb{T}_{[t]}^{\epsilon_t}(\widehat{\boldsymbol{\xi}}_{[t]}), \xi_{[t+1,T]}\right)\right]\cdots\right\}\right]\right] \\
&= \mathbb{E}_{\widehat{\mathbb{P}}_{[t]}}\left[\sum_{s\in[t]} \boldsymbol{c}_s^\top \boldsymbol{x}_s(\mathbb{T}_s^\epsilon(\widehat{\boldsymbol{\xi}}_{[s]})) + \mathbb{E}_{\widehat{\mathbb{P}}_{t+1|\widehat{\boldsymbol{\xi}}_{[t]}}}\left[\sup_{\xi_{t+1}\in\Xi_{t+1}:\|\widehat{\boldsymbol{\xi}}_{t+1}-\xi_{t+1}\|\le\vartheta} V_{t+1}^{\widehat{\boldsymbol{\xi}}_{[t+1]}}\left(\mathbb{T}_{[t]}^{\epsilon[t]}(\widehat{\boldsymbol{\xi}}_{[t]}), \xi_{t+1}\right)\right]\right].
\end{aligned}$$

Therefore by induction, we can show that

$$V_1^{\widehat{\xi}_1} - \sum_{t\in[T]} \epsilon_t \le \sum_{t\in[T]} \mathbb{E}_{\widehat{\mathbb{P}}}\left[f_t\left(\mathbb{T}_{[t]}^{\epsilon_t}(\widehat{\boldsymbol{\xi}}_{[t]})\right)\right] = \mathbb{E}_{\mathbb{P}^\epsilon}\left[\sum_{t\in[T]} \boldsymbol{c}_t^\top \boldsymbol{x}_t(\boldsymbol{\xi}_{[t]})\right].$$

Letting $\epsilon_2, \ldots, \epsilon_T \downarrow 0$ completes the proof. □

### EC.2.2. Proof of Theorem 2

PROPOSITION EC.2. *Let $p \in [1, \infty)$. Then for any feasible policy $(\boldsymbol{x}_1, \ldots, \boldsymbol{x}_T)$, it holds that*

$$\begin{aligned}
&\sup_{\mathbb{P}\in\mathcal{P}(\Xi)}\left\{\mathbb{E}_{\mathbb{P}}\left[\sum_{t\in[T]} \boldsymbol{c}_t^\top \boldsymbol{x}_t(\boldsymbol{\xi}_{[t]})\right] - \lambda \mathsf{C}_p^p(\widehat{\mathbb{P}}, \mathbb{P})\right\} \\
&= \sup_{\gamma_{[t]}\in\Gamma_{[t]}^{\mathsf{C}}} \mathbb{E}_{\gamma_{[t]}}\left[\sum_{s\in[t]} \boldsymbol{c}_s^\top \boldsymbol{x}_s(\boldsymbol{\xi}_{[s]}) + \mathbb{E}_{\widehat{\mathbb{P}}_{t+1|\widehat{\boldsymbol{\xi}}_{[t]}}}\left[\sup_{\xi_{t+1}\in\Xi_{t+1}}\left\{V_{t+1}^{\widehat{\boldsymbol{\xi}}_{[t+1]}}(\boldsymbol{\xi}_{[t]}, \xi_{t+1}) - \lambda\|\xi_{t+1} - \widehat{\boldsymbol{\xi}}_{t+1}\|^p\right\}\right]\right].
\end{aligned}$$

*Proof of Proposition EC.2.* We prove it by induction. By definition of the causal transport distance, the case of $t=T$ holds trivially. Suppose we have proved for the case of $t$, where $t=2,\dots,T$, now we prove the result also holds for $t-1$. According to the induction hypothesis, we have that

$$\begin{aligned}
V &:= \sup_{\mathbb{P}\in\mathfrak{M}^{\mathsf{C}}} \mathbb{E}_{\mathbb{P}}\left[\sum_{s\in[T]} \boldsymbol{c}_s^\top \boldsymbol{x}_s(\boldsymbol{\xi}_{[s]})\right] \\
&= \sup_{\gamma_{[t]}\in\Gamma^{\mathsf{C}}_{[t]}} \mathbb{E}_{\gamma_{[t]}}\left[\sum_{s\in[t]} \boldsymbol{c}_s^\top \boldsymbol{x}_s(\boldsymbol{\xi}_{[s]}) - \lambda\sum_{s\in[t]}\|\boldsymbol{\xi}_s-\widehat{\boldsymbol{\xi}}_s\|^p + \mathbb{E}_{\widehat{\mathbb{P}}_{t+1|\widehat{\xi}_{[t]}}}\left[\sup_{\xi_{t+1}\in\Xi_{t+1}}\left\{V_{t+1}^{\widehat{\boldsymbol{\xi}}_{[t+1]}}(\boldsymbol{\xi}_{[t]},\xi_{t+1}) - \lambda\|\xi_{t+1}-\widehat{\boldsymbol{\xi}}_{t+1}\|^p\right\}\right]\right] \\
&=: \sup_{\gamma_{[t]}\in\Gamma^{\mathsf{C}}_{[t]}} \mathbb{E}_{\gamma_{[t]}}\left[\sum_{s\in[t]} \boldsymbol{c}_s^\top \boldsymbol{x}_s(\boldsymbol{\xi}_{[s]}) - \lambda\sum_{s\in[t]}\|\boldsymbol{\xi}_s-\widehat{\boldsymbol{\xi}}_s\|^p + \bar{V}_{t+1}^{\widehat{\boldsymbol{\xi}}_{[t]}}(\boldsymbol{\xi}_{[t]})\right].
\end{aligned}$$

For any causal transport plan $\gamma_{[t]}$ from $\widehat{\mathbb{P}}_{[t]}$ to $\mathbb{P}_{[t]}$, by the tower property of conditional expectation, it holds that

$$\begin{aligned}
&\mathbb{E}_{\gamma_{[t]}}\left[\sum_{s\in[t]} \boldsymbol{c}_s^\top \boldsymbol{x}_s(\boldsymbol{\xi}_{[s]}) - \lambda\sum_{s\in[t]}\|\boldsymbol{\xi}_s-\widehat{\boldsymbol{\xi}}_s\|^p + \bar{V}_{t+1}^{\widehat{\boldsymbol{\xi}}_{[t]}}(\boldsymbol{\xi}_{[t]})\right] \\
&= \mathbb{E}_{\gamma_{[t-1]}}\left[\sum_{s\in[t-1]} \boldsymbol{c}_s^\top \boldsymbol{x}_s(\boldsymbol{\xi}_{[s]}) - \lambda\sum_{s\in[t-1]}\|\boldsymbol{\xi}_s-\widehat{\boldsymbol{\xi}}_s\|^p\right] + \mathbb{E}_{\gamma_{[t]}}\left[\boldsymbol{c}_t^\top \boldsymbol{x}_t(\boldsymbol{\xi}_{[t]}) - \lambda\|\boldsymbol{\xi}_t-\widehat{\boldsymbol{\xi}}_t\|^p + \bar{V}_{t+1}^{\widehat{\boldsymbol{\xi}}_{[t]}}(\boldsymbol{\xi}_{[t]})\right] \\
&= \mathbb{E}_{\gamma_{[t-1]}}\left[\sum_{s\in[t-1]} \boldsymbol{c}_s^\top \boldsymbol{x}_s(\boldsymbol{\xi}_{[s]}) - \lambda\sum_{s\in[t-1]}\|\boldsymbol{\xi}_s-\widehat{\boldsymbol{\xi}}_s\|^p\right] \\
&\quad + \mathbb{E}_{\gamma_{[t-1]}}\left[\mathbb{E}_{\gamma_{\widehat{\xi}_t|(\widehat{\xi}_{[t-1]},\xi_{[t-1]})}}\left[\mathbb{E}_{\gamma_{\xi_t|(\widehat{\xi}_{[t]},\xi_{[t-1]})}}\left[\boldsymbol{c}_t^\top \boldsymbol{x}_t(\boldsymbol{\xi}_{[t]}) - \lambda\|\boldsymbol{\xi}_t-\widehat{\boldsymbol{\xi}}_t\|^p + \bar{V}_{t+1}^{\widehat{\boldsymbol{\xi}}_{[t]}}(\boldsymbol{\xi}_{[t]})\right]\right]\right] \\
&= \mathbb{E}_{\gamma_{[t-1]}}\left[\sum_{s\in[t-1]} \boldsymbol{c}_s^\top \boldsymbol{x}_s(\boldsymbol{\xi}_{[s]}) - \lambda\sum_{s\in[t-1]}\|\boldsymbol{\xi}_s-\widehat{\boldsymbol{\xi}}_s\|^p\right] \\
&\quad + \mathbb{E}_{\gamma_{[t-1]}}\left[\mathbb{E}_{\widehat{\mathbb{P}}_{t|\widehat{\xi}_{[t-1]}}}\left[\mathbb{E}_{\gamma_{\xi_t|(\widehat{\xi}_{[t]},\xi_{[t-1]})}}\left[\boldsymbol{c}_t^\top \boldsymbol{x}_t(\boldsymbol{\xi}_{[t]}) - \lambda\|\boldsymbol{\xi}_t-\widehat{\boldsymbol{\xi}}_t\|^p + \bar{V}_{t+1}^{\widehat{\boldsymbol{\xi}}_{[t]}}(\boldsymbol{\xi}_{[t]})\right]\right]\right],
\end{aligned} \tag{EC.3}$$

where, in the last equations above, we have used the non-anticipativity of the causal transport plan. Note that given the distribution $\widehat{\mathbb{P}}_{[t]}$, the joint distribution $\gamma_{[t]}$ is determined completely by $\gamma_{[t-1]}$ and $\gamma_{\xi_t|(\widehat{\xi}_{[t]},\xi_{[t-1]})}$. Thereby maximizing over $\gamma_{[t]}\in\Gamma^{\mathsf{C}}_{[t]}$ is equivalent to maximizing over $\gamma_{[t-1]}\in\Gamma^{\mathsf{C}}_{[t-1]}$ and $\gamma_{\xi_t|(\widehat{\xi}_{[t]},\xi_{[t-1]})}\in\mathcal{P}(\Xi_t)$. Observe that

$$\begin{aligned}
&\sup_{\gamma_{\xi_t|(\widehat{\xi}_{[t]},\xi_{[t-1]})}\in\mathcal{P}(\Xi_t)} \mathbb{E}_{\gamma_{\xi_t|(\widehat{\xi}_{[t]},\xi_{[t-1]})}}\left[\boldsymbol{c}_t^\top \boldsymbol{x}_t(\boldsymbol{\xi}_{[t]}) - \lambda\|\boldsymbol{\xi}_t-\widehat{\boldsymbol{\xi}}_t\|^p + \bar{V}_{t+1}^{\widehat{\boldsymbol{\xi}}_{[t]}}(\boldsymbol{\xi}_{[t]})\right] \\
&= \sup_{\xi_t\in\Xi_t}\left\{c_t^\top \boldsymbol{x}_t(\boldsymbol{\xi}_{[t-1]},\xi_t) + \bar{V}_{t+1}^{\widehat{\boldsymbol{\xi}}_{[t]}}(\boldsymbol{\xi}_{[t-1]},\xi_t) - \lambda\|\xi_t-\widehat{\boldsymbol{\xi}}_t\|^p\right\} \\
&= \sup_{\xi_t\in\Xi_t}\left\{V_t^{\widehat{\boldsymbol{\xi}}_{[t]}}(\boldsymbol{\xi}_{[t-1]},\xi_t) - \lambda\|\xi_t-\widehat{\boldsymbol{\xi}}_t\|^p\right\}.
\end{aligned}$$

Hence, together with (EC.3), it follows that

$$V = \sup_{\gamma_{[t-1]}\in\Gamma^{\mathsf{C}}_{[t-1]}} \mathbb{E}_{\gamma_{[t-1]}}\left[\sum_{s\in[t-1]} \boldsymbol{c}_s^\top \boldsymbol{x}_s(\boldsymbol{\xi}_{[s]}) + \mathbb{E}_{\widehat{\mathbb{P}}_{t|\widehat{\xi}_{[t-1]}}}\left[\sup_{\xi_t\in\Xi_t}\left\{V_t^{\widehat{\boldsymbol{\xi}}_{[t]}}(\boldsymbol{\xi}_{[t-1]},\xi_t) - \lambda\|\xi_t-\widehat{\boldsymbol{\xi}}_t\|^p\right\}\right]\right],$$

which completes the induction. □

*Proof of Theorem 2.* We first show that

$$\sup_{\mathbb{P}\in\mathcal{P}(\Xi)}\left\{\mathbb{E}_{\mathbb{P}}\Big[\sum_{t\in[T]}\boldsymbol{c}_t^\top\boldsymbol{x}_t(\boldsymbol{\xi}_{[t]})\Big]-\lambda\mathsf{C}_p^p(\widehat{\mathbb{P}},\mathbb{P})\right\}=\sup_{\mathbb{P}\in\mathcal{P}(\Xi)}\left\{\mathbb{E}_{\mathbb{P}}\Big[\sum_{t\in[T]}\boldsymbol{c}_t^\top\boldsymbol{x}_t(\boldsymbol{\xi}_{[t]})\Big]-\lambda\mathsf{D}_p^p(\widehat{\mathbb{P}},\mathbb{P})\right\}. \tag{EC.4}$$

Since $\mathfrak{M}\subset\mathfrak{M}^{\mathsf{C}}$, the left-hand side is less than or equal to the right-hand side. It remains to prove the other direction. Using Proposition EC.2, the right-hand side above equals

$$c_1^\top x_1+\mathbb{E}_{\widehat{\mathbb{P}}_{\widehat{\xi}_2}}\Bigg[\sup_{\xi_2\in\Xi_2}\Big\{c_2^\top x_2(\xi_2)+\mathbb{E}_{\widehat{\mathbb{P}}_{3|\widehat{\xi}_2}}\Big[\sup_{\xi_3\in\Xi_3}\Big\{c_3^\top x_3(\xi_{[3]})+\cdots+$$
$$\mathbb{E}_{\widehat{\mathbb{P}}_{T|\widehat{\xi}_{[T-1]}}}\Big[\sup_{\xi_T\in\Xi_T}\big\{c_T^\top x_T(\xi_{[T]})-\lambda\|\xi_T-\widehat{\boldsymbol{\xi}}_T\|^p\big\}\Big]-\cdots-\lambda\|\xi_3-\widehat{\boldsymbol{\xi}}_3\|^p\Big\}\Big]-\lambda\|\xi_2-\widehat{\boldsymbol{\xi}}_2\|^p\Big\}\Bigg].$$

Denote $f_t(\xi_{[t]})=c_t^\top x_t(\xi_{[t]})$. Let $\mathbb{T}_1^{\epsilon_1}$ be the identity transport map. For $t=2,\ldots,T$, for any $\epsilon_t>0$ that is sufficiently small, there exists a transport map $\mathbb{T}_t^{\epsilon_t}$, where $\mathbb{T}_t^{\epsilon_t}:\widehat{\Xi}_{[t]}\to\Xi_t$, such that

$$f_t(\mathbb{T}_{[t]}^{\epsilon_t}(\widehat{\boldsymbol{\xi}}_{[t]}))+\mathbb{E}_{\widehat{\mathbb{P}}_{t+1|\widehat{\xi}_{[t]}}}\Bigg[\sup_{\xi_{t+1}\in\Xi_{t+1}}\Big\{f_{t+1}\Big(\mathbb{T}_{[t]}^{\epsilon_t}(\widehat{\boldsymbol{\xi}}_{[t]}),\xi_{t+1}\Big)+\cdots+\mathbb{E}_{\widehat{\mathbb{P}}_{T|\widehat{\xi}_{[T-1]}}}\Big[\sup_{\xi_T\in\Xi_T}\Big\{f_T\Big(\mathbb{T}_{[t]}^{\epsilon_t}(\widehat{\boldsymbol{\xi}}_{[t]}),\xi_{[t+1,T]}\Big)$$
$$-\lambda\|\xi_T-\widehat{\boldsymbol{\xi}}_T\|^p\Big\}\Big]-\cdots-\lambda\|\xi_{t+1}-\widehat{\boldsymbol{\xi}}_{t+1}\|^p\Big\}\Bigg]-\lambda\|\mathbb{T}_t^{\epsilon_t}(\widehat{\boldsymbol{\xi}}_{[t]})-\widehat{\boldsymbol{\xi}}_t\|^p$$
$$\geq\sup_{\xi_t\in\Xi_t}\Big\{V_t^{\widehat{\xi}_{[t]}}\Big(\mathbb{T}_{[t-1]}^{\epsilon_{[t-1]}}(\widehat{\boldsymbol{\xi}}_{[t-1]}),\xi_t\Big)-\lambda\|\xi_t-\widehat{\boldsymbol{\xi}}_t\|^p\Big\}-\epsilon_t.$$

Let $\mathbb{T}^\epsilon=(\mathbb{T}_1^{\epsilon_1},\ldots,\mathbb{T}_T^{\epsilon_T})$. Define

$$\mathbb{P}^\epsilon:=\mathbb{T}_\#^\epsilon\widehat{\mathbb{P}}.$$

Similar to the reasoning in the proof of Theorem 1, $\mathbb{T}_t^{\epsilon_t}$ can be chosen to be a bi-causal transport map, so that

$$\mathbb{E}_{\mathbb{P}^\epsilon}\Big[\sum_{t\in[T]}\boldsymbol{c}_t^\top\boldsymbol{x}_t(\boldsymbol{\xi}_{[t]})\Big]-\lambda\mathsf{D}_p^p(\widehat{\mathbb{P}},\mathbb{P}^\epsilon)\geq V_1^{\widehat{\xi}_1}-\sum_{t\in[T]}\epsilon_T.$$

Letting $\epsilon_2,\ldots,\epsilon_T\downarrow 0$ completes the proof.

To prove the second part, define

$$\mathcal{R}(\theta):=\begin{cases}\sup_{\mathbb{P}\in\mathcal{P}(\Xi)}\Big\{\mathbb{E}_{\mathbb{P}}\big[\sum_{s\in[T]}\boldsymbol{c}_s^\top\boldsymbol{x}_s(\boldsymbol{\xi}_{[s]})\big]:\ \mathsf{C}_p(\widehat{\mathbb{P}},\mathbb{P})^p\leq\theta\Big\}, & \theta\geq 0,\\ -\infty, & \theta<0.\end{cases}$$

Note that the right-hand side equals to (EC.1) when $\theta=\vartheta^p$. It is not difficult to see that $\mathcal{R}(\cdot)$ is upper semi-continuous and concave. Taking the Legendre transform of $-\mathcal{R}$ gives that for any $\lambda>0$,

$$\begin{aligned}(-\mathcal{R})^*(-\lambda)&=\sup_{\theta\geq 0}\{-\lambda\theta+\mathcal{R}(\theta)\}\\&=\sup_{\theta\geq 0}\sup_{\mathbb{P}\in\mathcal{P}(\Xi)}\left\{\mathbb{E}_{\mathbb{P}}\Big[\sum_{s\in[T]}\boldsymbol{c}_s^\top\boldsymbol{x}_s(\boldsymbol{\xi}_{[s]})\Big]-\lambda\theta:\ \mathsf{C}_p(\widehat{\mathbb{P}},\mathbb{P})^p\leq\theta\right\}\\&=\sup_{\mathbb{P}\in\mathcal{P}(\Xi)}\sup_{\theta\geq 0}\left\{\mathbb{E}_{\mathbb{P}}\Big[\sum_{s\in[T]}\boldsymbol{c}_s^\top\boldsymbol{x}_s(\boldsymbol{\xi}_{[s]})\Big]-\lambda\theta:\mathsf{C}_p(\widehat{\mathbb{P}},\mathbb{P})^p\leq\theta\right\}\\&=\sup_{\mathbb{P}\in\mathcal{P}(\Xi)}\left\{\mathbb{E}_{\mathbb{P}}\Big[\sum_{s\in[T]}\boldsymbol{c}_s^\top\boldsymbol{x}_s(\boldsymbol{\xi}_{[s]})\Big]-\lambda\mathsf{C}_p(\widehat{\mathbb{P}},\mathbb{P})^p\right\}.\end{aligned}$$

Note that for $\lambda = 0$, it also holds that $(-\mathcal{R})^*(0) = \sup_{\mathbb{P}\in\mathcal{P}(\Xi)} \mathbb{E}_{\mathbb{P}}\left[\sum_{s\in[T]} \boldsymbol{c}_s^\top \boldsymbol{x}_s(\boldsymbol{\xi}_{[s]})\right]$ and $(-\mathcal{R})^*(-\lambda)$ is $+\infty$ when $\lambda < 0$. By Fenchel biconjugation theorem, we have that,

$$\mathcal{R}(\theta) = -(-\mathcal{R})^{**}(\theta) = -\sup_{\lambda\geq 0}\{-\lambda\theta - \mathcal{R}^*(-\lambda)\} = \min_{\lambda\geq 0}\left\{\lambda\theta + \sup_{\mathbb{P}\in\mathcal{P}(\Xi)}\left\{\mathbb{E}_{\mathbb{P}}\left[\sum_{s\in[T]} \boldsymbol{c}_s^\top \boldsymbol{x}_s(\boldsymbol{\xi}_{[s]})\right] - \lambda \mathrm{C}_p(\widehat{\mathbb{P}},\mathbb{P})^p\right\}\right\}.$$

Using Proposition EC.2, we obtain that

$$\begin{aligned}
&\sup_{\mathbb{P}\in\mathfrak{M}^{\mathrm{C}}} \mathbb{E}_{\mathbb{P}}\left[\sum_{t\in[T]} \boldsymbol{c}_t^\top \boldsymbol{x}_t(\boldsymbol{\xi}_{[t]})\right] \\
&= \min_{\lambda\geq 0}\Bigg\{\lambda\vartheta^p + c_1^\top x_1 + \mathbb{E}_{\widehat{\mathbb{P}}_2}\Bigg[\sup_{\xi_2\in\Xi_2}\Bigg\{c_2^\top \boldsymbol{x}_2(\xi_2) + \mathbb{E}_{\widehat{\mathbb{P}}_{3|\widehat{\xi}_2}}\Bigg[\sup_{\xi_3\in\Xi_3}\Bigg\{c_3^\top \boldsymbol{x}_3(\xi_{[3]}) + \cdots \\
&\qquad + \mathbb{E}_{\widehat{\mathbb{P}}_{T|\widehat{\boldsymbol{\xi}}_{[T-1]}}}\Big[\sup_{\xi_T\in\Xi_T}\left\{c_T^\top \boldsymbol{x}_T(\xi_{[T]}) - \lambda\|\xi_T - \widehat{\boldsymbol{\xi}}_T\|^p\right\}\Big]\cdots - \lambda\|\xi_3 - \widehat{\boldsymbol{\xi}}_3\|^p\Bigg\}\Bigg] - \lambda\|\xi_2 - \widehat{\boldsymbol{\xi}}_2\|^p\Bigg\}\Bigg]\Bigg\}.
\end{aligned}$$

Finally, we show that

$$\sup_{\mathbb{P}\in\mathfrak{M}} \mathbb{E}_{\mathbb{P}}\left[\sum_{s\in[T]} \boldsymbol{c}_s^\top \boldsymbol{x}_s(\boldsymbol{\xi}_{[s]})\right] = \sup_{\mathbb{P}\in\mathfrak{M}^{\mathrm{C}}} \mathbb{E}_{\mathbb{P}}\left[\sum_{s\in[T]} \boldsymbol{c}_s^\top \boldsymbol{x}_s(\boldsymbol{\xi}_{[s]})\right].$$

The remaining argument is similar to the first part of the proof. □

## EC.3. Proof of Theorem 3

The proof is structured similar to the proof of Theorem 1 and 2. We first consider a relaxed problem where the nested ambiguity set is replaced by the causal ambiguity set, and prove its equivalence to ($\mathsf{P}_{\mathsf{dynamic}}$) in Proposition EC.3 and EC.4, for $p = \infty$ and $p \in [1,\infty)$, respectively. We then show in Propositions EC.5 that the optimal robust cost-to-go function is lower semicontinuous, a condition analogous to the assumption imposed in Theorem 1 and 2. Finally, we complete the proof by establishing the equivalence of the relaxed problem and ($\mathsf{P}_{\mathsf{static}}$) or ($\mathsf{P}_{\mathsf{static}}$-$\mathsf{soft}$), respectively.

Proposition EC.3. *For $p = \infty$, under the setup in Theorem 3, the problem*

$$\inf_{\boldsymbol{x}\in\mathcal{X}} \sup_{\mathbb{P}\in\mathfrak{M}^{\mathrm{C}}} \mathbb{E}_{\mathbb{P}}\left[\sum_{t\in[T]} \boldsymbol{c}_t^\top \boldsymbol{x}_t(\boldsymbol{\xi}_{[t]})\right] \tag{EC.5}$$

*is equivalent to* ($\mathsf{P}_{\mathsf{dynamic}}$).

Proposition EC.4. *For $p \in [1,\infty)$, under the setup in Theorem 3, the problem*

$$\inf_{\boldsymbol{x}\in\mathcal{X}} \sup_{\mathbb{P}\in\mathcal{P}(\Xi_{[T]})} \left\{\mathbb{E}_{\mathbb{P}}\left[\sum_{t\in[T]} \boldsymbol{c}_t^\top \boldsymbol{x}_t(\boldsymbol{\xi}_{[t]})\right] - \lambda \mathrm{C}_p^p(\widehat{\mathbb{P}},\mathbb{P})\right\} \tag{EC.6}$$

*is equivalent to* ($\mathsf{P}_{\mathsf{dynamic}}$).

Proposition EC.5. *The value function*

$$V_t^\star(x_{t-1},\xi_t) := \inf_{x_t\in\mathcal{X}_t(x_{t-1},\xi_t)} \left\{c_t^\top x_t + \mathbb{E}_{\widehat{\mathbb{P}}_{t+1}}\left[Q_{t+1}(x_t,\widehat{\boldsymbol{\xi}}_{t+1})\right]\right\}$$

*is lower semicontinuous in $\xi_t$.*

Note that by expressing $Q_t$ in terms of $V_t$, we have that

$$V_t^\star(x_{t-1},\xi_t) := \begin{cases} \inf_{x_t\in\mathcal{X}_t(x_{t-1},\xi_t)}\left\{c_t^\top x_t + \mathbb{E}_{\widehat{\mathbb{P}}_{t+1}}\left[\sup_{\|\xi_{t+1}-\widehat{\xi}_{t+1}\|\leq\vartheta} V_{t+1}^\star(x_t,\xi_{t+1})\right]\right\}, & p = \infty, \\ \inf_{x_t\in\mathcal{X}_t(x_{t-1},\xi_t)}\left\{c_t^\top x_t + \mathbb{E}_{\widehat{\mathbb{P}}_{t+1}}\left[\sup_{\xi_{t+1}} V_{t+1}^\star(x_t,\xi_{t+1}) - \lambda\|\xi_{t+1} - \widehat{\boldsymbol{\xi}}_{t+1}\|^p\right]\right\}, & p \in [1,\infty). \end{cases}$$

### EC.3.1. Proof of Proposition EC.3

Using Proposition EC.1 and Remark 2, problem (EC.5) is equivalent to

$$\inf_{\boldsymbol{x}_t\in\mathcal{X}_t(\boldsymbol{x}_{t-1}),\forall t\in[T]}\left\{c_1^\top x_1+\mathbb{E}_{\widehat{\mathbb{P}}_2}\left[\sup_{\|\xi_2-\widehat{\xi}_2\|\le\vartheta}\left\{c_2^\top\boldsymbol{x}_2(\xi_2)+\mathbb{E}_{\widehat{\mathbb{P}}_3}\left[\sup_{\|\xi_3-\widehat{\xi}_3\|\le\vartheta}\left\{c_3^\top\boldsymbol{x}_3(\xi_{[3]})+\cdots+\mathbb{E}_{\widehat{\mathbb{P}}_T}\left[\sup_{\|\xi_T-\widehat{\xi}_T\|\le\vartheta}c_T^\top\boldsymbol{x}_T(\xi_{[T]})\right]\cdots\right\}\right]\right\}\right]\right\}.$$

With slight abuse of notation, define

$$\bar V_t(\boldsymbol{x}_{[t-1]}):=\inf_{\boldsymbol{x}_s\in\mathcal{X}_s(\boldsymbol{x}_{s-1}),s=t,\ldots T}\left\{c_1^\top x_1+\mathbb{E}_{\widehat{\mathbb{P}}_2}\left[\sup_{\|\xi_2-\widehat{\xi}_2\|\le\vartheta}\left\{c_2^\top\boldsymbol{x}_2(\xi_2)+\mathbb{E}_{\widehat{\mathbb{P}}_3}\left[\sup_{\|\xi_3-\widehat{\xi}_3\|\le\vartheta}\left\{c_3^\top\boldsymbol{x}_3(\xi_{[3]})+\cdots+\mathbb{E}_{\widehat{\mathbb{P}}_T}\left[\sup_{\|\xi_T-\widehat{\xi}_T\|\le\vartheta}c_T^\top\boldsymbol{x}_T(\xi_{[T]})\right]\cdots\right\}\right]\right\}\right]\right\}.$$

Let us show recursively that for fixed $\boldsymbol{x}_{[t-1]}$,

$$\begin{aligned}&\bar V_t(\boldsymbol{x}_{[t-1]})\\&=c_1^\top x_1+\mathbb{E}_{\widehat{\mathbb{P}}_2}\left[\sup_{\|\xi_2-\widehat{\xi}_2\|\le\vartheta}\left\{c_2^\top\boldsymbol{x}_2(\xi_2)+\mathbb{E}_{\widehat{\mathbb{P}}_3}\left[\sup_{\|\xi_3-\widehat{\xi}_3\|\le\vartheta}\left\{c_3^\top\boldsymbol{x}_3(\xi_{[3]})+\cdots+\right.\right.\right.\right.\\&\left.\left.\left.\left.\mathbb{E}_{\widehat{\mathbb{P}}_{t-1}}\left[\sup_{\|\xi_{t-1}-\widehat{\xi}_{t-1}\|\le\vartheta}\left\{c_{t-1}^\top\boldsymbol{x}_{t-1}(\xi_{[t-1]})+\mathcal{Q}_t(\boldsymbol{x}_{t-1}(\xi_{[t-1]}),\xi_{[t-1]})\right\}\right]\cdots\right\}\right]\right\}\right] &\text{(EC.7)}\\&=c_1^\top x_1+\mathbb{E}_{\widehat{\mathbb{P}}_2}\left[\sup_{\|\xi_2-\widehat{\xi}_2\|\le\vartheta}\left\{c_2^\top\boldsymbol{x}_2(\xi_2)+\mathbb{E}_{\widehat{\mathbb{P}}_3}\left[\sup_{\|\xi_3-\widehat{\xi}_3\|\le\vartheta}\left\{c_3^\top\boldsymbol{x}_3(\xi_{[3]})+\cdots+\right.\right.\right.\right.\\&\left.\left.\left.\left.\mathbb{E}_{\widehat{\mathbb{P}}_{t-1}}\left[\sup_{\|\xi_{t-1}-\widehat{\xi}_{t-1}\|\le\vartheta}\left\{c_{t-1}^\top\boldsymbol{x}_{t-1}(\xi_{[t-1]})+\mathbb{E}_{\widehat{\mathbb{P}}_t}\left[Q_t(\boldsymbol{x}_{t-1}(\xi_{[t-1]}),\widehat{\boldsymbol{\xi}}_t)\right]\right\}\right]\cdots\right\}\right]\right\}\right]. &\text{(EC.8)}\end{aligned}$$

We first prove (EC.8). The case of $t=T+1$ holds trivially. Suppose it holds for some $t+1$, $t=2,\ldots,T$, then

$$\bar V_t(\boldsymbol{x}_{[t-1]})=\inf_{\boldsymbol{x}_t\in\mathcal{X}_t(\boldsymbol{x}_{t-1})}\left\{c_1^\top x_1+\mathbb{E}_{\widehat{\mathbb{P}}_2}\left[\sup_{\|\xi_2-\widehat{\xi}_2\|\le\vartheta}\left\{c_2^\top\boldsymbol{x}_2(\xi_2)+\mathbb{E}_{\widehat{\mathbb{P}}_3}\left[\sup_{\|\xi_3-\widehat{\xi}_3\|\le\vartheta}\left\{c_3^\top\boldsymbol{x}_3(\xi_{[3]})+\cdots+\mathbb{E}_{\widehat{\mathbb{P}}_t}\left[\sup_{\|\xi_t-\widehat{\xi}_t\|\le\vartheta}\left\{c_t^\top\boldsymbol{x}_t(\xi_{[t]})+\mathbb{E}_{\widehat{\mathbb{P}}_{t+1}}\left[Q_{t+1}(\boldsymbol{x}_t(\xi_{[t]}),\widehat{\boldsymbol{\xi}}_{t+1})\right]\right\}\right]\cdots\right\}\right]\right\}\right]\right\}.$$

We next show it holds also for $t$. To this end, define $g_t(x_t,\xi_{[t]}):\mathbb{R}^{d_{x_t}}\times\Xi_{[t]}\to\mathbb{R}\cup\{+\infty\}$ as

$$g_t(x_t,\xi_{[t]}):=\begin{cases}c_t^\top x_t+\mathbb{E}_{\widehat{\mathbb{P}}_{t+1}}\left[Q_{t+1}(x_t,\widehat{\boldsymbol{\xi}}_{t+1})\right], & \text{if } x_t\in\mathcal{X}_t(\boldsymbol{x}_{t-1}(\xi_{[t-1]}),\xi_t),\\ +\infty, & \text{otherwise.}\end{cases}$$

Then $g$ is random lower semi-continuous (Shapiro et al. 2021, Definition 9.47). Denote by $\mathcal{Z}_{[t]}$ the space of measurable functions from $\Xi_{[t]}$ to $\mathbb{R}$ and define $\mathcal{R}_{[t]}:\mathcal{Z}_{[t]}\to\mathbb{R}\cup\{\infty\}$ as

$$\mathcal{R}_{[t]}(\nu):=c_1^\top x_1+\mathbb{E}_{\widehat{\mathbb{P}}_2}\left[\sup_{\|\xi_2-\widehat{\xi}_2\|\le\vartheta}\left\{c_2^\top\boldsymbol{x}_2(\xi_2)+\mathbb{E}_{\widehat{\mathbb{P}}_3}\left[\sup_{\|\xi_3-\widehat{\xi}_3\|\le\vartheta}\left\{c_3^\top\boldsymbol{x}_3(\xi_{[3]})+\cdots+\mathbb{E}_{\widehat{\mathbb{P}}_t}\left[\sup_{\|\xi_t-\widehat{\xi}_t\|\le\vartheta}\nu(\xi_{[t]})\right]\cdots\right\}\right]\right\}\right].$$

By definition, $\mathcal{R}_{[t]}$ is monotone and continuous with respect to the $L^\infty$-norm. Let $X_t$ be a subset of $\mathbb{R}^{d_{x_t}}$ containing $\bigcup_{x_{t-1}\in\mathbb{R}^{d_{x_{t-1}}},\xi_t\in\Xi_t}\mathcal{X}_t(x_{t-1},\xi_t)$ and $\mathfrak{X}_t$ be the set of measurable functions from $\Xi_{[t]}$ to $X_t$. It holds that

$$\begin{aligned}
&\inf_{\chi_t\in\mathfrak{X}_t}\mathcal{R}_{[t]}(g_t(\chi_t(\cdot),\cdot))\\
&=\inf_{\chi_t\in\mathfrak{X}_t}\left\{c_1^\top x_1+\mathbb{E}_{\widehat{\mathbb{P}}_2}\left[\sup_{\|\xi_2-\widehat{\xi}_2\|\le\vartheta}\left\{c_2^\top \boldsymbol{x}_2(\xi_2)+\mathbb{E}_{\widehat{\mathbb{P}}_3}\left[\sup_{\|\xi_3-\widehat{\xi}_3\|\le\vartheta}\left\{c_3^\top \boldsymbol{x}_3(\xi_{[3]})+\cdots+\right.\right.\right.\right.\right.\\
&\qquad\left.\left.\left.\left.\left.\mathbb{E}_{\widehat{\mathbb{P}}_t}\left[\sup_{\|\xi_t-\widehat{\xi}_t\|\le\vartheta}g_t(\chi_t(\xi_{[t]}),\xi_{[t]})\right]\cdots\right\}\right]\right\}\right]\right\}\\
&=\inf_{\substack{\chi_t(\xi_{[t]})\in\mathcal{X}_t(\boldsymbol{x}_{t-1}(\xi_{[t-1]}),\xi_t)\\ \forall\xi_{[t]}}}\left\{c_1^\top x_1+\mathbb{E}_{\widehat{\mathbb{P}}_2}\left[\sup_{\|\xi_2-\widehat{\xi}_2\|\le\vartheta}\left\{c_2^\top \boldsymbol{x}_2(\xi_2)+\mathbb{E}_{\widehat{\mathbb{P}}_3}\left[\sup_{\|\xi_3-\widehat{\xi}_3\|\le\vartheta}\left\{c_3^\top \boldsymbol{x}_3(\xi_{[3]})+\cdots+\right.\right.\right.\right.\right.\\
&\qquad\left.\left.\left.\left.\left.\mathbb{E}_{\widehat{\mathbb{P}}_t}\left[\sup_{\|\xi_t-\widehat{\xi}_t\|\le\vartheta}\left\{c_t^\top\chi_t(\xi_{[t]})+\mathbb{E}_{\widehat{\mathbb{P}}_{t+1}}\left[Q_{t+1}(\chi_t(\xi_{[t]}),\widehat{\xi}_{t+1})\right]\right\}\right]\cdots\right\}\right]\right\}\right]\right\}\\
&=\inf_{\boldsymbol{x}_t\in\mathcal{X}_t(\boldsymbol{x}_{t-1})}\left\{c_1^\top x_1+\mathbb{E}_{\widehat{\mathbb{P}}_2}\left[\sup_{\|\xi_2-\widehat{\xi}_2\|\le\vartheta}\left\{c_2^\top \boldsymbol{x}_2(\xi_2)+\mathbb{E}_{\widehat{\mathbb{P}}_3}\left[\sup_{\|\xi_3-\widehat{\xi}_3\|\le\vartheta}\left\{c_3^\top \boldsymbol{x}_3(\xi_{[3]})+\cdots+\right.\right.\right.\right.\right.\\
&\qquad\left.\left.\left.\left.\left.\mathbb{E}_{\widehat{\mathbb{P}}_t}\left[\sup_{\|\xi_t-\widehat{\xi}_t\|\le\vartheta}\left\{c_t^\top\boldsymbol{x}_t(\xi_{[t]})+\mathbb{E}_{\widehat{\mathbb{P}}_{t+1}}\left[Q_{t+1}(\boldsymbol{x}_t(\xi_{[t]}),\widehat{\xi}_{t+1})\right]\right\}\right]\cdots\right\}\right]\right\}\right]\right\}\\
&=\bar{V}_t(\boldsymbol{x}_{[t-1]}).
\end{aligned}$$

Define $\nu_t\in\mathcal{Z}_{[t]}$ as

$$\nu_t(\xi_{[t]}):=\inf_{x_t\in\mathcal{X}_t(\boldsymbol{x}_{t-1}(\xi_{[t-1]}),\xi_t)}g_t(x_t,\xi_{[t]}),\quad \xi_{[t]}\in\Xi_{[t]}.$$

Using the interchangeability principle (Shapiro et al. 2021, Theorem 9.110), we have

$$\bar{V}_t(\boldsymbol{x}_{[t-1]})=\inf_{\chi_t\in\mathfrak{X}_t}\mathcal{R}(g_t(\chi_t(\cdot),\cdot))=\mathcal{R}(\nu_t),$$

and observe that

$$\mathcal{R}(\nu_t)=c_1^\top x_1+\mathbb{E}_{\widehat{\mathbb{P}}_2}\left[\sup_{\|\xi_2-\widehat{\xi}_2\|\le\vartheta}\left\{c_2^\top \boldsymbol{x}_2(\xi_2)+\mathbb{E}_{\widehat{\mathbb{P}}_3}\left[\sup_{\|\xi_3-\widehat{\xi}_3\|\le\vartheta}\left\{c_3^\top \boldsymbol{x}_3(\xi_{[3]})+\cdots+\mathbb{E}_{\widehat{\mathbb{P}}_t}\left[Q_t(\boldsymbol{x}_{t-1},\widehat{\xi}_t)\right]\cdots\right\}\right]\right\}\right].$$

This completes the induction for (EC.8), which also shows that the optimal value of (1) is

$$\inf_{x_1\in\mathcal{X}_1}\left\{c_1^\top x_1+\mathbb{E}_{\widehat{\mathbb{P}}_2}[Q_2(x_1,\widehat{\xi}_2)]\right\}.$$

Next, we prove (EC.7), or equivalently,

$$\mathbb{E}_{\widehat{\mathbb{P}}_t}\left[Q_t(x_{t-1},\widehat{\xi}_t)\right]=\mathcal{Q}_t(x_{t-1},\xi_{[t-1]}),$$

which, together with (EC.8), implies (EC.7). We prove it by induction. The base case $t=T+1$ trivially holds. Suppose (EC.7) holds for some $t+1$, $t=2,\ldots,T$, then

$$\begin{aligned}
\bar{V}_t(\boldsymbol{x}_{[t-1]})=&\inf_{\boldsymbol{x}_t\in\mathcal{X}_t(\boldsymbol{x}_{t-1})}\left\{c_1^\top x_1+\mathbb{E}_{\widehat{\mathbb{P}}_2}\left[\sup_{\|\xi_2-\widehat{\xi}_2\|\le\vartheta}\left\{c_2^\top \boldsymbol{x}_2(\xi_2)+\mathbb{E}_{\widehat{\mathbb{P}}_3}\left[\sup_{\|\xi_3-\widehat{\xi}_3\|\le\vartheta}\left\{c_3^\top \boldsymbol{x}_3(\xi_{[3]})+\cdots+\right.\right.\right.\right.\right.\\
&\qquad\left.\left.\left.\left.\left.\mathbb{E}_{\widehat{\mathbb{P}}_t}\left[\sup_{\|\xi_t-\widehat{\xi}_t\|\le\vartheta}\left\{c_t^\top\boldsymbol{x}_t(\xi_{[t]})+\mathcal{Q}_{t+1}(\boldsymbol{x}_t(\xi_{[t]}),\xi_{[t]})\right\}\right]\cdots\right\}\right]\right\}\right]\right\}.
\end{aligned}$$

Exchanging $\inf_{\boldsymbol{x}_t}$ and $\mathbb{E}_{\widehat{\mathbb{P}}_s}[\sup_{\xi_s}]$, $s=2,\dots,t-1$, it follows that

$$\begin{aligned}
&\bar{V}_t(\boldsymbol{x}_{[t-1]}) \\
&\geq c_1^\top x_1 + \mathbb{E}_{\widehat{\mathbb{P}}_2}\Bigg[\sup_{\|\xi_2-\widehat{\xi}_2\|\leq\vartheta}\Big\{c_2^\top \boldsymbol{x}_2(\xi_2) + \mathbb{E}_{\widehat{\mathbb{P}}_3}\Big[\sup_{\|\xi_3-\widehat{\xi}_3\|\leq\vartheta}\Big\{c_3^\top \boldsymbol{x}_3(\xi_{[3]}) + \cdots + \\
&\qquad \inf_{\boldsymbol{x}_t\in\mathcal{X}_t(\boldsymbol{x}_{t-1})} \mathbb{E}_{\widehat{\mathbb{P}}_t}\Big[\sup_{\|\xi_t-\widehat{\xi}_t\|\leq\vartheta}\big\{c_t^\top \boldsymbol{x}_t(\xi_{[t]}) + \mathcal{Q}_{t+1}(\boldsymbol{x}_t(\xi_{[t]}),\xi_{[t]})\big\}\Big]\cdots\Big\}\Big]\Big\}\Bigg] \\
&= c_1^\top x_1 + \mathbb{E}_{\widehat{\mathbb{P}}_2}\Bigg[\sup_{\|\xi_2-\widehat{\xi}_2\|\leq\vartheta}\Big\{c_2^\top \boldsymbol{x}_2(\xi_2) + \mathbb{E}_{\widehat{\mathbb{P}}_3}\Big[\sup_{\|\xi_3-\widehat{\xi}_3\|\leq\vartheta}\Big\{c_3^\top \boldsymbol{x}_3(\xi_{[3]}) + \cdots + \\
&\qquad \inf_{\boldsymbol{x}_t\in\mathcal{X}_t(\boldsymbol{x}_{t-1})} \sup_{\mathbb{P}_t\in\mathcal{P}(\Xi_t):\mathcal{W}_\infty(\widehat{\mathbb{P}}_t,\mathbb{P}_t)\leq\vartheta} \mathbb{E}_{\mathbb{P}_t}\Big[\boldsymbol{c}_t^\top \boldsymbol{x}_t(\boldsymbol{\xi}_{[t]}) + \mathcal{Q}_{t+1}(\boldsymbol{x}_t(\boldsymbol{\xi}_{[t]}),\boldsymbol{\xi}_{[t]})\Big]\cdots\Big\}\Big]\Big\}\Bigg] \\
&= c_1^\top x_1 + \mathbb{E}_{\widehat{\mathbb{P}}_2}\Bigg[\sup_{\|\xi_2-\widehat{\xi}_2\|\leq\vartheta}\Big\{c_2^\top \boldsymbol{x}_2(\xi_2) + \mathbb{E}_{\widehat{\mathbb{P}}_3}\Big[\sup_{\|\xi_3-\widehat{\xi}_3\|\leq\vartheta}\Big\{c_3^\top \boldsymbol{x}_3(\xi_{[3]}) + \cdots + \mathcal{Q}_t(\boldsymbol{x}_{t-1}(\boldsymbol{\xi}_{[t-1]}),\boldsymbol{\xi}_{[t-1]})\cdots\Big\}\Big]\Big\}\Bigg],
\end{aligned}$$

where we have used the duality for Wasserstein DRO in the first equality above. This, together with (EC.8), shows that

$$\mathbb{E}_{\widehat{\mathbb{P}}_t}\big[Q_t(x_{t-1},\widehat{\boldsymbol{\xi}}_t)\big] \geq \mathcal{Q}_t(x_{t-1},\xi_{[t-1]}).$$

To prove the other direction, by definition of $Q_t$, we have that

$$\begin{aligned}
&\mathbb{E}_{\widehat{\mathbb{P}}_t}\Big[Q_t(x_{t-1},\widehat{\boldsymbol{\xi}}_t)\Big] \\
&= \mathbb{E}_{\widehat{\mathbb{P}}_t}\Big[\sup_{\|\xi_t-\widehat{\xi}_t\|\leq\vartheta}\ \inf_{x_t\in\mathcal{X}_t(x_{t-1},\xi_t)}\big\{c_t^\top x_t + \mathbb{E}_{\widehat{\mathbb{P}}_{t+1}}\big[Q_{t+1}(x_t,\widehat{\boldsymbol{\xi}}_{t+1})\big]\big\}\Big] \\
&= \sup_{\mathbb{P}_t\in\mathcal{P}(\Xi_t):\mathcal{W}_\infty(\widehat{\mathbb{P}}_t,\mathbb{P}_t)\leq\vartheta} \mathbb{E}_{\mathbb{P}_t}\Big[\inf_{x_t\in\mathcal{X}_t(x_{t-1},\xi_t)}\big\{c_t^\top x_t + \mathbb{E}_{\widehat{\mathbb{P}}_{t+1}}\big[Q_{t+1}(x_t,\widehat{\boldsymbol{\xi}}_{t+1})\big]\big\}\Big] \\
&= \sup_{\mathbb{P}_t\in\mathcal{P}(\Xi_t):\mathcal{W}_\infty(\widehat{\mathbb{P}}_t,\mathbb{P}_t)\leq\vartheta}\ \inf_{\boldsymbol{x}_t(\xi_{[t-1]},\cdot)\in\mathcal{X}_t(x_{t-1},\cdot)} \mathbb{E}_{\mathbb{P}_t}\Big[c_t^\top \boldsymbol{x}_t(\xi_{[t-1]},\boldsymbol{\xi}_t) + \mathbb{E}_{\widehat{\mathbb{P}}_{t+1}}\big[Q_{t+1}\big(\boldsymbol{x}_t(\xi_{[t-1]},\boldsymbol{\xi}_t),\widehat{\boldsymbol{\xi}}_{t+1}\big)\big]\Big] \\
&\leq \inf_{\boldsymbol{x}_t(\xi_{[t-1]},\cdot)\in\mathcal{X}_t(x_{t-1},\cdot)}\ \sup_{\mathbb{P}_t\in\mathcal{P}(\Xi_t):\mathcal{W}_\infty(\widehat{\mathbb{P}}_t,\mathbb{P}_t)\leq\vartheta} \mathbb{E}_{\mathbb{P}_t}\Big[c_t^\top \boldsymbol{x}_t(\xi_{[t-1]},\boldsymbol{\xi}_t) + \mathbb{E}_{\widehat{\mathbb{P}}_{t+1}}\big[Q_{t+1}\big(\boldsymbol{x}_t(\xi_{[t-1]},\boldsymbol{\xi}_t),\widehat{\boldsymbol{\xi}}_{t+1}\big)\big]\Big] \\
&= \mathcal{Q}_t(x_{t-1},\xi_{[t-1]}).
\end{aligned}$$

where the second equality follows from the duality of Wasserstein DRO (Zhang et al. 2022, Section 5), where we use our Assumption 3 to validate Assumption 1 in Zhang et al. (2022). The third equality follows from the interchangeability principle (Shapiro et al. 2021, Theorem 9.108). □

### EC.3.2. Proof of Proposition EC.4

The proof is similar to the case of $p=\infty$, where we adopt notations such as $X_t$, $\mathfrak{X}_t$ and $\mathcal{Z}_{[t]}$. Using Remark 2, problem ($\mathsf{P}_{\mathsf{static}}$-soft) is equivalent to

$$\begin{aligned}
\inf_{\boldsymbol{x}_t\in\mathcal{X}_t(\boldsymbol{x}_{t-1}),\forall t\in[T]}\Bigg\{&c_1^\top x_1 + \mathbb{E}_{\widehat{\mathbb{P}}_2}\Bigg[\sup_{\xi_2\in\Xi_2}\Bigg\{c_2^\top \boldsymbol{x}_2(\xi_2) + \mathbb{E}_{\widehat{\mathbb{P}}_3}\Big[\sup_{\xi_3\in\Xi_3}\Big\{c_3^\top \boldsymbol{x}_3(\xi_{[3]}) + \\
&\cdots + \mathbb{E}_{\widehat{\mathbb{P}}_T}\Big[\sup_{\xi_T\in\Xi_T}\Big\{c_T^\top \boldsymbol{x}_T(\xi_{[T]}) - \lambda\|\xi_T-\widehat{\boldsymbol{\xi}}_T\|^p\Big\}\Big]\cdots - \lambda\|\xi_3-\widehat{\boldsymbol{\xi}}_3\|^p\Big\}\Big] - \lambda\|\xi_2-\widehat{\boldsymbol{\xi}}_2\|^p\Bigg\}\Bigg]\Bigg\}.
\end{aligned}$$

Let us show recursively that for fixed $\boldsymbol{x}_{[t-1]}$,

$$\begin{aligned}
&V_t(\boldsymbol{x}_{[t-1]})\\
&:= \inf_{\boldsymbol{x}_s\in\mathcal{X}_s(\boldsymbol{x}_{s-1}),s=t,\dots T}\Bigg\{c_1^\top x_1 + \mathbb{E}_{\widehat{\mathbb{P}}_2}\Bigg[\sup_{\xi_2\in\Xi_2}\Bigg\{c_2^\top \boldsymbol{x}_2(\xi_2) + \mathbb{E}_{\widehat{\mathbb{P}}_3}\Bigg[\sup_{\xi_3\in\Xi_3}\Bigg\{c_3^\top \boldsymbol{x}_3(\xi_{[3]}) +\\
&\qquad \cdots + \mathbb{E}_{\widehat{\mathbb{P}}_T}\Big[\sup_{\xi_T\in\Xi_T}\Big\{c_T^\top \boldsymbol{x}_T(\xi_{[T]}) - \lambda\|\xi_T - \widehat{\boldsymbol{\xi}}_T\|^p\Big\}\Big]\cdots - \lambda\|\xi_3 - \widehat{\boldsymbol{\xi}}_3\|^p\Bigg\}\Bigg] - \lambda\|\xi_2 - \widehat{\boldsymbol{\xi}}_2\|^p\Bigg\}\Bigg]\Bigg\}\\
&= c_1^\top x_1 + \mathbb{E}_{\widehat{\mathbb{P}}_2}\Bigg[\sup_{\xi_2\in\Xi_2}\Bigg\{c_2^\top \boldsymbol{x}_2(\xi_2) + \mathbb{E}_{\widehat{\mathbb{P}}_3}\Bigg[\sup_{\xi_3\in\Xi_3}\Bigg\{c_3^\top \boldsymbol{x}_3(\xi_{[3]}) + \cdots + \mathbb{E}_{\widehat{\mathbb{P}}_{t-1}}\Big[\sup_{\xi_{t-1}\in\Xi_{t-1}}\Big\{c_{t-1}^\top \boldsymbol{x}_{t-1}(\xi_{[t-1]}) +\\
&\qquad \mathcal{Q}_t(\boldsymbol{x}_{t-1}(\xi_{[t-1]}),\xi_{[t-1]}) - \lambda\|\xi_{t-1} - \widehat{\boldsymbol{\xi}}_{t-1}\|^p\Big\}\Big]\cdots - \lambda\|\xi_3 - \widehat{\boldsymbol{\xi}}_3\|^p\Bigg\}\Bigg] - \lambda\|\xi_2 - \widehat{\boldsymbol{\xi}}_2\|^p\Bigg\}\Bigg]
\end{aligned} \tag{EC.9}$$

$$\begin{aligned}
&= c_1^\top x_1 + \mathbb{E}_{\widehat{\mathbb{P}}_2}\Bigg[\sup_{\xi_2\in\Xi_2}\Bigg\{c_2^\top \boldsymbol{x}_2(\xi_2) + \mathbb{E}_{\widehat{\mathbb{P}}_3}\Bigg[\sup_{\xi_3\in\Xi_3}\Bigg\{c_3^\top \boldsymbol{x}_3(\xi_{[3]}) + \cdots + \mathbb{E}_{\widehat{\mathbb{P}}_{t-1}}\Big[\sup_{\xi_{t-1}\in\Xi_{t-1}}\Big\{c_{t-1}^\top \boldsymbol{x}_{t-1}(\xi_{[t-1]}) +\\
&\qquad \mathbb{E}_{\widehat{\mathbb{P}}_t}\big[\mathcal{Q}_t(\boldsymbol{x}_{t-1}(\xi_{[t-1]}),\widehat{\boldsymbol{\xi}}_t)\big] - \lambda\|\xi_{t-1} - \widehat{\boldsymbol{\xi}}_{t-1}\|^p\Big\}\Big]\cdots - \lambda\|\xi_3 - \widehat{\boldsymbol{\xi}}_3\|^p\Bigg\}\Bigg] - \lambda\|\xi_2 - \widehat{\boldsymbol{\xi}}_2\|^p\Bigg\}\Bigg].
\end{aligned} \tag{EC.10}$$

The case of $t=T+1$ holds trivially. Suppose (EC.10) holds for some $t+1$, $t=2,\dots,T$, now we prove for $t$. Using the induction hypothesis, it holds that

$$\begin{aligned}
&V_t(\boldsymbol{x}_{[t-1]})\\
&= \inf_{\boldsymbol{x}_t\in\mathcal{X}_t(\boldsymbol{x}_{t-1})}\Bigg\{c_1^\top x_1 + \mathbb{E}_{\widehat{\mathbb{P}}_2}\Bigg[\sup_{\xi_2\in\Xi_2}\Bigg\{c_2^\top \boldsymbol{x}_2(\xi_2) + \mathbb{E}_{\widehat{\mathbb{P}}_3}\Bigg[\sup_{\xi_3\in\Xi_3}\Bigg\{c_3^\top \boldsymbol{x}_3(\xi_{[3]}) + \cdots +\\
&\quad \mathbb{E}_{\widehat{\mathbb{P}}_t}\Big[\sup_{\xi_t\in\Xi_t}\Big\{c_t^\top \boldsymbol{x}_t(\xi_{[t]}) + \mathbb{E}_{\widehat{\mathbb{P}}_{t+1}}\big[\mathcal{Q}_{t+1}(\boldsymbol{x}_t(\xi_{[t]}),\widehat{\boldsymbol{\xi}}_{t+1})\big] - \lambda\|\xi_t - \widehat{\boldsymbol{\xi}}_t\|^p\Big\}\Big]\cdots - \lambda\|\xi_3 - \widehat{\boldsymbol{\xi}}_3\|^p\Bigg\}\Bigg] - \lambda\|\xi_2 - \widehat{\boldsymbol{\xi}}_2\|^p\Bigg\}\Bigg]\Bigg\}.
\end{aligned}$$

Define

$$g_t(x_t,\xi_{[t]}) := \begin{cases} c_t^\top x_t + \mathbb{E}_{\widehat{\mathbb{P}}_{t+1}}\big[\mathcal{Q}_{t+1}(x_t,\widehat{\boldsymbol{\xi}}_{t+1})\big] - \lambda\|\xi_t - \widehat{\boldsymbol{\xi}}_t\|^p, & \text{if } x_t \in \mathcal{X}_t(\boldsymbol{x}_{t-1}(\xi_{[t-1]}),\xi_t),\\ +\infty, & \text{otherwise,}\end{cases}$$

$$\begin{aligned}
\mathcal{R}_{[t]}(\nu) := c_1^\top x_1 + \mathbb{E}_{\widehat{\mathbb{P}}_2}\Bigg[\sup_{\xi_2\in\Xi_2}\Bigg\{c_2^\top \boldsymbol{x}_2(\xi_2) + \mathbb{E}_{\widehat{\mathbb{P}}_3}\Bigg[\sup_{\xi_3\in\Xi_3}\Bigg\{c_3^\top \boldsymbol{x}_3(\xi_{[3]}) + \cdots +\\
\mathbb{E}_{\widehat{\mathbb{P}}_t}\Big[\sup_{\xi_t\in\Xi_t}\Big\{\nu(\xi_{[t]}) - \lambda\|\xi_t - \widehat{\boldsymbol{\xi}}_t\|^p\Big\}\Big]\cdots - \lambda\|\xi_3 - \widehat{\boldsymbol{\xi}}_3\|^p\Bigg\}\Bigg] - \lambda\|\xi_2 - \widehat{\boldsymbol{\xi}}_2\|^p\Bigg\}\Bigg],
\end{aligned}$$

and

$$\nu_t(\xi_{[t]}) := \inf_{x_t\in\mathcal{X}_t(\boldsymbol{x}_{t-1}(\xi_{[t-1]}),\xi_t)} g_t(x_t,\xi_{[t]}), \quad \xi_{[t]}\in\Xi_{[t]}.$$

Then $g$ is random lower semi-continuous, $\mathcal{R}_{[t]}$ is monotone and continuous with respect to the $L^\infty$-norm.

Using the interchangeability principle (Shapiro et al. 2021, Theorem 9.110), we have that

$$\begin{aligned}V_t(\boldsymbol{x}_{[t-1]}) &= \inf_{\chi_t \in \mathfrak{X}_t} \mathcal{R}_{[t]}(g_t(\chi_t(\cdot),\cdot))\\ &= \mathcal{R}_{[t]}(\nu_t)\\ &= c_1^\top x_1 + \mathbb{E}_{\widehat{\mathbb{P}}_2}\left[\sup_{\xi_2\in\Xi_2}\left\{c_2^\top \boldsymbol{x}_2(\xi_2) + \mathbb{E}_{\widehat{\mathbb{P}}_3}\left[\sup_{\xi_3\in\Xi_3}\left\{c_3^\top \boldsymbol{x}_3(\xi_{[3]}) + \cdots + \mathbb{E}_{\widehat{\mathbb{P}}_t}\left[Q_t(\boldsymbol{x}_{t-1},\widehat{\boldsymbol{\xi}}_t)\right]\cdots\right.\right.\right.\right.\\ &\qquad\left.\left.\left.\left. - \lambda\|\xi_3 - \widehat{\boldsymbol{\xi}}_3\|^p\right\}\right] - \lambda\|\xi_2-\widehat{\boldsymbol{\xi}}_2\|^p\right\}\right].\end{aligned}$$

This completes the induction for (EC.8).

Next, we prove (EC.9) by showing that

$$\mathbb{E}_{\widehat{\mathbb{P}}_t}\left[Q_t(x_{t-1},\widehat{\boldsymbol{\xi}}_t)\right] = \mathcal{Q}_t(x_{t-1},\xi_{[t-1]}).$$

The base case $t=T+1$ trivially holds. Suppose we have shown the case for some $t+1$, $t=2,\dots,T$. Exchanging $\inf_{\boldsymbol{x}_t}$ and $\mathbb{E}_{\widehat{\mathbb{P}}_s}[\sup_{\xi_s}]$, $s=2,\dots,t-1$, we obtain that

$$\begin{aligned}&V_t(\boldsymbol{x}_{[t-1]})\\ &\geq c_1^\top x_1 + \mathbb{E}_{\widehat{\mathbb{P}}_2}\left[\sup_{\xi_2\in\Xi_2}\left\{c_2^\top \boldsymbol{x}_2(\xi_2) + \mathbb{E}_{\widehat{\mathbb{P}}_3}\left[\sup_{\xi_3\in\Xi_3}\left\{c_3^\top \boldsymbol{x}_3(\xi_{[3]}) + \cdots + \inf_{\boldsymbol{x}_t\in\mathcal{X}_t(\boldsymbol{x}_{t-1})}\right.\right.\right.\right.\\ &\qquad \mathbb{E}_{\widehat{\mathbb{P}}_t}\left[\sup_{\xi_t\in\Xi_t}\left\{c_t^\top \boldsymbol{x}_t(\xi_{[t]}) + \mathcal{Q}_{t+1}(\boldsymbol{x}_t(\xi_{[t]}),\xi_{[t]}) - \lambda\|\xi_t - \widehat{\boldsymbol{\xi}}_t\|^p\right\}\right]\cdots \left.\left.\left.\left.- \lambda\|\xi_3-\widehat{\boldsymbol{\xi}}_3\|^p\right\}\right] - \lambda\|\xi_2-\widehat{\boldsymbol{\xi}}_2\|^p\right\}\right]\\ &= c_1^\top x_1 + \mathbb{E}_{\widehat{\mathbb{P}}_2}\left[\sup_{\xi_2\in\Xi_2}\left\{c_2^\top \boldsymbol{x}_2(\xi_2) + \mathbb{E}_{\widehat{\mathbb{P}}_3}\left[\sup_{\xi_3\in\Xi_3}\left\{c_3^\top \boldsymbol{x}_3(\xi_{[3]}) + \cdots + \inf_{\boldsymbol{x}_t\in\mathcal{X}_t(\boldsymbol{x}_{t-1})}\sup_{\mathbb{P}_t\in\mathcal{P}(\Xi_t)}\left\{\mathbb{E}_{\mathbb{P}_t}\left[\boldsymbol{c}_t^\top \boldsymbol{x}_t(\xi_{[t-1]},\boldsymbol{\xi}_t)\right.\right.\right.\right.\right.\right.\\ &\qquad \left.+\mathcal{Q}_{t+1}\big(\boldsymbol{x}_t(\xi_{[t-1]},\boldsymbol{\xi}_t),(\xi_{[t-1]},\boldsymbol{\xi}_t)\big)\right] - \lambda\mathcal{W}_p^p(\widehat{\mathbb{P}}_t,\mathbb{P}_t)\Big\}\cdots \left.\left.\left.\left.- \lambda\|\xi_3-\widehat{\boldsymbol{\xi}}_3\|^p\right\}\right] - \lambda\|\xi_2-\widehat{\boldsymbol{\xi}}_2\|^p\right\}\right].\\ &= c_1^\top x_1 + \mathbb{E}_{\widehat{\mathbb{P}}_2}\left[\sup_{\xi_2\in\Xi_2}\left\{c_2^\top \boldsymbol{x}_2(\xi_2) + \mathbb{E}_{\widehat{\mathbb{P}}_3}\left[\sup_{\xi_3\in\Xi_3}\left\{c_3^\top \boldsymbol{x}_3(\xi_{[3]}) + \cdots + \mathcal{Q}_t(\boldsymbol{x}_{t-1}(\xi_{[t-1]}),\xi_{t-1})\right.\right.\right.\right.\\ &\qquad \cdots \left.\left.\left.\left.- \lambda\|\xi_3-\widehat{\boldsymbol{\xi}}_3\|^p\right\}\right] - \lambda\|\xi_2-\widehat{\boldsymbol{\xi}}_2\|^p\right\}\right],\end{aligned} \tag{EC.11}$$

where we have used the reformulation for soft penalized Wasserstein DRO (Zhang et al. 2022) in the first equality, where Assumption 1 in Zhang et al. (2022) is satisfied due to Assumption 3. Together with (EC.10), this shows $\mathbb{E}_{\widehat{\mathbb{P}}_t}[Q_t(x_{t-1},\widehat{\boldsymbol{\xi}}_t)] \geq \mathcal{Q}_t(x_{t-1},\xi_{[t-1]})$. To show that the above inequality holds as equality, by definition of $Q_t$, we have that

$$\begin{aligned}&\mathbb{E}_{\widehat{\mathbb{P}}_t}\left[Q_t(x_{t-1},\widehat{\boldsymbol{\xi}}_t)\right]\\ &= \mathbb{E}_{\widehat{\mathbb{P}}_t}\left[\sup_{\xi_t\in\Xi_t}\inf_{x_t\in\mathcal{X}_t(x_{t-1},\xi_t)}\left\{c_t^\top x_t + \mathbb{E}_{\widehat{\mathbb{P}}_{t+1}}\left[Q_{t+1}(x_t,\widehat{\boldsymbol{\xi}}_{t+1})\right] - \lambda\|\xi_t-\widehat{\boldsymbol{\xi}}_t\|^p\right\}\right]\\ &= \sup_{\mathbb{P}_t\in\mathcal{P}(\Xi_t)}\left\{\mathbb{E}_{\mathbb{P}_t}\left[\inf_{x_t\in\mathcal{X}_t(x_{t-1},\xi_t)}\left\{c_t^\top x_t + \mathbb{E}_{\widehat{\mathbb{P}}_{t+1}}\left[Q_{t+1}(x_t,\widehat{\boldsymbol{\xi}}_{t+1})\right]\right\}\right] - \lambda\mathcal{W}_p^p(\widehat{\mathbb{P}}_t,\mathbb{P}_t)\right\}\\ &= \sup_{\mathbb{P}_t\in\mathcal{P}(\Xi_t)}\inf_{\boldsymbol{x}_t(\xi_{[t-1]},\cdot)\in\mathcal{X}_t(x_{t-1},\cdot)}\left\{\mathbb{E}_{\mathbb{P}_t}\left[c_t^\top \boldsymbol{x}_t(\xi_{[t-1]},\boldsymbol{\xi}_t) + \mathbb{E}_{\widehat{\mathbb{P}}_{t+1}}\left[Q_{t+1}\big(\boldsymbol{x}_t(\xi_{[t-1]},\boldsymbol{\xi}_t),\widehat{\boldsymbol{\xi}}_{t+1}\big)\right]\right] - \lambda\mathcal{W}_p^p(\widehat{\mathbb{P}}_t,\mathbb{P}_t)\right\},\end{aligned}$$

where the second equality follows from the Wasserstein DRO reformulation (Zhang et al. 2022), and the third equality follows from the interchangeability principle (Shapiro et al. 2021, Theorem 9.108). Exchanging sup and min, we obtain

$$
\begin{aligned}
&\mathbb{E}_{\widehat{\mathbb{P}}_t}\left[Q_t(x_{t-1},\widehat{\boldsymbol{\xi}}_t)\right] \\
&\le \inf_{\boldsymbol{x}_t(\xi_{[t-1]},\cdot)\in\mathcal{X}_t(x_{t-1},\cdot)} \sup_{\mathbb{P}_t\in\mathcal{P}(\Xi_t)} \left\{\mathbb{E}_{\mathbb{P}_t}\left[c_t^\top \boldsymbol{x}_t(\xi_{[t-1]},\boldsymbol{\xi}_t)+\mathbb{E}_{\widehat{\mathbb{P}}_{t+1}}\left[Q_{t+1}\left(\boldsymbol{x}_t(\xi_{[t-1]},\boldsymbol{\xi}_t),\widehat{\boldsymbol{\xi}}_{t+1}\right)\right]\right]-\lambda\mathcal{W}_p^p(\widehat{\mathbb{P}}_t,\mathbb{P}_t)\right\} \\
&= Q_t(x_{t-1},\xi_{[t-1]}).
\end{aligned}
$$

This completes the proof of (EC.9). □

### EC.3.3. Proof of Proposition EC.5

Set

$$F(x_t,\xi_t)=c_t^\top x_t+\mathbb{E}_{\widehat{\mathbb{P}}_{t+1}}\left[Q_{t+1}(x_t,\widehat{\boldsymbol{\xi}}_{t+1})\right].$$

By Assumption 3, $F$ is bounded from below. In addition, the constraint set $\mathcal{X}_t(x_{t-1},\xi_t)$ satisfies the relaxed Slater condition. Thereby, using convex programming duality,

$$
\begin{aligned}
V_t^\star(x_{t-1},\xi_t) &= \inf_{x_t\in\mathcal{X}_t(x_{t-1},\xi_t)} F(x_t,\xi_t) \\
&= \inf_{x_t}\left\{F(x_t,\xi_t)\,:\, x_t\ge 0,\; A_t x_t=b_t-B_t x_{t-1}\right\} \\
&= \sup_{\alpha\ge 0,\beta}\inf_{x_t}\left\{F(x_t,\xi_t)-\alpha^\top x_t+\beta^\top\left(A_t x_t-b_t+B_t x_{t-1}\right)\right\} \\
&=: \sup_{\alpha\ge 0,\beta} G_{\alpha,\beta}(\xi_t).
\end{aligned}
$$

Observe that $F$ is concave and finite-valued in $\xi_t$ and $G_{\alpha,\beta}(\xi_t)$ is expressed as the infimum of concave and finite-valued functions of $\xi_t$, hence $G_{\alpha,\beta}(\xi_t)$ is also concave and thus lower semicontinuous in $\xi_t$. Since the pointwise supremum preserves lower semicontinuity, $V_t^\star(x_{t-1},\xi_t)$ is also lower semicontinuous in $\xi_t$. □

*Proof of Theorem 3.* We first consider $p=\infty$. We prove the result by showing that

$$\inf_{\boldsymbol{x}\in\mathcal{X}}\sup_{\mathbb{P}\in\mathfrak{M}}\mathbb{E}_{\mathbb{P}}\left[\sum_{t\in[T]}\boldsymbol{c}_t^\top\boldsymbol{x}_t(\boldsymbol{\xi}_{[t]})\right]=\inf_{\boldsymbol{x}\in\mathcal{X}}\sup_{\mathbb{P}\in\mathfrak{M}^{\mathsf{C}}}\mathbb{E}_{\mathbb{P}}\left[\sum_{t\in[T]}\boldsymbol{c}_t^\top\boldsymbol{x}_t(\boldsymbol{\xi}_{[t]})\right].$$

Since $\mathfrak{M}\subset\mathfrak{M}^{\mathsf{C}}$, the left-hand side is less than or equal to the right-hand side. It remains to prove the other direction. Using Proposition EC.3, the right-hand side above equals

$$
\begin{aligned}
&\inf_{x_1\in\mathcal{X}_1} c_1^\top x_1+\mathbb{E}_{\widehat{\mathbb{P}}_2}\Bigg[\sup_{\xi_2\in\Xi_2:\|\widehat{\xi}_2-\xi_2\|\le\vartheta}\;\inf_{x_2\in\mathcal{X}_2(x_1,\xi_2)}\Bigg\{c_2^\top x_2(\xi_2)+\mathbb{E}_{\widehat{\mathbb{P}}_3}\Bigg[\sup_{\xi_3\in\Xi_3:\|\widehat{\xi}_3-\xi_3\|\le\vartheta}\;\inf_{x_3\in\mathcal{X}_3(x_2,\xi_3)}\Bigg\{ \\
&\qquad c_3^\top x_3(\xi_{[3]})+\cdots+\mathbb{E}_{\widehat{\mathbb{P}}_T}\Bigg[\sup_{\xi_T\in\Xi_T:\|\widehat{\xi}_T-\xi_T\|\le\vartheta}\;\inf_{x_T\in\mathcal{X}_T(x_{T-1},\xi_T)} c_T^\top x_T(\xi_{[T]})\Bigg]\cdots\Bigg\}\Bigg]\Bigg\}\Bigg].
\end{aligned}
$$

For any $\varepsilon>0$, there exists a feasible policy $\boldsymbol{x}^\varepsilon=(\boldsymbol{x}_1^\varepsilon,\boldsymbol{x}_2^\varepsilon,\ldots,\boldsymbol{x}_T^\varepsilon)$ such that

$$\sup_{\mathbb{P}\in\mathfrak{M}}\mathbb{E}_{\mathbb{P}}\left[\sum_{t\in[T]}\boldsymbol{c}_t^\top\boldsymbol{x}_t^\varepsilon(\boldsymbol{\xi}_{[t]})\right]\le\inf_{\boldsymbol{x}\in\mathcal{X}}\sup_{\mathbb{P}\in\mathfrak{M}}\mathbb{E}_{\mathbb{P}}\left[\sum_{t\in[T]}\boldsymbol{c}_t^\top\boldsymbol{x}_t(\boldsymbol{\xi}_{[t]})\right]+\varepsilon.$$

Denote $f_t(x_t,\xi_t)=c_t^\top x_t$. Let $\mathbb{T}_1^{\epsilon_1}$ be the identity transport map. For $t=2,\ldots,T$, for any $\epsilon_t>0$ that is sufficiently small, there exists a transport map $\mathbb{T}_t^{\epsilon_t}:\widehat{\Xi}_{[t]}\to\Xi_t$, such that $\|\mathbb{T}_t^{\epsilon_t}(\widehat{\boldsymbol{\xi}}_{[t]})-\widehat{\boldsymbol{\xi}}_t\|\le\vartheta$, and

$$
\begin{aligned}
&f_t(\boldsymbol{x}_t^{\varepsilon}(\mathbb{T}_t^{\epsilon}(\widehat{\xi}_{[t]})),\mathbb{T}_t^{\epsilon}(\widehat{\xi}_{[t]}))+\mathbb{E}_{\widehat{\mathbb{P}}_{t+1}}\Bigg[\sup_{\xi_{t+1}\in\Xi_{t+1}:\|\widehat{\xi}_{t+1}-\xi_{t+1}\|\le\vartheta}\ \inf_{x_{t+1}\in\mathcal{X}_{t+1}(x_t,\xi_{t+1})}\Bigg\{\\
&\qquad f_{t+1}\Big(x_{t+1},\xi_{t+1}\Big)+\cdots+\mathbb{E}_{\widehat{\mathbb{P}}_T}\Bigg[\sup_{\xi_T\in\Xi_T:\|\widehat{\xi}_T-\xi_T\|\le\vartheta}\ \inf_{x_T\in\mathcal{X}_T(x_{T-1},\xi_T)} f_T\Big(x_T,\xi_T\Big)\Bigg]\cdots\Bigg\}\Bigg]\\
&\ge\sup_{\xi_t\in\Xi_t:\|\widehat{\xi}_t-\xi_t\|\le\vartheta} V_t^\star\Big(\boldsymbol{x}_{t-1}^{\varepsilon}(\mathbb{T}_{[t-1]}^{\epsilon}(\widehat{\boldsymbol{\xi}}_{[t-1]})),\xi_t\Big)-\epsilon_t.
\end{aligned}
$$

Let $\mathbb{T}^\epsilon=(\mathbb{T}_1^{\epsilon_1},\ldots,\mathbb{T}_T^{\epsilon_T})$ and define $\mathbb{P}^\epsilon:=\mathbb{T}_\#^\epsilon\widehat{\mathbb{P}}$. Using the construction in the proof of Theorem 1, $\mathbb{T}_t^{\epsilon_t}$ can be chosen to be a bi-causal transport map. It follows that $\mathbb{P}^\epsilon$ is a feasible distribution, and

$$
\sup_{\mathbb{P}\in\mathfrak{M}}\mathbb{E}_{\mathbb{P}}\left[\sum_{t\in[T]}\boldsymbol{c}_t^\top\boldsymbol{x}_t^{\varepsilon}(\boldsymbol{\xi}_{[t]})\right]\ge\mathbb{E}_{\mathbb{P}^\epsilon}\left[\sum_{t\in[T]}\boldsymbol{c}_t^\top\boldsymbol{x}_t^{\varepsilon}(\boldsymbol{\xi}_{[t]})\right]=\sum_{t\in[T]}\mathbb{E}_{\widehat{\mathbb{P}}}\Big[f_t(\boldsymbol{x}_t^{\varepsilon},\mathbb{T}_t^{\epsilon}(\widehat{\boldsymbol{\xi}}_{[t]}))\Big]\ge V_1^\star-\sum_{t\in[T]}\epsilon_T.
$$

Letting $\varepsilon\downarrow 0$ and $\epsilon_2,\ldots,\epsilon_T\downarrow 0$ yields the desired result.

Next, we consider $p=[1,\infty)$. We prove the result by showing that

$$
\inf_{\boldsymbol{x}\in\mathcal{X}}\sup_{\mathbb{P}\in\mathcal{P}(\Xi_{[T]})}\left\{\mathbb{E}_{\mathbb{P}}\left[\sum_{t\in[T]}\boldsymbol{c}_t^\top\boldsymbol{x}_t(\boldsymbol{\xi}_{[t]})\right]-\lambda\mathsf{D}_p^p(\widehat{\mathbb{P}},\mathbb{P})\right\}=\inf_{\boldsymbol{x}\in\mathcal{X}}\sup_{\mathbb{P}\in\mathcal{P}(\Xi_{[T]})}\left\{\mathbb{E}_{\mathbb{P}}\left[\sum_{t\in[T]}\boldsymbol{c}_t^\top\boldsymbol{x}_t(\boldsymbol{\xi}_{[t]})\right]-\lambda\mathsf{C}_p^p(\widehat{\mathbb{P}},\mathbb{P})\right\}.
$$

Since $\mathfrak{M}\subset\mathfrak{M}^{\mathsf{C}}$, the left-hand side is less than or equal to the right-hand side. It remains to prove the other direction. Using Proposition EC.3, the right-hand side above equals

$$
\begin{aligned}
&\inf_{x_1\in\mathcal{X}_1} c_1^\top x_1+\mathbb{E}_{\widehat{\mathbb{P}}_2}\Bigg[\sup_{\xi_2\in\Xi_2}\inf_{x_2\in\mathcal{X}_2(x_1,\xi_2)}\Bigg\{c_2^\top x_2(\xi_2)+\mathbb{E}_{\widehat{\mathbb{P}}_3}\Bigg[\sup_{\xi_3\in\Xi_3}\inf_{x_3\in\mathcal{X}_3(x_2,\xi_3)}\Bigg\{c_3^\top x_3(\xi_{[3]})+\cdots+\\
&\quad\mathbb{E}_{\widehat{\mathbb{P}}_T}\Bigg[\sup_{\xi_T\in\Xi_T}\inf_{x_T\in\mathcal{X}_T(x_{T-1},\xi_T)}\Big\{c_T^\top x_T(\xi_{[T]})-\lambda\|\xi_T-\widehat{\boldsymbol{\xi}}_T\|^p\Big\}\Bigg]\cdots-\lambda\|\xi_3-\widehat{\boldsymbol{\xi}}_3\|^p\Bigg\}\Bigg]-\lambda\|\xi_2-\widehat{\boldsymbol{\xi}}_2\|^p\Bigg\}\Bigg].
\end{aligned}
$$

For any $\varepsilon>0$, there exists a feasible policy $\boldsymbol{x}^\varepsilon=(\boldsymbol{x}_1^\varepsilon,\boldsymbol{x}_2^\varepsilon,\ldots,\boldsymbol{x}_T^\varepsilon)$ such that

$$
\sup_{\mathbb{P}\in\mathcal{P}(\Xi_{[T]})}\left\{\mathbb{E}_{\mathbb{P}}\left[\sum_{t\in[T]}\boldsymbol{c}_t^\top\boldsymbol{x}_t^{\varepsilon}(\boldsymbol{\xi}_{[t]})\right]-\lambda\mathsf{C}_p^p(\widehat{\mathbb{P}},\mathbb{P})\right\}\le\inf_{\boldsymbol{x}\in\mathcal{X}}\sup_{\mathbb{P}\in\mathcal{P}(\Xi_{[T]})}\left\{\mathbb{E}_{\mathbb{P}}\left[\sum_{t\in[T]}\boldsymbol{c}_t^\top\boldsymbol{x}_t(\boldsymbol{\xi}_{[t]})\right]-\lambda\mathsf{C}_p^p(\widehat{\mathbb{P}},\mathbb{P})\right\}+\varepsilon.
$$

Denote $f_t(x_t,\xi_t)=c_t^\top x_t$. Let $\mathbb{T}_1^{\epsilon_1}$ be the identity transport map. For $t=2,\ldots,T$, for any $\epsilon_t>0$ that is sufficiently small, there exists a transport map $\mathbb{T}_t^{\epsilon_t}:\widehat{\Xi}_{[t]}\to\Xi_t$, such that

$$
\begin{aligned}
&f_t(\boldsymbol{x}_t^{\varepsilon}(\mathbb{T}_t^{\epsilon}(\widehat{\xi}_{[t]})),\mathbb{T}_t^{\epsilon}(\widehat{\xi}_{[t]}))+\mathbb{E}_{\widehat{\mathbb{P}}_{t+1}}\Bigg[\sup_{\xi_{t+1}\in\Xi_{t+1}}\inf_{x_{t+1}\in\mathcal{X}_{t+1}(x_t,\xi_{t+1})}\Bigg\{f_{t+1}\Big(x_{t+1},\xi_{t+1}\Big)+\cdots+\\
&\qquad\mathbb{E}_{\widehat{\mathbb{P}}_T}\Bigg[\sup_{\xi_T\in\Xi_T}\inf_{x_T\in\mathcal{X}_T(x_{T-1},\xi_T)}\Big\{f_T\Big(x_T,\xi_T\Big)\cdots-\lambda\|\xi_T-\widehat{\boldsymbol{\xi}}_T\|^p\Big\}\Bigg]-\lambda\|\xi_{t+1}-\widehat{\boldsymbol{\xi}}_{t+1}\|^p\Bigg\}\Bigg]\\
&\ge\sup_{\xi_t\in\Xi_t}\Big\{V_t^\star\Big(\boldsymbol{x}_{t-1}^{\varepsilon}(\mathbb{T}_{[t-1]}^{\epsilon}(\widehat{\boldsymbol{\xi}}_{[t-1]})),\xi_t\Big)-\lambda\|\xi_t-\widehat{\boldsymbol{\xi}}_t\|^p\Big\}-\epsilon_t.
\end{aligned}
$$

Let $\mathbb{T}^{\epsilon} = (\mathbb{T}_1^{\epsilon_1}, \ldots, \mathbb{T}_T^{\epsilon_T})$ and define $\mathbb{P}^{\epsilon} := \mathbb{T}^{\epsilon}_{\#}\widehat{\mathbb{P}}$. Using the construction in the proof of Theorem 2, $\mathbb{T}_t^{\epsilon_t}$ can be chosen to be a bi-causal transport map. It follows that $\mathbb{P}^{\epsilon}$ is a feasible distribution, and

$$\begin{aligned}
&\sup_{\mathbb{P}\in\mathcal{P}(\Xi_{[T]})} \left\{\mathbb{E}_{\mathbb{P}}\left[\sum_{t\in[T]} \boldsymbol{c}_t^\top \boldsymbol{x}_t^{\varepsilon}(\boldsymbol{\xi}_{[t]})\right] - \lambda \mathrm{C}_p^p(\widehat{\mathbb{P}},\mathbb{P})\right\} \\
&\geq \mathbb{E}_{\mathbb{P}^{\epsilon}}\left[\sum_{t\in[T]} \boldsymbol{c}_t^\top \boldsymbol{x}_t^{\varepsilon}(\boldsymbol{\xi}_{[t]})\right] - \lambda \mathrm{C}_p^p(\widehat{\mathbb{P}},\mathbb{P}^{\epsilon}) \\
&= \sum_{t\in[T]} \mathbb{E}_{\widehat{\mathbb{P}}}\left[f_t(\boldsymbol{x}_t^{\varepsilon}, \mathbb{T}_t^{\epsilon}(\widehat{\boldsymbol{\xi}}_{[t]}))\right] - \lambda \mathrm{C}_p^p(\widehat{\mathbb{P}},\mathbb{P}^{\epsilon}) \\
&\geq V_1^\star - \sum_{t\in[T]} \epsilon_T.
\end{aligned}$$

Letting $\varepsilon \downarrow 0$ and $\epsilon_2, \ldots, \epsilon_T \downarrow 0$ yields the desired result. □

## EC.4. Proof for Proposition 2

By Theorem 3, we have

$$\begin{aligned}
&\inf_{\boldsymbol{x}\in\mathcal{X}} \sup_{\mathbb{P}\in\mathfrak{M}} \mathbb{E}_{\mathbb{P}}\left[\sum_{t\in[T]} \boldsymbol{c}_t^\top \boldsymbol{x}_t(\boldsymbol{\xi}_{[t]})\right] \\
&= \inf_{x_1\in\mathcal{X}_1, \lambda\geq 0} \lambda\vartheta^p + c_1^\top x_1 + \mathbb{E}_{\widehat{\mathbb{P}}_2}\Bigg[\sup_{\xi_2\in\Xi_2} \inf_{x_2\in\mathcal{X}_2(x_1,\xi_2)} \Bigg\{c_2^\top x_2(\xi_2) + \mathbb{E}_{\widehat{\mathbb{P}}_3}\Bigg[\sup_{\xi_3\in\Xi_3} \inf_{x_3\in\mathcal{X}_3(x_2,\xi_3)} \Bigg\{c_3^\top x_3(\xi_{[3]}) + \cdots + \\
&\qquad \mathbb{E}_{\widehat{\mathbb{P}}_T}\Bigg[\sup_{\xi_T\in\Xi_T} \inf_{x_T\in\mathcal{X}_T(x_{T-1},\xi_T)} \Big\{c_T^\top x_T(\xi_{[T]}) - \lambda\|\xi_T - \widehat{\boldsymbol{\xi}}_T\|^p\Big\}\Bigg] \cdots - \lambda\|\xi_3 - \widehat{\boldsymbol{\xi}}_3\|^p\Bigg\}\Bigg] - \lambda\|\xi_2 - \widehat{\boldsymbol{\xi}}_2\|^p\Bigg\}\Bigg].
\end{aligned}$$

Let us derive a lower bound of the problem above by allowing $\lambda$ to be different at each stage. Specifically, introducing auxiliary variables $\lambda_2 = \cdots = \lambda_T$, the problem above is equivalent to

$$\begin{aligned}
&\inf_{\substack{x_1\in\mathcal{X}_1 \\ \lambda_2=\cdots=\lambda_T\geq 0}} \frac{\vartheta^p}{T-1}\sum_{t=2}^{T} \lambda_t + c_1^\top x_1 + \mathbb{E}_{\widehat{\mathbb{P}}_2}\Bigg[\sup_{\xi_2\in\Xi_2} \inf_{x_2\in\mathcal{X}_2(x_1,\xi_2)} \Bigg\{c_2^\top x_2(\xi_2) + \mathbb{E}_{\widehat{\mathbb{P}}_3}\Bigg[\sup_{\xi_3\in\Xi_3} \inf_{x_3\in\mathcal{X}_3(x_2,\xi_3)} \Bigg\{c_3^\top x_3(\xi_{[3]}) + \cdots + \\
&\qquad \mathbb{E}_{\widehat{\mathbb{P}}_T}\Bigg[\sup_{\xi_T\in\Xi_T} \inf_{x_T\in\mathcal{X}_T(x_{T-1},\xi_T)} \Big\{c_T^\top x_T(\xi_{[T]}) - \lambda_T\|\xi_T - \widehat{\boldsymbol{\xi}}_T\|^p\Big\}\Bigg] \cdots - \lambda_3\|\xi_3 - \widehat{\boldsymbol{\xi}}_3\|^p\Bigg\}\Bigg] - \lambda_2\|\xi_2 - \widehat{\boldsymbol{\xi}}_2\|^p\Bigg\}\Bigg].
\end{aligned}$$

Dropping the alignment constraints $\lambda_2 = \cdots = \lambda_T$ and setting $\vartheta_t = (\frac{\vartheta^p}{T-1})^{1/p}$, we arrive at a lower bound

$$\begin{aligned}
&\inf_{\substack{x_1\in\mathcal{X}_1, \\ \lambda_2,\ldots\lambda_T\geq 0}} \sum_{t=2}^{T} \lambda_t\vartheta_t^p + c_1^\top x_1 + \mathbb{E}_{\widehat{\mathbb{P}}_2}\Bigg[\sup_{\xi_2\in\Xi_2} \inf_{x_2\in\mathcal{X}_2(x_1,\xi_2)} \Bigg\{c_2^\top x_2(\xi_2) + \mathbb{E}_{\widehat{\mathbb{P}}_3}\Bigg[\sup_{\xi_3\in\Xi_3} \inf_{x_3\in\mathcal{X}_3(x_2,\xi_3)} \Bigg\{c_3^\top x_3(\xi_{[3]}) + \cdots + \\
&\qquad \mathbb{E}_{\widehat{\mathbb{P}}_T}\Bigg[\sup_{\xi_T\in\Xi_T} \inf_{x_T\in\mathcal{X}_T(x_{T-1},\xi_T)} \Big\{c_T^\top x_T(\xi_{[T]}) - \lambda_T\|\xi_T - \widehat{\boldsymbol{\xi}}_T\|^p\Big\}\Bigg] \cdots - \lambda_3\|\xi_3 - \widehat{\boldsymbol{\xi}}_3\|^p\Bigg\}\Bigg] - \lambda_2\|\xi_2 - \widehat{\boldsymbol{\xi}}_2\|^p\Bigg\}\Bigg] \\
&= \inf_{x_1\in\mathcal{X}_1} c_1^\top x_1 + \sup_{\mathbb{P}_2:\mathcal{W}(\widehat{\mathbb{P}}_2,\mathbb{P}_2)\leq\vartheta_2} \mathbb{E}_{\mathbb{P}_2}\Bigg[\inf_{x_2\in\mathcal{X}_2(x_1,\boldsymbol{\xi}_2)} c_2^\top x_2 + \\
&\qquad \sup_{\mathbb{P}_3:\mathcal{W}(\widehat{\mathbb{P}}_3,\mathbb{P}_3)\leq\vartheta_3} \mathbb{E}_{\mathbb{P}_3}\Bigg[\inf_{x_3\in\mathcal{X}_3(x_2,\boldsymbol{\xi}_3)} c_3^\top x_3 + \ldots + \sup_{\mathbb{P}_T:\mathcal{W}(\widehat{\mathbb{P}}_T,\mathbb{P}_T)\leq\vartheta_T} \mathbb{E}_{\mathbb{P}_T}\Bigg[\inf_{x_T\in\mathcal{X}_T(x_{T-1},\boldsymbol{\xi}_T)} c_T^\top x_T\Bigg]\ldots\Bigg]\Bigg],
\end{aligned}$$

where the last equality uses duality result for Wasserstein DRO (Zhang et al. 2022).

Next, applying finite sample guarantee result for Wasserstein DRO (Kuhn et al. 2019) for each stage independently, when $N \geq \frac{\log(C/\beta)}{c}$, by setting

$$\vartheta_{t,N}(\beta) = \left(\frac{\log(C/\beta)}{cN}\right)^{\min(\frac{p}{d},\frac{1}{2})},$$

we have with probability at least $(1-\beta)^{T-1}$, it holds that

$$\inf_{\boldsymbol{x}\in\mathcal{X}} \mathbb{E}_{\mathbb{P}^*}\left[\sum_{t\in[T]} \boldsymbol{c}_t^\top \boldsymbol{x}_t(\boldsymbol{\xi}_{[t]})\right] \leq \inf_{x_1\in\mathcal{X}_1} c_1^\top x_1 + \sup_{\mathbb{P}_2:\mathcal{W}(\widehat{\mathbb{P}}_2,\mathbb{P}_2)\leq\vartheta_2} \mathbb{E}_{\mathbb{P}_2}\Bigg[\inf_{x_2\in\mathcal{X}_2(x_1,\boldsymbol{\xi}_2)} c_2^\top x_2 + \sup_{\mathbb{P}_3:\mathcal{W}(\widehat{\mathbb{P}}_3,\mathbb{P}_3)\leq\vartheta_3} \mathbb{E}_{\mathbb{P}_3}\Bigg[$$
$$\inf_{x_3\in\mathcal{X}_3(x_2,\boldsymbol{\xi}_3)} c_3^\top x_3 + \cdots + \sup_{\mathbb{P}_T:\mathcal{W}(\widehat{\mathbb{P}}_T,\mathbb{P}_T)\leq\vartheta_T} \mathbb{E}_{\mathbb{P}_T}\left[\inf_{x_T\in\mathcal{X}_T(x_{T-1},\boldsymbol{\xi}_T)} c_T^\top x_T\right]\ldots\Bigg]\Bigg],$$

Note that the probability $(1-\beta)^{T-1}$ results from the fact that for $\mathbb{P}^*$ to be in the nested ambiguity set, it is sufficient to have each marginal distribution $\mathbb{P}_t^*$ to be in the Wasserstein ambiguity set at each stage, each of which has $1-\beta$ probability by the finite sample guarantee result for Wasserstein DRO. Using the elementary inequality $(1-\beta)^{T-1} \geq 1-(T-1)\beta$, we replace $\beta$ with $\beta/(T-1)$ in the expression of $\vartheta_{t,N}(\beta)$ yields that the above holds with probability at least $1-\beta$ by setting

$$\vartheta_{t,N}(\beta) = \left(\frac{\log(C(T-1)/\beta)}{cN}\right)^{\min(\frac{p}{d},\frac{1}{2})},$$

Finally, setting $\vartheta_N(\beta) = \sup_{t=2,\ldots,T}\left((T-1)\vartheta_{t,N}(\beta)^p\right)^{1/p}$ yields the result. □

## EC.5. Proofs for Section 4

### EC.5.1. Proof of Corollary 1

Applying Theorem 3, we have that for $p=\infty$,

$$\begin{aligned}
Q_t(x_{t-1},\widehat{\boldsymbol{c}}_t) &= \sup_{\|c_t-\widehat{\boldsymbol{c}}_t\|\leq\vartheta} \inf_{x_t\in\mathcal{X}_t(x_{t-1})} \left\{c_t^\top x_t + \mathbb{E}_{\widehat{\mathbb{P}}_{t+1}}\left[Q_{t+1}(x_t,\widehat{\boldsymbol{c}}_{t+1})\right]\right\} \\
&= \inf_{x_t\in\mathcal{X}_t(x_{t-1})} \sup_{\|c_t-\widehat{\boldsymbol{c}}_t\|\leq\vartheta} \left\{c_t^\top x_t + \mathbb{E}_{\widehat{\mathbb{P}}_{t+1}}\left[Q_{t+1}(x_t,\widehat{\boldsymbol{c}}_{t+1})\right]\right\} \\
&= \inf_{x_t\in\mathcal{X}_t(x_{t-1})} \left\{\widehat{\boldsymbol{c}}_t^\top x_t + \vartheta\|x_t\|_* + \mathbb{E}_{\widehat{\mathbb{P}}_{t+1}}\left[Q_{t+1}(x_t,\widehat{\boldsymbol{c}}_{t+1})\right]\right\}.
\end{aligned}$$

The second equality is due to Sion's minimax theorem: the objective is convex and continuous in $x_t$ (Shapiro et al. 2021, Section 3.2.1), as well as linear in $c_t$. In addition, both of the feasible sets are convex with the feasible set of $c_t$ being compact.

For $p\in[1,\infty)$, define

$$\mathcal{L}(\delta^p) = \sup_{c_t\in\Xi_t:\|c_t-\widehat{\boldsymbol{c}}_t\|^p\leq\delta^p} \inf_{x_t\in\mathcal{X}_t(x_{t-1})} \left\{c_t^\top x_t + \mathbb{E}_{\widehat{\mathbb{P}}_{t+1}}\left[Q_{t+1}(x_t,\widehat{\boldsymbol{c}}_{t+1})\right]\right\}.$$

Observe that $\mathcal{L}(\delta^p)$ is a concave function. Consider the Legendre transform of $-\mathcal{L}$

$$\begin{aligned}
(-\mathcal{L})^*(-\lambda) &= \sup_{\delta\geq 0}\left\{(-\lambda)\delta^p - (-\mathcal{L})(\delta^p)\right\} \\
&= \sup_{\delta\geq 0} \sup_{c_t\in\Xi_t:\|c_t-\widehat{\boldsymbol{c}}_t\|^p\leq\delta^p} \inf_{x_t\in\mathcal{X}_t(x_{t-1})} \left\{-\lambda\delta^p + c_t^\top x_t + \mathbb{E}_{\widehat{\mathbb{P}}_{t+1}}\left[Q_{t+1}(x_t,\widehat{\boldsymbol{c}}_{t+1})\right]\right\} \\
&= \sup_{c_t\in\Xi_t} \inf_{x_t\in\mathcal{X}_t(x_{t-1})} \left\{-\lambda\|c_t-\widehat{\boldsymbol{c}}_t\|^p + c_t^\top x_t + \mathbb{E}_{\widehat{\mathbb{P}}_{t+1}}\left[Q_{t+1}(x_t,\widehat{\boldsymbol{c}}_{t+1})\right]\right\} \\
&= Q_t(x_{t-1},\widehat{\boldsymbol{c}}_t).
\end{aligned}$$

On the other hand, similar to the case of $p=\infty$, using Sion's minimax theorem,

$$\mathcal{L}(\delta^p)=\inf_{x_t\in\mathcal{X}_t(x_{t-1})}\ \sup_{c_t\in\Xi_t:\|c_t-\widehat{\boldsymbol{c}}_t\|^p\le\delta^p}\left\{c_t^\top x_t+\mathbb{E}_{\widehat{\mathbb{P}}_{t+1}}\left[Q_{t+1}(x_t,\widehat{\boldsymbol{c}}_{t+1})\right]\right\}.$$

It follows that

$$\begin{aligned}Q_t(x_{t-1},\widehat{\boldsymbol{c}}_t)=(-\mathcal{L})^*(-\lambda)&=\sup_{\delta\ge0}\left\{(-\lambda)\delta^p-(-\mathcal{L})(\delta^p)\right\}\\&=\sup_{\delta\ge0}\ \inf_{x_t\in\mathcal{X}_t(x_{t-1})}\ \sup_{c_t\in\Xi_t:\|c_t-\widehat{\boldsymbol{c}}_t\|^p\le\delta^p}\left\{-\lambda\delta^p+c_t^\top x_t+\mathbb{E}_{\widehat{\mathbb{P}}_{t+1}}\left[Q_{t+1}(x_t,\widehat{\boldsymbol{c}}_{t+1})\right]\right\}\\&=\sup_{\delta\ge0}\ \inf_{x_t\in\mathcal{X}_t(x_{t-1})}\left\{-\lambda\delta^p+\delta\|x_t\|_*+\widehat{\boldsymbol{c}}_t^\top x_t+\mathbb{E}_{\widehat{\mathbb{P}}_{t+1}}\left[Q_{t+1}(x_t,\widehat{\boldsymbol{c}}_{t+1})\right]\right\}.\end{aligned}$$

Define

$$h(\delta):=\inf_{x_t\in\mathcal{X}_t(x_{t-1})}\left\{\delta\|x_t\|_*+\widehat{\boldsymbol{c}}_t^\top x_t+\mathbb{E}_{\widehat{\mathbb{P}}_{t+1}}\left[Q_{t+1}(x_t,\widehat{\boldsymbol{c}}_{t+1})\right]\right\}.$$

Since $h$ is the pointwise infimum of a family of affine functions of $\delta$, it is concave and upper semicontinuous. This implies that $h$ is continuous, and particularly, $\lim_{\delta\downarrow0}h(\delta)=h(0)$. By Assumption 3, we have $h(0)>-\infty$. Define $H(\delta):=-\lambda\delta^p+h(\delta)$. Then $H(0)=h(0)$. With this definition, we have

$$Q_t(x_{t-1},\widehat{\boldsymbol{c}}_t)=\sup_{\delta\ge0}H(\delta).$$

We consider two separate cases, $Q_t(x_{t-1},\widehat{\boldsymbol{c}}_t)>h(0)$ and $Q_t(x_{t-1},\widehat{\boldsymbol{c}}_t)=h(0)$.

For the first case, we have that for all sufficiently small $\delta_0>0$, $Q_t(x_{t-1},\widehat{\boldsymbol{c}}_t)=\sup_{\delta\ge\delta_0}H(\delta)$. By Assumption 2, there exists some $x_t^0\in\mathcal{X}_t(x_{t-1})$. Define

$$D_0:=\|x_t^0\|_*+\frac{\widehat{\boldsymbol{c}}_t^\top x_t^0+\mathbb{E}_{\widehat{\mathbb{P}}_{t+1}}\left[Q_{t+1}(x_t^0,\widehat{\boldsymbol{c}}_{t+1})\right]-h(0)}{\delta_0},$$

Then it follows that for any $\delta\ge\delta_0$ and $x_t$ with $\|x_t\|_*>D_0$,

$$\begin{aligned}&\delta\|x_t\|_*+\widehat{\boldsymbol{c}}_t^\top x_t+\mathbb{E}_{\widehat{\mathbb{P}}_{t+1}}\left[Q_{t+1}(x_t,\widehat{\boldsymbol{c}}_{t+1})\right]\\&\ge\delta_0\|x_t\|_*+\widehat{\boldsymbol{c}}_t^\top x_t+\mathbb{E}_{\widehat{\mathbb{P}}_{t+1}}\left[Q_{t+1}(x_t,\widehat{\boldsymbol{c}}_{t+1})\right]\\&>\delta_0\|x_t^0\|_*+\widehat{\boldsymbol{c}}_t^\top x_t^0+\mathbb{E}_{\widehat{\mathbb{P}}_{t+1}}\left[Q_{t+1}(x_t^0,\widehat{\boldsymbol{c}}_{t+1})\right]-h(0)+\widehat{\boldsymbol{c}}_t^\top x_t+\mathbb{E}_{\widehat{\mathbb{P}}_{t+1}}\left[Q_{t+1}(x_t,\widehat{\boldsymbol{c}}_{t+1})\right]\\&\ge\delta_0\|x_t^0\|_*+\widehat{\boldsymbol{c}}_t^\top x_t^0+\mathbb{E}_{\widehat{\mathbb{P}}_{t+1}}\left[Q_{t+1}(x_t^0,\widehat{\boldsymbol{c}}_{t+1})\right]-h(0)+h(0)\\&\ge\delta_0\|x_t^0\|_*+\widehat{\boldsymbol{c}}_t^\top x_t^0+\mathbb{E}_{\widehat{\mathbb{P}}_{t+1}}\left[Q_{t+1}(x_t^0,\widehat{\boldsymbol{c}}_{t+1})\right],\end{aligned}$$

which shows that such $x_t$ cannot be a minimizer of the minimization problem defining $h(\delta)$. Therefore we have for any $D\ge D_0$,

$$\begin{aligned}Q_t(x_{t-1},\widehat{\boldsymbol{c}}_t)&=\sup_{\delta\ge\delta_0}H(\delta)\\&=\sup_{\delta\ge\delta_0}\ \inf_{x_t\in\mathcal{X}_t(x_{t-1}),\|x_t\|_*\le D}\left\{-\lambda\delta^p+\delta\|x_t\|_*+\widehat{\boldsymbol{c}}_t^\top x_t+\mathbb{E}_{\widehat{\mathbb{P}}_{t+1}}\left[Q_{t+1}(x_t,\widehat{\boldsymbol{c}}_{t+1})\right]\right\}.\end{aligned}$$

Letting $\delta_0\downarrow0$, it follows that

$$\begin{aligned}Q_t(x_{t-1},\widehat{\boldsymbol{c}}_t)&=\lim_{\delta_0\downarrow0}\ \sup_{\delta>\delta_0}\ \inf_{x_t\in\mathcal{X}_t(x_{t-1}),\|x_t\|_*\le D}\left\{-\lambda\delta^p+\delta\|x_t\|_*+\widehat{\boldsymbol{c}}_t^\top x_t+\mathbb{E}_{\widehat{\mathbb{P}}_{t+1}}\left[Q_{t+1}(x_t,\widehat{\boldsymbol{c}}_{t+1})\right]\right\}\\&=\sup_{\delta>0}\ \inf_{x_t\in\mathcal{X}_t(x_{t-1}),\|x_t\|_*\le D}\left\{-\lambda\delta^p+\delta\|x_t\|_*+\widehat{\boldsymbol{c}}_t^\top x_t+\mathbb{E}_{\widehat{\mathbb{P}}_{t+1}}\left[Q_{t+1}(x_t,\widehat{\boldsymbol{c}}_{t+1})\right]\right\}\\&=\inf_{x_t\in\mathcal{X}_t(x_{t-1}),\|x_t\|_*\le D}\ \sup_{\delta>0}\left\{-\lambda\delta^p+\delta\|x_t\|_*+\widehat{\boldsymbol{c}}_t^\top x_t+\mathbb{E}_{\widehat{\mathbb{P}}_{t+1}}\left[Q_{t+1}(x_t,\widehat{\boldsymbol{c}}_{t+1})\right]\right\},\end{aligned}$$

where the third equality follows from Sion's minimax theorem because of the boundedness of $x_t$. Viewing the inner supremum as a function of $D$, it is monotone decreasing in $D$ and bounded from below by $Q_t(x_{t-1},\widehat{\boldsymbol{c}}_t)$. Therefore, letting $D\to\infty$ yields that

$$Q_t(x_{t-1},\widehat{\boldsymbol{c}}_t)=\inf_{x_t\in\mathcal{X}_t(x_{t-1})}\sup_{\delta>0}\Big\{-\lambda\delta^p+\delta\|x_t\|_*+\widehat{\boldsymbol{c}}_t^\top x_t+\mathbb{E}_{\widehat{\mathbb{P}}_{t+1}}\big[Q_{t+1}(x_t,\widehat{\boldsymbol{c}}_{t+1})\big]\Big\}$$
$$=\begin{cases}\inf_{x_t\in\mathcal{X}_t(x_{t-1})}\Big\{\widehat{\boldsymbol{c}}_t^\top x_t+(1-1/p)(\frac{1}{p\lambda})^{\frac{1}{p-1}}\|x_t\|_*^{\frac{p}{p-1}}+\mathbb{E}_{\widehat{\mathbb{P}}_{t+1}}\big[Q_{t+1}(x_t,\widehat{\boldsymbol{c}}_{t+1})\big]\Big\}, & \text{if } p\in(1,\infty),\\ \inf_{x_t\in\mathcal{X}_t(x_{t-1}),\|x_t\|_*\le\lambda}\Big\{\widehat{\boldsymbol{c}}_t^\top x_t+\mathbb{E}_{\widehat{\mathbb{P}}_{t+1}}\big[Q_{t+1}(x_t,\widehat{\boldsymbol{c}}_{t+1})\big]\Big\}, & \text{if } p=1.\end{cases}$$

For the second case, it holds that

$$Q_t(x_{t-1},\widehat{\boldsymbol{c}}_t)=h(0)=\inf_{x_t\in\mathcal{X}_t(x_{t-1})}\Big\{\widehat{\boldsymbol{c}}_t^\top x_t+\mathbb{E}_{\widehat{\mathbb{P}}_{t+1}}\big[Q_{t+1}(x_t,\widehat{\boldsymbol{c}}_{t+1})\big]\Big\},$$

therefore we have for every $\delta>0$, $H(\delta)\le H(0)$, which implies that $h(\delta)\le h(0)+\lambda\delta^p$. Moreover, by Assumption 3, $h(0)>-\infty$. Consider the minimization problem defining $h(\delta)$, $\delta>0$. As $\|x_t\|_*\to\infty$, the first term goes to infinity, and the sum of the second and third terms is bounded from below. Hence, the minimizer is attained at some finite point, denoted as $x_t(\delta)$. Thus, we have

$$h(0)+\lambda\delta^p\ge h(\delta)=\delta\|x_t(\delta)\|_*+\widehat{\boldsymbol{c}}_t^\top x_t(\delta)+\mathbb{E}_{\widehat{\mathbb{P}}_{t+1}}\big[Q_{t+1}(x_t(\delta),\widehat{\boldsymbol{c}}_{t+1})\big]\ge\delta\|x_t(\delta)\|_*+h(0),$$

which implies

$$\|x_t(\delta)\|_*\le\lambda\delta^{p-1}. \tag{EC.12}$$

When $p\in(1,\infty)$, letting $\delta\downarrow 0$, $x_t(\delta)\to 0\in\mathcal{X}_t(x_{t-1})$. It follows that

$$Q_t(x_{t-1},\widehat{\boldsymbol{c}}_t)=h(0)=\lim_{\delta\downarrow 0}h(\delta)=\lim_{\delta\downarrow 0}\ \delta\|x_t(\delta)\|_*+\widehat{\boldsymbol{c}}_t^\top x_t(\delta)+\mathbb{E}_{\widehat{\mathbb{P}}_{t+1}}\big[Q_{t+1}(x_t(\delta),\widehat{\boldsymbol{c}}_{t+1})\big]=\mathbb{E}_{\widehat{\mathbb{P}}_{t+1}}\big[Q_{t+1}(0,\widehat{\boldsymbol{c}}_{t+1})\big].$$

Meanwhile,

$$h(0)\le\inf_{x_t\in\mathcal{X}_t(x_{t-1})}\Big\{\widehat{\boldsymbol{c}}_t^\top x_t+(1-1/p)(\frac{1}{p\lambda})^{\frac{1}{p-1}}\|x_t\|_*^{\frac{p}{p-1}}+\mathbb{E}_{\widehat{\mathbb{P}}_{t+1}}\big[Q_{t+1}(x_t,\widehat{\boldsymbol{c}}_{t+1})\big]\Big\}\le\mathbb{E}_{\widehat{\mathbb{P}}_{t+1}}\big[Q_{t+1}(0,\widehat{\boldsymbol{c}}_{t+1})\big]=h(0).$$

This verifies the expression for $Q_t(x_{t-1},\widehat{\boldsymbol{c}}_t)$. When $p=1$, (EC.12) means $\|x_t(\delta)\|_*\le\lambda$. It follows that

$$\begin{aligned}Q_t(x_{t-1},\widehat{\boldsymbol{c}}_t)=H(0)&=\lim_{\delta\downarrow 0}-\lambda\delta^p+h(\delta)\\&=\lim_{\delta\downarrow 0}-\lambda\delta^p+\delta\|x_t(\delta)\|_*+\widehat{\boldsymbol{c}}_t^\top x_t(\delta)+\mathbb{E}_{\widehat{\mathbb{P}}_{t+1}}\big[Q_{t+1}(x_t(\delta),\widehat{\boldsymbol{c}}_{t+1})\big]\\&\ge\inf_{x_t\in\mathcal{X}_t(x_{t-1}),\|x_t\|_*\le\lambda}\Big\{\widehat{\boldsymbol{c}}_t^\top x_t+\mathbb{E}_{\widehat{\mathbb{P}}_{t+1}}\big[Q_{t+1}(x_t,\widehat{\boldsymbol{c}}_{t+1})\big]\Big\}.\end{aligned}$$

But $H(0)=h(0)=\inf_{x_t\in\mathcal{X}_t(x_{t-1})}\big\{\widehat{\boldsymbol{c}}_t^\top x_t+\mathbb{E}_{\widehat{\mathbb{P}}_{t+1}}[Q_{t+1}(x_t,\widehat{\boldsymbol{c}}_{t+1})]\big\}$. So the above inequality holds as equality, which verifies the expression for $Q_t(x_{t-1},\widehat{\boldsymbol{c}}_t)$.

Combining the two cases completes the proof. □

### EC.5.2. Proof of Proposition 3

We first prove for the case of $p=\infty$. Consider

$$\sup_{\|c_t-\widehat{\boldsymbol{c}}_t\|\le\vartheta}\Big\{c_t^\top\bar{\boldsymbol{x}}_t(c_{[t]})+\mathbb{E}_{\widehat{\mathbb{P}}_{t+1}}\big[Q_{t+1}(\bar{\boldsymbol{x}}_t(c_{[t]}),\widehat{\boldsymbol{c}}_{t+1})\big]\Big\}.$$

By definition of $\bar{\boldsymbol{x}}_t$, we have

$$c_t^\top\bar{\boldsymbol{x}}_t(c_{[t]})+\mathbb{E}_{\widehat{\mathbb{P}}_{t+1}}\big[Q_{t+1}(\bar{\boldsymbol{x}}_t(c_{[t]}),\widehat{\boldsymbol{c}}_{[t]})\big]=c_t^\top\widehat{\boldsymbol{x}}_t(\widehat{\boldsymbol{c}}^c_{[t]})+\mathbb{E}_{\widehat{\mathbb{P}}_{t+1}}\big[Q_{t+1}(\widehat{\boldsymbol{x}}_t(\widehat{\boldsymbol{c}}^c_{[t]}),\widehat{\boldsymbol{c}}_{t+1})\big].$$

Decompose the right-hand side as

$$\begin{aligned}
&c_t^\top \widehat{\boldsymbol{x}}_t\big(\widehat{\boldsymbol{c}}^c_{[t]}\big) + \mathbb{E}_{\widehat{\mathbb{P}}_{t+1}}\big[Q_{t+1}\big(\widehat{\boldsymbol{x}}_t(\widehat{\boldsymbol{c}}^c_{[t]}), \widehat{\boldsymbol{c}}_{t+1}\big)\big] \\
&= (\widehat{\boldsymbol{c}}^c_t)^\top \widehat{\boldsymbol{x}}_t(\widehat{\boldsymbol{c}}^c_{[t]}) + (c_t - \widehat{\boldsymbol{c}}^c_t)^\top \widehat{\boldsymbol{x}}_t(\widehat{\boldsymbol{c}}^c_{[t]}) + \mathbb{E}_{\widehat{\mathbb{P}}_{t+1}}\big[Q_{t+1}\big(\widehat{\boldsymbol{x}}_t(\widehat{\boldsymbol{c}}^c_{[t]}), \widehat{\boldsymbol{c}}_{t+1}\big)\big] \\
&\le (\widehat{\boldsymbol{c}}^c_t)^\top \widehat{\boldsymbol{x}}_t(\widehat{\boldsymbol{c}}^c_{[t]}) + \vartheta \|\widehat{\boldsymbol{x}}_t(\widehat{\boldsymbol{c}}^c_{[t]})\|_* + \mathbb{E}_{\widehat{\mathbb{P}}_{t+1}}\big[Q_{t+1}\big(\widehat{\boldsymbol{x}}_t(\widehat{\boldsymbol{c}}^c_{[t]}), \widehat{\boldsymbol{c}}_{t+1}\big)\big].
\end{aligned}$$

By definition of $\widehat{\boldsymbol{c}}^c_{[t]}$, the last expression is upper bounded by

$$\widehat{\boldsymbol{c}}_t^\top \widehat{\boldsymbol{x}}_t(\widehat{\boldsymbol{c}}^c_{[t-1]}, \widehat{\boldsymbol{c}}_t) + \vartheta \|\widehat{\boldsymbol{x}}_t(\widehat{\boldsymbol{c}}^c_{[t-1]}, \widehat{\boldsymbol{c}}_t)\|_* + \mathbb{E}_{\widehat{\mathbb{P}}_{t+1}}\Big[Q_{t+1}(\widehat{\boldsymbol{x}}_t(\widehat{\boldsymbol{c}}^c_{[t-1]}, \widehat{\boldsymbol{c}}_t), \widehat{\boldsymbol{c}}_{t+1})\Big] = Q_t(\widehat{\boldsymbol{x}}_{t-1}(\widehat{\boldsymbol{c}}^c_{[t-1]}), \widehat{\boldsymbol{c}}_t).$$

Therefore, we have shown that

$$\sup_{\|c_t - \widehat{\boldsymbol{c}}_t\| \le \vartheta} \Big\{ c_t^\top \bar{\boldsymbol{x}}_t(c_{[t]}) + \mathbb{E}_{\widehat{\mathbb{P}}_{t+1}}\big[Q_{t+1}(\bar{\boldsymbol{x}}_t(c_{[t]}), \widehat{\boldsymbol{c}}_{t+1})\big] \Big\} \le Q_t(\bar{\boldsymbol{x}}_{t-1}(c_{[t-1]}), \widehat{\boldsymbol{c}}_t).$$

Taking the expectation over $\widehat{\boldsymbol{c}}_t \sim \widehat{\mathbb{P}}_t$ and using Theorem 3, we obtain that

$$\mathbb{E}_{\widehat{\mathbb{P}}_t}\left[ \sup_{\|c_t - \widehat{\boldsymbol{c}}_t\| \le \vartheta} \Big\{ c_t^\top \bar{\boldsymbol{x}}_t(c_{[t]}) + \mathbb{E}_{\widehat{\mathbb{P}}_{t+1}}\big[Q_{t+1}(\bar{\boldsymbol{x}}_t(c_{[t]}), \widehat{\boldsymbol{c}}_{t+1})\big] \Big\} \right] \le \mathcal{Q}_t(\bar{\boldsymbol{x}}_{t-1}(c_{[t-1]})).$$

Hence we have

$$\bar{\boldsymbol{x}}_t(c_{[t-1]}, \cdot) \in \operatorname*{arg\,min}_{\boldsymbol{x}_t(c_{[t-1]}, \cdot) \in \mathcal{X}_t(\bar{\boldsymbol{x}}_{t-1}(c_{[t-1]}), \cdot)} \rho_t\big[\boldsymbol{c}_t^\top \boldsymbol{x}_t(c_{[t-1]}, \boldsymbol{c}_t) + \mathcal{Q}_{t+1}\big(\boldsymbol{x}_t(c_{[t-1]}, \boldsymbol{c}_t)\big)\big],$$

which shows the optimality of the policy $(\bar{\boldsymbol{x}}_1, \ldots, \bar{\boldsymbol{x}}_T)$.

Next, we prove for the case of $p \in [1, \infty)$. Consider

$$\sup_{c_t \in \Xi_t} \Big\{ c_t^\top \bar{\boldsymbol{x}}_t(c_{[t]}) + \mathbb{E}_{\widehat{\mathbb{P}}_{t+1}}\big[Q_{t+1}(\bar{\boldsymbol{x}}_t(c_{[t]}), \widehat{\boldsymbol{c}}_{t+1})\big] - \lambda \|c_t - \widehat{\boldsymbol{c}}_t\|^p \Big\}.$$

Using the expression for $\bar{\boldsymbol{x}}_t$, we obtain

$$c_t^\top \bar{\boldsymbol{x}}_t(c_{[t]}) + \mathbb{E}_{\widehat{\mathbb{P}}_{t+1}}\big[Q_{t+1}\big(\bar{\boldsymbol{x}}_t(c_{[t]}), \widehat{\boldsymbol{c}}_{[t]}\big)\big] = c_t^\top \widehat{\boldsymbol{x}}_t\big(\widehat{\boldsymbol{c}}^c_{[t]}\big) + \mathbb{E}_{\widehat{\mathbb{P}}_{t+1}}\big[Q_{t+1}\big(\widehat{\boldsymbol{x}}_t(\widehat{\boldsymbol{c}}^c_{[t]}), \widehat{\boldsymbol{c}}_{t+1}\big)\big].$$

Decompose the right-hand side as

$$\begin{aligned}
&c_t^\top \widehat{\boldsymbol{x}}_t\big(\widehat{\boldsymbol{c}}^c_t\big) + \mathbb{E}_{\widehat{\mathbb{P}}_{t+1}}\big[Q_{t+1}\big(\widehat{\boldsymbol{x}}_t(\widehat{\boldsymbol{c}}^c_{[t]}), \widehat{\boldsymbol{c}}_{t+1}\big)\big] \\
&= (\widehat{\boldsymbol{c}}^c_{[t]})^\top \widehat{\boldsymbol{x}}_t(\widehat{\boldsymbol{c}}^c_{[t]}) + (c_t - \widehat{\boldsymbol{c}}^c_t)^\top \widehat{\boldsymbol{x}}_t(\widehat{\boldsymbol{c}}^c_{[t]}) + \mathbb{E}_{\widehat{\mathbb{P}}_{t+1}}\big[Q_{t+1}\big(\widehat{\boldsymbol{x}}_t(\widehat{\boldsymbol{c}}^c_{[t]}), \widehat{\boldsymbol{c}}_{t+1}\big)\big] \\
&\le (\widehat{\boldsymbol{c}}^c_t)^\top \widehat{\boldsymbol{x}}_t(\widehat{\boldsymbol{c}}^c_{[t]}) + \frac{1}{p}\big((p\lambda)^{1/p}\|c_t - \widehat{\boldsymbol{c}}^c_t\|\big)^p + (1 - \frac{1}{p})\big((p\lambda)^{-1/p}\|\widehat{\boldsymbol{x}}_t(\widehat{\boldsymbol{c}}^c_{[t]})\|_*\big)^{\frac{p}{p-1}} + \mathbb{E}_{\widehat{\mathbb{P}}_{t+1}}\big[Q_{t+1}\big(\widehat{\boldsymbol{x}}_t(\widehat{\boldsymbol{c}}^c_{[t]}), \widehat{\boldsymbol{c}}_{t+1}\big)\big],
\end{aligned}$$

where the last inequality follows from Young's inequality. By definition of $\widehat{\boldsymbol{c}}^c_{[t]}$, the last expression is upper bounded by

$$\widehat{\boldsymbol{c}}_t^\top \widehat{\boldsymbol{x}}_t(\widehat{\boldsymbol{c}}^c_{[t-1]}, \widehat{\boldsymbol{c}}_t) + \lambda\|c_t - \widehat{\boldsymbol{c}}_t\|^p + (1 - \frac{1}{p})\big((p\lambda)^{-1/p}\|\widehat{\boldsymbol{x}}_t(\widehat{\boldsymbol{c}}^c_{[t-1]}, \widehat{\boldsymbol{c}}_t)\|_*\big)^{\frac{p}{p-1}} + \mathbb{E}_{\widehat{\mathbb{P}}_{t+1}}\big[Q_{t+1}\big(\widehat{\boldsymbol{x}}_t(\widehat{\boldsymbol{c}}^c_{[t-1]}, \widehat{\boldsymbol{c}}_t), \widehat{\boldsymbol{c}}_{t+1}\big)\big].$$

Thereby we have

$$\begin{aligned}
&\sup_{c_t \in \Xi_t} \Big\{ c_t^\top \bar{\boldsymbol{x}}_t(c_{[t]}) + \mathbb{E}_{\widehat{\mathbb{P}}_{t+1}}\big[\widehat{Q}_{t+1}(\bar{\boldsymbol{x}}_t(c_{[t]}), \widehat{\boldsymbol{c}}_{t+1})\big] - \lambda\|c_t - \widehat{\boldsymbol{c}}_t\|^p \Big\} \\
&\le \widehat{\boldsymbol{c}}_t^\top \widehat{\boldsymbol{x}}_t(\widehat{\boldsymbol{c}}^c_{[t-1]}, \widehat{\boldsymbol{c}}_t) + (1 - \frac{1}{p})\big((p\lambda)^{-1/p}\|\widehat{\boldsymbol{x}}_t(\widehat{\boldsymbol{c}}^c_{[t-1]}, \widehat{\boldsymbol{c}}_t)\|_*\big)^{\frac{p}{p-1}} + \mathbb{E}_{\widehat{\mathbb{P}}_{t+1}}\big[Q_{t+1}\big(\widehat{\boldsymbol{x}}_t(\widehat{\boldsymbol{c}}^c_{[t-1]}, \widehat{\boldsymbol{c}}_t), \widehat{\boldsymbol{c}}_{t+1}\big)\big]. \\
&= Q_t(\widehat{\boldsymbol{x}}_{t-1}(\widehat{\boldsymbol{c}}^c_{[t-1]}), \widehat{\boldsymbol{c}}_t) \\
&= Q_t(\bar{\boldsymbol{x}}_{t-1}(c_{[t-1]}), \widehat{\boldsymbol{c}}_t),
\end{aligned}$$

which shows the optimality of the policy $(\bar{\boldsymbol{x}}_1, \ldots, \bar{\boldsymbol{x}}_T)$. □

### EC.5.3. Proof of Corollary 2

We first compute

$$Q_t(x_{t-1},\widehat{\boldsymbol{b}}_t) := \sup_{\|b_t-\widehat{\boldsymbol{b}}_t\|\le\vartheta} \inf_{x_t\in\mathcal{X}_t(x_{t-1},b_t)} \left\{c_t^\top x_t + \mathbb{E}_{\widehat{\mathbb{P}}_{t+1}}\left[Q_{t+1}(x_t,\widehat{\boldsymbol{b}}_{t+1})\right]\right\},$$

Observe that the objective $c_t^\top x_t + \mathbb{E}_{\widehat{\mathbb{P}}_{t+1}}\left[Q_{t+1}(x_t,\widehat{\boldsymbol{b}}_{t+1})\right]$ is convex in $x_t$ (Shapiro et al. 2021, Section 3.2.1). Since the feasible set $\mathcal{X}_t(x_{t-1},b_t)$ consists of only linear constraints, by Assumption 2, the relaxed Slater condition is satisfied. In addition, by Assumption 3 the inner objective is bounded from below. Therefore using convex programming duality, we obtain that

$$\begin{aligned}
&\inf_{x_t\in\mathcal{X}_t(x_{t-1},b_t)} \left\{c_t^\top x_t + \mathbb{E}_{\widehat{\mathbb{P}}_{t+1}}\left[Q_{t+1}(x_t,\widehat{\boldsymbol{b}}_{t+1})\right]\right\}\\
&= \max_{y_t\in\mathbb{R}^{d_t}} \left\{y_t^\top(b_t - B_t x_{t-1}) + \inf_{x_t\ge 0}\left\{(c_t - A_t^\top y_t)^\top x_t + \mathbb{E}_{\widehat{\mathbb{P}}_{t+1}}\left[Q_{t+1}(x_t,\widehat{\boldsymbol{b}}_{t+1})\right]\right\}\right\}.
\end{aligned}$$

It follows that

$$\begin{aligned}
&Q_t(x_{t-1},\widehat{\boldsymbol{b}}_t)\\
&= \sup_{\|b_t-\widehat{\boldsymbol{b}}_t\|\le\vartheta} \max_{y_t\in\mathbb{R}^{d_t}} \left\{y_t^\top(b_t - B_t x_{t-1}) + \inf_{x_t\ge 0}\left\{(c_t - A_t^\top y_t)^\top x_t + \mathbb{E}_{\widehat{\mathbb{P}}_{t+1}}\left[Q_{t+1}(x_t,\widehat{\boldsymbol{b}}_{t+1})\right]\right\}\right\}\\
&= \max_{y_t\in\mathbb{R}^{d_t}} \sup_{\|b_t-\widehat{\boldsymbol{b}}_t\|\le\vartheta} \left\{y_t^\top(b_t - B_t x_{t-1}) + \inf_{x_t\ge 0}\left\{(c_t - A_t^\top y_t)^\top x_t + \mathbb{E}_{\widehat{\mathbb{P}}_{t+1}}\left[Q_{t+1}(x_t,\widehat{\boldsymbol{b}}_{t+1})\right]\right\}\right\}\\
&= \max_{y_t\in\mathbb{R}^{d_t}} \left\{y_t^\top(\widehat{\boldsymbol{b}}_t - B_t x_{t-1}) + \vartheta\|y_t\|_* + \inf_{x_t\ge 0}\left\{(c_t - A_t^\top y_t)^\top x_t + \mathbb{E}_{\widehat{\mathbb{P}}_{t+1}}\left[Q_{t+1}(x_t,\widehat{\boldsymbol{b}}_{t+1})\right]\right\}\right\}\\
&= \max_{y_t\in\mathbb{R}^{d_t}} \left\{(\widehat{\boldsymbol{b}}_t - B_t x_{t-1})^\top y_t + \vartheta\|y_t\|_* - \psi_t^*(A_t^\top y_t)\right\}.
\end{aligned}$$

When $\|\cdot\| = \|\cdot\|_1$, we have

$$\vartheta\|y_t\|_* = \vartheta\|y_t\|_\infty = \max_{j\in[d_t],\delta\in\{1,-1\}} \vartheta\delta e_j^\top y_t.$$

Therefore we have

$$\begin{aligned}
&\max_{y_t\in\mathbb{R}^{d_t}} \left\{(\widehat{\boldsymbol{b}}_t - B_t x_{t-1})^\top y_t + \vartheta\|y_t\|_* - \psi_t^*(A_t^\top y_t)\right\}\\
&= \max_{y_t\in\mathbb{R}^{d_t}} \left\{(\widehat{\boldsymbol{b}}_t - B_t x_{t-1})^\top y_t + \max_{j\in[d_t],\delta\in\{1,-1\}} \vartheta\delta e_j^\top y_t - \psi_t^*(A_t^\top y_t)\right\}\\
&= \max_{j\in[d_t],\delta\in\{1,-1\}} \max_{y_t\in\mathbb{R}^{d_t}} \left\{(\widehat{\boldsymbol{b}}_t - B_t x_{t-1} + \vartheta\delta e_j)^\top y_t - \psi_t^*(A_t^\top y_t)\right\}\\
&= \max_{j\in[d_t],\delta\in\{1,-1\}} \inf_{x_t\ge 0} \left\{c_t^\top x_t + \mathbb{E}_{\widehat{\mathbb{P}}_{t+1}}\left[Q_{t+1}(x_t,\widehat{\boldsymbol{b}}_{t+1})\right] : A_t x_t = \widehat{\boldsymbol{b}}_t - B_t x_{t-1} + \vartheta\delta e_j\right\},
\end{aligned}$$

where the last equality applies convex programming duality as the previous usage.

The case of uncertainty in $\boldsymbol{B}_t$ can be dealt with in a similar way. Using convex programming duality we obtain that

$$\begin{aligned}
&Q_t(x_{t-1}, \widehat{\boldsymbol{B}}_t) \\
&= \sup_{\|B_t - \widehat{\boldsymbol{B}}_t\| \le \vartheta} \inf_{x_t \in \mathcal{X}_t(x_{t-1}, B_t)} \left\{ c_t^\top x_t + \mathbb{E}_{\widehat{\mathbb{P}}_{t+1}} \left[ Q_{t+1}(x_t, \widehat{\boldsymbol{B}}_{t+1}) \right] \right\} \\
&= \sup_{\|B_t - \widehat{\boldsymbol{B}}_t\| \le \vartheta} \max_{y_t \in \mathbb{R}^{d_t}} \left\{ y_t^\top (b_t - B_t x_{t-1}) + \inf_{x_t \ge 0} \left\{ (c_t - A_t^\top y_t)^\top x_t + \mathbb{E}_{\widehat{\mathbb{P}}_{t+1}} \left[ Q_{t+1}\big(x_t, \widehat{\boldsymbol{B}}_{t+1}\big) \right] \right\} \right\} \\
&= \max_{y_t \in \mathbb{R}^{d_t}} \sup_{\|B_t - \widehat{\boldsymbol{B}}_t\| \le \vartheta} \left\{ y_t^\top (b_t - B_t x_{t-1}) + \inf_{x_t \ge 0} \left\{ (c_t - A_t^\top y_t)^\top x_t + \mathbb{E}_{\widehat{\mathbb{P}}_{t+1}} \left[ Q_{t+1}\big(x_t, \widehat{\boldsymbol{B}}_{t+1}\big) \right] \right\} \right\} \\
&= \max_{y_t \in \mathbb{R}^{d_t}} \left\{ y_t^\top (b_t - \widehat{\boldsymbol{B}}_t x_{t-1}) + \vartheta \|x_{t-1} y_t^\top\|_* + \inf_{x_t \ge 0} \left\{ (c_t - A_t^\top y_t)^\top x_t + \mathbb{E}_{\widehat{\mathbb{P}}_{t+1}} \left[ \widehat{Q}_{t+1}\big(x_t, \widehat{\boldsymbol{b}}_{t+1}\big) \right] \right\} \right\} \\
&= \max_{y_t \in \mathbb{R}^{d_t}} \left\{ (b_t - \widehat{\boldsymbol{B}}_t x_{t-1})^\top y_t + \vartheta \|x_{t-1} y_t^\top\|_* - \psi_t^*(A_t^\top y_t) \right\}.
\end{aligned}$$

When $\|\cdot\| = \|\cdot\|_{\text{op}}$, we claim that

$$\max_{\|\Delta\|_{\text{op}} \le \vartheta} y_t^\top \Delta x_{t-1} = \vartheta \|y_t\|_\infty \|x_{t-1}\| = \max_{j \in [m_t], \delta \in \{1,-1\}} \vartheta \delta y_t^\top e_j \|x_{t-1}\|.$$

Indeed, using Hölder's inequality and the definition of the operator norm, it holds that

$$y^\top \Delta x \le \|y\|_\infty \|\Delta x\|_1 \le \|y\|_\infty \|\Delta\|_{\text{op}} \|x\|.$$

Moreover, the inequality holds as equality at $\tilde{\Delta} = \vartheta \tilde{y} \tilde{x}^\top$, where $\tilde{x}$, $\tilde{y}$ are such that $\tilde{x}^\top x = \|x\|$, $\|\tilde{x}\|_* = 1$ and $\tilde{y}^\top y = \|y\|_\infty$, $\|\tilde{y}\|_1 = 1$:

$$\|\tilde{\Delta}\|_{\text{op}} = \sup_{\|v\| \le 1} \|\vartheta \tilde{y} \tilde{x}^\top v\|_1 = \vartheta \sup_{\|v\| \le 1} \|\tilde{y}\|_1 |\tilde{x}^\top v| = \vartheta, \quad y^\top \tilde{\Delta} x = \vartheta y^\top \tilde{y} \tilde{x}^\top x = \vartheta \|y\|_\infty \|x\|.$$

Therefore we have

$$\begin{aligned}
&\max_{y_t \in \mathbb{R}^{d_t}} \left\{ (\widehat{\boldsymbol{b}}_t - B_t x_{t-1})^\top y_t + \vartheta \|y_t\|_* - \psi_t^*(A_t^\top y_t) \right\} \\
&= \max_{y_t \in \mathbb{R}^{d_t}} \left\{ (\widehat{\boldsymbol{b}}_t - B_t x_{t-1})^\top y_t + \max_{j \in [m_t], \delta \in \{1,-1\}} \vartheta \delta y_t^\top e_j \|x_{t-1}\| - \psi_t^*(A_t^\top y_t) \right\} \\
&= \max_{j \in [d_t], \delta \in \{1,-1\}} \max_{y_t \in \mathbb{R}^{d_t}} \left\{ (\widehat{\boldsymbol{b}}_t - B_t x_{t-1} + \vartheta \delta e_j \|x_{t-1}\|)^\top y_t y_t - \psi_t^*(A_t^\top y_t) \right\} \\
&= \max_{j \in [d_t], \delta \in \{1,-1\}} \inf_{x_t \ge 0} \left\{ c_t^\top x_t + \mathbb{E}_{\widehat{\mathbb{P}}_{t+1}} \left[ Q_{t+1}(x_t, \widehat{\boldsymbol{b}}_{t+1}) \right] : A_t x_t = \widehat{\boldsymbol{b}}_t - B_t x_{t-1} + \vartheta \delta e_j \|x_{t-1}\| \right\}.
\end{aligned}$$

□

### EC.5.4. Proof of Corollary 3

We only show the proof for the case of uncertainty in $\boldsymbol{b}_t$. The proof for the uncertainty in $\boldsymbol{B}_t$ is similar.

Define a function $\ell_t : \Xi_t \to \mathbb{R}$ as

$$\ell_t(b_t) := \max_{y_t \in \mathbb{R}^{d_t}} \left\{ y_t^\top (b_t - B_t x_{t-1}) - \psi_t^*(A_t^\top y_t) \right\}.$$

We prove by induction that

$$Q_t(x_{t-1}, \widehat{\boldsymbol{b}}_t) = \inf_{x_t \in \mathcal{X}_t(x_{t-1}, b_t)} \left\{ c_t^\top x_t + \mathbb{E}_{\widehat{\mathbb{P}}_{t+1}} \left[ \widehat{Q}_{t+1}(x_t, \widehat{\boldsymbol{b}}_{t+1}) \right] \right\} + \infty \cdot \mathbf{1}\left\{ \lambda < \max_{s=t,\dots,T} \|\ell_t\|_{\text{Lip}} \right\}.$$

Note that the objective $c_t^\top x_t + \mathbb{E}_{\widehat{\mathbb{P}}_{t+1}}\big[\widehat{Q}_{t+1}(x_t, \widehat{\boldsymbol{c}}_{t+1})\big]$ is convex in $x_t$ (Shapiro et al. 2021, Section 3.2.1). Since the feasible set $\mathcal{X}_t(x_{t-1}, b_t)$ consists of only linear constraints, by Assumption 2, the relaxed Slater condition is satisfied. In addition, by Assumption 3 the inner objective is bounded from below. Therefore using convex programming duality, we obtain that

$$\begin{aligned}
&\inf_{x_t \in \mathcal{X}_t(x_{t-1}, b_t)} \left\{ c_t^\top x_t + \mathbb{E}_{\widehat{\mathbb{P}}_{t+1}}\left[\widehat{Q}_{t+1}\big(x_t, \widehat{\boldsymbol{b}}_{t+1}\big)\right]\right\} \\
&= \max_{y_t \in \mathbb{R}^{d_t}} \left\{ y_t^\top (b_t - B_t x_{t-1}) + \inf_{x_t \geq 0} \left\{ (c_t - A_t^\top y_t)^\top x_t + \mathbb{E}_{\widehat{\mathbb{P}}_{t+1}}\left[\widehat{Q}_{t+1}\big(x_t, \widehat{\boldsymbol{b}}_{t+1}\big)\right]\right\}\right\}.
\end{aligned}$$

It follows that

$$\begin{aligned}
&Q_t(x_{t-1}, \widehat{\boldsymbol{b}}_t) \\
&= \sup_{b_t \in \Xi_t} \left\{ \inf_{x_t \in \mathcal{X}_t(x_{t-1}, b_t)} \left\{ c_t^\top x_t + \mathbb{E}_{\widehat{\mathbb{P}}_{t+1}}\left[Q_{t+1}(x_t, \widehat{\boldsymbol{b}}_{t+1})\right]\right\} - \lambda \|b_t - \widehat{\boldsymbol{b}}_t\| \right\} + \infty \cdot \mathbf{1}\left\{\lambda < \max_{s=t+1,\ldots,T} \|\ell_s\|_{\mathrm{Lip}}\right\} \\
&= \sup_{b_t \in \Xi_t} \max_{y_t \in \mathbb{R}^{d_t}} \left\{ y_t^\top (b_t - B_t x_{t-1}) + \inf_{x_t \geq 0} \left\{ (c_t - A_t^\top y_t)^\top x_t + \mathbb{E}_{\widehat{\mathbb{P}}_{t+1}}\left[\widehat{Q}_{t+1}\big(x_t, \widehat{\boldsymbol{b}}_{t+1}\big)\right]\right\} - \lambda \|b_t - \widehat{\boldsymbol{b}}_t\| \right\} \\
&\qquad + \infty \cdot \mathbf{1}\left\{\lambda < \max_{s=t+1,\ldots,T} \|\ell_s\|_{\mathrm{Lip}}\right\} \\
&= \max_{y_t \in \mathbb{R}^{d_t}} \sup_{b_t \in \Xi_t} \left\{ y_t^\top (b_t - B_t x_{t-1}) + \inf_{x_t \geq 0} \left\{ (c_t - A_t^\top y_t)^\top x_t + \mathbb{E}_{\widehat{\mathbb{P}}_{t+1}}\left[\widehat{Q}_{t+1}\big(x_t, \widehat{\boldsymbol{b}}_{t+1}\big)\right]\right\} - \lambda \|b_t - \widehat{\boldsymbol{b}}_t\| \right\} \\
&\qquad + \infty \cdot \mathbf{1}\left\{\lambda < \max_{s=t+1,\ldots,T} \|\ell_s\|_{\mathrm{Lip}}\right\} \\
&= \max_{y_t \in \mathbb{R}^{d_t}} \sup_{b_t \in \Xi_t} \left\{ y_t^\top (b_t - B_t x_{t-1}) - \psi_t^*(A_t^\top y_t) - \lambda \|b_t - \widehat{\boldsymbol{b}}_t\| \right\} + \infty \cdot \mathbf{1}\left\{\lambda < \max_{s=t+1,\ldots,T} \|\ell_s\|_{\mathrm{Lip}}\right\}.
\end{aligned}$$

Observe that the function $\ell_t$ is convex and Lipschitz. It follows that

$$\begin{aligned}
&\sup_{b_t \in \Xi_t} \left\{ y_t^\top (b_t - B_t x_{t-1}) - \psi_t^*(A_t^\top y_t) - \lambda \|b_t - \widehat{\boldsymbol{b}}_t\| \right\} \\
&= \begin{cases} \sup_{b_t \in \Xi_t} \left\{ y_t^\top (b_t - B_t x_{t-1}) - \psi_t^*(A_t^\top y_t) - \lambda \|b_t - \widehat{\boldsymbol{b}}_t\| \right\}, & \text{if } \lambda \geq \|\ell_t\|_{\mathrm{Lip}}, \\ \infty, & \text{otherwise.} \end{cases}
\end{aligned}$$

Hence

$$\begin{aligned}
Q_t(x_{t-1}, \widehat{\boldsymbol{b}}_t) &= \max_{y_t \in \mathbb{R}^{d_t}} \left\{ (\widehat{\boldsymbol{b}}_t - B_t x_{t-1})^\top y_t - \psi_t^*(A_t^\top y_t) \right\} + \infty \cdot \mathbf{1}\left\{\lambda < \max_{s=t,\ldots,T} \|\ell_s\|_{\mathrm{Lip}}\right\} \\
&= \inf_{x_t \in \mathcal{X}_t(x_{t-1}, \widehat{\boldsymbol{b}}_t)} \left\{ c_t^\top x_t + \mathbb{E}_{\widehat{\mathbb{P}}_{t+1}}\left[\widehat{Q}_{t+1}(x_t, \widehat{\boldsymbol{b}}_{t+1})\right]\right\} + \infty \cdot \mathbf{1}\left\{\lambda < \max_{s=t,\ldots,T} \|\ell_s\|_{\mathrm{Lip}}\right\}.
\end{aligned}$$

Furthermore, we have the dual optimizer $\lambda^\star = \max_{s=2,\ldots,T} \|\ell_s\|_{\mathrm{Lip}}$.

It remains to prove

$$\|\ell_t\|_{\mathrm{Lip}} = \sup_{y \in \mathrm{dom}\, \psi_t^*(A_t^\top \cdot)} \|y\|_*, \quad \mathrm{dom}\, \psi_t^*(A_t^\top \cdot) = \mathcal{S}_t. \tag{EC.13}$$

When $t = T$, we have

$$\psi_T^*(A_T^\top y_T) = \max_{x \geq 0} \left\{ y_T^\top A_T x - c_T^\top x \right\},$$

which is zero when $A_T^\top y_T \le c_T$ and infinite otherwise. Thus $\operatorname{dom} \psi_T^*(A_T^\top \cdot) = \{y \in \mathbb{R}^{d_T} : A_T^\top y \le c_T\} = \mathcal{S}_T$. Now suppose (EC.13) holds for some $t+1$, let us prove the case of $t$, where $t \in [T-1]$. By definition of $\psi_t^*$, we have

$$
\begin{aligned}
&\psi_t^*(A_t^\top y_t) \\
&= \max_{x_t \ge 0} \left\{ y_t^\top A_t x_t - c_t^\top x_t - \mathbb{E}_{\widehat{\mathbb{P}}_{t+1}} \left[ \max_{y_{t+1} \in \mathbb{R}^{d_{t+1}}} \left\{ (\widehat{\boldsymbol{b}}_{t+1} - B_{t+1} x_t)^\top y_{t+1} - \psi_{t+1}^*(A_{t+1}^\top y_{t+1}) \right\} \right] \right\} \\
&= \max_{x_t \ge 0} \left\{ y_t^\top A_t x_t - c_t^\top x_t - \max_{\boldsymbol{y}_{t+1} \in \mathcal{Y}_{t+1}} \mathbb{E}_{\widehat{\mathbb{P}}_{t+1}} \left[ (\widehat{\boldsymbol{b}}_{t+1} - B_{t+1} x_t)^\top \boldsymbol{y}_{t+1}(\widehat{\boldsymbol{b}}_{t+1}) - \psi_{t+1}^*(A_{t+1}^\top \boldsymbol{y}_{t+1}(\widehat{\boldsymbol{b}}_{t+1})) \right] \right\} \\
&= \max_{x_t \ge 0} \min_{\boldsymbol{y}_{t+1} \in \mathcal{Y}_{t+1}} \left\{ y_t^\top A_t x_t - c_t^\top x_t - \mathbb{E}_{\widehat{\mathbb{P}}_{t+1}} \left[ (\widehat{\boldsymbol{b}}_{t+1} - B_{t+1} x_t)^\top \boldsymbol{y}_{t+1}(\widehat{\boldsymbol{b}}_{t+1}) - \psi_{t+1}^*(A_{t+1}^\top \boldsymbol{y}_{t+1}(\widehat{\boldsymbol{b}}_{t+1})) \right] \right\}.
\end{aligned}
$$

Observe that

$$
\begin{aligned}
&\max_{x_t \ge 0} \left\{ y_t^\top A_t x_t - c_t^\top x_t - \mathbb{E}_{\widehat{\mathbb{P}}_{t+1}} \left[ (\widehat{\boldsymbol{b}}_{t+1} - B_{t+1} x_t)^\top \boldsymbol{y}_{t+1}(\widehat{\boldsymbol{b}}_{t+1}) \right] \right\} \\
&= -\mathbb{E}_{\widehat{\mathbb{P}}_{t+1}} \left[ \widehat{\boldsymbol{b}}_{t+1}^\top \boldsymbol{y}_{t+1}(\widehat{\boldsymbol{b}}_{t+1}) \right] + \max_{x_t \ge 0} \left\{ (A_t^\top y_t - c_t)^\top x_t - \mathbb{E}_{\widehat{\mathbb{P}}_{t+1}} \left[ -x_t^\top B_{t+1}^\top \boldsymbol{y}_{t+1}(\widehat{\boldsymbol{b}}_{t+1}) \right] \right\} \\
&= \begin{cases} -\mathbb{E}_{\widehat{\mathbb{P}}_{t+1}} \left[ \widehat{\boldsymbol{b}}_{t+1}^\top \boldsymbol{y}_{t+1}(\widehat{\boldsymbol{b}}_{t+1}) \right], & \text{if } A_t^\top y_t + B_{t+1}^\top \mathbb{E}_{\widehat{\mathbb{P}}_{t+1}}[\boldsymbol{y}_{t+1}(\widehat{\boldsymbol{b}}_{t+1})] \le c_t, \\ \infty, & \text{otherwise.} \end{cases}
\end{aligned}
$$

Hence, $\psi_t^*(A_t^\top y_t)$ is finite if and only if there exists $\boldsymbol{y}_{t+1} \in \mathcal{Y}_{t+1}$ such that $A_t^\top y_t + B_{t+1}^\top \mathbb{E}_{\widehat{\mathbb{P}}_{t+1}}[\boldsymbol{y}_{t+1}(\widehat{\boldsymbol{b}}_{t+1})] \le c_t$. Thus,

$$
\operatorname{dom} \psi_t^*(A_t^\top y_t) = \left\{ y_t \in \mathbb{R}^{d_t} : \exists \boldsymbol{y}_{t+1} \in \mathcal{Y}_{t+1} \text{ s.t. } A_t^\top y_t + B_{t+1}^\top \mathbb{E}_{\widehat{\mathbb{P}}_{t+1}}[\boldsymbol{y}_{t+1}(\widehat{\boldsymbol{b}}_{[t+1]})] \le c_t \right\},
$$

which completes the induction for (EC.13). □

## EC.6. Additional Details for Section 5

### EC.6.1. Proof of Corollary 4

Introducing variables $\boldsymbol{x}_T^\pm, s_1, s_2 \ge 0$ we rewrite the problem as

$$
\begin{aligned}
\min_{\boldsymbol{x}_1,\dots,\boldsymbol{x}_{T-1} \ge 0, \boldsymbol{x}_T^\pm, W_T} \quad & \max_{\mathbb{P} \in \mathfrak{M}} \mathbb{E}[\boldsymbol{x}_T^+ - \boldsymbol{x}_T^-] \\
\text{s.t.} \quad & \mathbf{1}^\top \boldsymbol{x}_1 = W_1, \\
& \mathbf{1}^\top \boldsymbol{x}_{t-1} = \boldsymbol{\xi}_{t-1}^\top \boldsymbol{x}_{t-2}, \quad t = 2, \dots, T, \\
& W_T = \boldsymbol{\xi}_T^\top \boldsymbol{x}_{T-1}, \\
& \boldsymbol{x}_T^+ - \boldsymbol{x}_T^- - s_1 = -a_0 - r_0 \boldsymbol{\xi}_T^\top x_{T-1}, \\
& \boldsymbol{x}_T^+ - \boldsymbol{x}_T^- - s_2 = -a_1 - r_1 \boldsymbol{\xi}_T^\top x_{T-1}.
\end{aligned}
$$

Then the result follows from substituting the following parameter values in Corollary 2:

$$
\begin{aligned}
&A_1 = \mathbf{1}^\top,\ B_1 = 0,\ b_1 = W_1,\ c_1 = \mathbf{0}, \\
&A_t = \mathbf{1}^\top,\ B_t = -\boldsymbol{\xi}_t^\top,\ b_t = 0,\ c_t = \mathbf{0}, \\
&A_T = \begin{pmatrix} 1 & -1 & -1 & 0 \\ 1 & -1 & 0 & -1 \end{pmatrix}, B_T = (r_0 \boldsymbol{\xi}_T^\top, r_1 \boldsymbol{\xi}_T^\top)^\top,\ b_T = (-a_0, -a_1)^\top,\ c_T = (1, -1, 0, 0)^\top.
\end{aligned}
$$

□

### EC.6.2. Data

The estimated mean vector and covariance matrix are as the following.

$$\begin{array}{rccccc} & \text{EEM} & \text{TLT} & \text{SCHP} & \text{XES} & \text{SKF} \\ \widehat{\mu} = & [0.005142 & 0.006166 & 0.004729 & -0.005437 & -0.026391] \end{array}$$

$$\widehat{\Sigma} = \begin{array}{c} \begin{array}{ccccc} \text{EEM} & \text{TLT} & \text{SCHP} & \text{XES} & \text{SKF} \end{array} \\ \begin{bmatrix} 0.002862 & -0.000634 & 0.000211 & 0.007217 & -0.003882 \\ -0.000634 & 0.001495 & 0.000231 & -0.002881 & 0.001658 \\ 0.000211 & 0.000231 & 0.000111 & 0.000426 & -0.000244 \\ 0.007217 & -0.002881 & 0.000426 & 0.030046 & -0.013416 \\ -0.003882 & 0.001658 & -0.000244 & -0.013416 & 0.010840 \end{bmatrix} \end{array} \begin{array}{c} \\ \text{EEM} \\ \text{TLT} \\ \text{SCHP} \\ \text{XES} \\ \text{SKF} \end{array}$$

### EC.6.3. Algorithms

**Algorithm 1** SDDP Algorithm for Robust Reformulation with Uncertainty in $\boldsymbol{B}_t$

**Initialize:** $\{\mathfrak{Q}_t\}_{t\in[T]}$ (initial lower approximation), $\bar{x}_0^k = 0$
**while** not converge **do**
  Sample $K$ scenario paths $\{\widehat{\boldsymbol{B}}^k_{[T]}\}_{k=1}^K$
  **for** $t = 1, \dots, T-1$ **do** ▹ *Forward pass*
    **for** $k = 1, \dots, K$ **do**
      **for** $(\delta, j) \in \{1, -1\} \times [m_t]$ **do**
        $(x_{tk}^{\delta j}, v_{tk}^{\delta j}) \leftarrow \inf_{x_t \geq 0} \left\{ c_t x_t + \mathfrak{Q}_{t+1}(x_t) : \; A_t x_t + \widehat{\boldsymbol{B}}_t^k \bar{x}_{t-1}^k = b_t + \vartheta\delta \|\bar{x}_{t-1}^k\| e_j \right\}$
        ▹ *Optimal solution and optimal value*
      **end for**
      $(\delta^*, j^*) \leftarrow \arg\max_{(\delta,j)} v_{tk}^{\delta j}$;   $\bar{x}_t^k \leftarrow x_{tk}^{\delta^* j^*}$
    **end for**
  **end for**
  **for** $t = T, \dots, 2$ **do** ▹ *Backward pass*
    **for** $k = 1, \dots, K$ **do**
      **for** $i = 1, \dots, \widehat{N}_t$ **do**
        **for** $(\delta, j) \in \{1, -1\} \times [m_t]$ **do**
          $(v_{ti}^{\delta j}(\bar{x}_{t-1}^k); y_{tki}^{\delta j}) \leftarrow \inf_{x_t \geq 0} \left\{ c_t x_t + \mathfrak{Q}_{t+1}(x_t) : \; A_t x_t + \widehat{\boldsymbol{B}}_t^i \bar{x}_{t-1}^k = b_t + \vartheta\delta \|\bar{x}_{t-1}^k\| e_j \right\}$
          ▹ *Optimal value and optimal dual solution*
        **end for**
        $(\delta^\star, j^\star) \leftarrow \arg\max_{(\delta,j)} v_{ti}^{\delta j}$;   $\tilde{Q}_{ti}(\bar{x}_{t-1}^k) \leftarrow v_{ti}^{\delta^\star j^\star}$;   $\pi_{ti}^k \leftarrow y_{tik}^{\delta^\star j^\star}$
      **end for**
      $\tilde{\mathfrak{Q}}_t(\bar{x}_{t-1}^k) \leftarrow \frac{1}{N_t} \sum_{i=1}^{N_t} \tilde{Q}_{ti}(\bar{x}_{t-1}^k)$;   $\tilde{g}_t^k \leftarrow \frac{1}{N_t} \sum_{i=1}^{N_t} \vartheta \delta^\star {\pi_{ti}^k}^\top e_{j^\star} \nabla \|\bar{x}_{t-1}^k\| - (\widehat{\boldsymbol{B}}_t^i)^\top \pi_{ti}^k$
      $\mathfrak{Q}_t(\cdot) \leftarrow \max\left( \mathfrak{Q}_t(\cdot), \; \tilde{\mathfrak{Q}}_t(\bar{x}_{t-1}^k) + {\tilde{g}_t^k}^\top (\cdot - \bar{x}_{t-1}^k) \right)$
    **end for**
  **end for**
**end while**

**Algorithm 2** Out-of-sample Test

**Input:** A training tree $(\widehat{\boldsymbol{\xi}}_1, \{\widehat{\boldsymbol{\xi}}_2^{i_2}\}_{i_2=1}^{\widehat{N}_2}, \ldots, \{\widehat{\boldsymbol{\xi}}_1^{i_T}\}_{i_T=1}^{\widehat{N}_T})$, a testing tree $(\boldsymbol{\xi}_1, \{\boldsymbol{\xi}_2^{i_2}\}_{i_2=1}^{N_1}, \ldots, \{\boldsymbol{\xi}_1^{i_T}\}_{i_T=1}^{N_T})$, number of testing paths $M$

**Output:** Average out-of-sample value $V_{\text{Avg}}$

```
Solve for the optimal first-stage decision x_1 using the training tree
for m = 1 : M do
    V_m ← 0
    Sample a path (ξ_1, ξ_2^{i_2}, ..., ξ_T^{i_T}) from the testing tree
    for t = 2 : T do
        Solve for the optimal decision x_t using x_{t-1} and the training sub-tree ({ξ̂_t^{i_t}}_{i_t=1}^{N̂_t}, ..., {ξ̂_1^{i_T}}_{i_T=1}^{N̂_T})
        Observe a testing scenario ξ_t^{i_t}
        Evaluate the per-stage out-of-sample cost C_t at stage t using x_t and ξ_t^{i_t}
        V_m ← V_m + C_t
    end for
end for
V_Avg ← (1/M) Σ_{m=1}^{M} V_m
```

Where:
- Step 1: $\boldsymbol{x}_1$
- Step 3: $V_m \leftarrow 0$
- Step 4: $(\boldsymbol{\xi}_1, \boldsymbol{\xi}_2^{i_2}, \ldots, \boldsymbol{\xi}_T^{i_T})$
- Step 6: $(\{\widehat{\boldsymbol{\xi}}_t^{i_t}\}_{i_t=1}^{\widehat{N}_t}, \ldots, \{\widehat{\boldsymbol{\xi}}_1^{i_T}\}_{i_T=1}^{\widehat{N}_T})$
- Step 12: $V_{\text{Avg}} \leftarrow \frac{1}{M}\sum_{m=1}^{M} V_m$